\documentclass[11pt, a4paper]{article}

\RequirePackage{amsthm,amsmath,amsfonts,amssymb}
\RequirePackage[authoryear]{natbib} 
\RequirePackage[colorlinks,citecolor=blue,urlcolor=blue,breaklinks]{hyperref}

\usepackage{bm}
\usepackage{mathrsfs}  
\usepackage{appendix}
\usepackage{commath}
\usepackage{dsfont}
\usepackage{subfig}
\usepackage{xcolor}
\usepackage{mathtools}
\usepackage{enumitem}
\usepackage{fancyhdr}
\usepackage[a4paper, margin=1in, headheight=14pt]{geometry}
\usepackage{caption}
\usepackage{xr}
\usepackage{booktabs}
\usepackage{multirow}
\usepackage{float}
\usepackage{graphicx}
\usepackage{pdflscape}

\usepackage{setspace}
\newcommand{\N}{\mathbb{N}}
\newcommand{\R}{\mathbb{R}}
\newcommand{\eps}{\varepsilon}
\renewcommand{\P}{\mathbb{P}}
\newcommand\E{\mathbb{E}}
\newcommand{\mc}{\mathcal}

\theoremstyle{plain}
\newtheorem{theorem}{Theorem}[section]
\newtheorem{lemma}[theorem]{Lemma}
\newtheorem{corollary}[theorem]{Corollary}

\theoremstyle{remark}
\newtheorem{remark}{Remark}[section]

\newtheorem{assumption}{Assumption}[section]

\begin{document}

\title{Robustness for free: asymptotic size and power of max-tests in high dimensions\footnote{We thank Moritz Jirak for helpful comments and suggestions.}}

\author{
\begin{tabular}{c}
Anders Bredahl Kock\footnote{Kock's research was supported by the European Research Council (ERC) grant number 101124535 -- HIDI (UKRI EP/Z002222/1).  He is also a member of, and grateful for support from, i) the Aarhus Center for Econometrics (ACE), funded by the Danish National Research Foundation grant number DNRF186,  and ii) the Center for Research in Energy: Economics and Markets (CoRE).} \\ 
\footnotesize	University of Oxford \\
\footnotesize Department of Economics\\
\footnotesize	10 Manor Rd, Oxford OX1 3UQ
\\
\footnotesize	{\footnotesize	\href{mailto:anders.kock@economics.ox.ac.uk}{anders.kock@economics.ox.ac.uk}} 
\end{tabular}
\begin{tabular}{c}
David Preinerstorfer \\ 
{\footnotesize WU Vienna University of Economics and Business} \\
{\footnotesize Institute for Statistics and Mathematics} \\
{\footnotesize Welthandelsplatz 1, 1020 Vienna} \\ 
{\footnotesize	 \href{mailto:david.preinerstorfer@wu.ac.at}{david.preinerstorfer@wu.ac.at}}
\end{tabular}
}

\date{First version: October, 2025 \\
This version: October 2026}

\maketitle	

\begin{abstract}
Allowing for adversarial contamination and heavy tails, we study testing whether the mean of a high-dimensional random vector equals zero. Because standard max-tests based on sample averages are highly non-robust, we propose a max-test based on quantile-winsorized observations. The test controls asymptotic size under adversarial contamination and only requires~$m>2$ moments, while allowing dimension to grow exponentially with sample size. We fully characterize its asymptotic power function. Comparing with the standard max-test, for which we also derive a power characterization as a benchmark, we show that robustness is obtained for free at first order: under the stronger conditions that make the standard max-test valid, our robust test has identical asymptotic power. We show that the decay condition on the contamination-rate imposed is necessary for this no-first-order-loss property. Finally, we also study the role of bootstrap critical values, showing that their use never decreases power, can strictly improve asymptotic power in extremely correlated designs, but often has no first-order asymptotic effect.
\end{abstract}

\section{Introduction}

Let~$X_{1,n},\hdots,X_{n,n}$ be i.i.d.~random vectors with mean~$\mu_n \in \R^d$ and consider the testing problem\footnote{It is without loss of generality to consider testing whether~$\mu_n=0_d$. Testing whether~$\mu_n=\mu_{0,n}$ for a given vector~$\mu_{0,n}$ in~$\R^d$ is equivalent to testing whether~$X_{1,n}-\mu_{0,n},\hdots,X_{n,n}-\mu_{0,n}$ has a mean of~$0_d$.}
\begin{equation}\label{eq:testingproblem}
\mathsf{H}_{0,n}:\mu_n=0_d\qquad\text{vs.}\qquad \mathsf{H}_{1,n}:\mu_n\neq 0_d.
\end{equation}
Together with its multi-sample variants, it has a long history in statistics, modern contributions typically focusing on the ``high-dimensional'' case where dimension~$d=d_n$ is large compared to sample size~$n$; cf., e.g.,~\cite{hotelling1931}, \cite{dempster1958high}, \cite{bai1996effect}, \cite{tony2014two}, \cite{wang2015high}, \cite{xu2016adaptive}, \cite{xue2020distribution}, \cite{zhang2020simple}, \cite{kock2023consistency, kock2024enhanced}, \cite{yang2024new}, \cite{jiang2024nonparametric}, \cite{qiu2025self}, and the recent overview in~\cite{huang2022overview}. Quite naturally, many of the existing tests for~\eqref{eq:testingproblem} are based on the (suitably scaled) sample mean~$S^{\dagger}_n=n^{-1/2}\sum_{i=1}^nX_{i,n}$, say. In high-dimensional settings, which is the main focus of the present article, one often works with~$S_n$ defined coordinate-wise via
\begin{equation*}
S_{n,j}=\frac{S^{\dagger}_{n,j}}{s_{n,j}}=\frac{1}{\sqrt{n}s_{n,j}}\sum_{i=1}^nX_{i,n, j},~~ \text{for} ~~ s_{n,j}^2=\frac{1}{n}\sum_{i=1}^n(X_{i,n,j}-\overline{X}_{n,j})^2,\quad\overline{X}_{n,j}=\frac{1}{n}\sum_{i=1}^nX_{i,n,j}.
\end{equation*}
Due to their scale-invariance, and to ensure that their values are not completely determined by the coordinate of the~$X_{i,n}$ with the highest variance, tests for~\eqref{eq:testingproblem} based on~$S_n$ are preferable to tests based on~$S^{\dagger}_n$. To aggregate the evidence against the null hypothesis inherent in the coordinates of~$S_n$, one may use its maximum norm as a test statistic. One then rejects the null hypothesis in~\eqref{eq:testingproblem} if~$\|S_n\|_{\infty}$ exceeds a suitable critical value~$c_{n, 1-\alpha}$, say, with~$\alpha \in (0, 1)$ the significance level. We refer to the resulting test~$\varphi_n$, say, as a ``max-test.'' To obtain critical values~$c_{n, 1-\alpha}$ ensuring asymptotic size control, one can build upon the recent progress on Gaussian and bootstrap approximations to the distribution of the \emph{centered} statistic~$S^{\dagger}_n$, i.e.,
\begin{equation*}
n^{-1/2} \sum_{i = 1}^n (X_{i,n} - \mu_n) = S^{\dagger}_n - n^{1/2}\mu_n,
\end{equation*}
over hyperrectangles, often imposing only minimal assumptions on the structure of the covariance matrix~$\Sigma_n$ of the~$X_{i,n}$, cf., e.g., \cite{chernozhukov2017central}, \cite{deng2020beyond}, \cite{lopes2020bootstrapping}, \cite{kuchibhotla2020high}, \cite{das2021central}, \cite{koike2021notes}, \cite{kuchibhotla2021high}, \cite{lopes2022central}, \cite{chernozhuokov2022improved}, \cite{fang2023high}, \cite{chernozhukov2023nearly}, \cite{koike2024high}. We refer to the review in~\cite{chernozhukov2023high} for further references. Such results can be used to readily obtain Gaussian approximations to the distribution of the \emph{centered} statistic~$S_n$, i.e., to 
\begin{equation}\label{eqn:centeredSN}
\frac{1}{\sqrt{n}s_{n,j}} \sum_{i = 1}^n (X_{i,n, j} - \mu_{n,j}) = S_{n,j} - n^{1/2}\mu_{n,j}/s_{n,j},
\end{equation}
see, e.g., our Lemma~\ref{lem:GAStanfardizedSums}. Because~$S_n$ clearly coincides with its centered version under the null hypothesis~$\mu_n = 0_d$ (but not under the alternative), this then allows one to control asymptotic size of the test through a suitable choice of~$c_{n, 1-\alpha}$ based on asymptotic approximations for the centered statistic~$S_n$. Power properties (mainly consistency) of such tests against \emph{particular} alternatives, e.g., sparse ones, have been established in the literature, cf., e.g., \cite{chang2017simulation}. Beyond that, the asymptotic power \emph{function} of such tests is less understood.  As the first contribution of the present article, we completely characterize the asymptotic power function of a max-test based on~$S_n$ in Theorem~\ref{thm:sums}. This result mainly serves as a benchmark for our results on tests based on winsorized statistics, but is of some independent interest. 

Let us emphasize that in contrast to asymptotic size, the asymptotic power of tests based on~$S_n$ does not follow from asymptotic approximations for the \emph{centered} statistic~$S_n$ as in~\eqref{eqn:centeredSN}: $S_n$ does \emph{not} coincide with its centered version under the alternative~$\mu_n \neq 0_d$, but differs in each coordinate~$j$ by the data-dependent term~$n^{1/2}\mu_{n,j}/s_{n,j}$.\footnote{Clearly, if the difference was data-independent, one could simply ``absorb it into the hyperrectangle,'' but this is not the case.\label{foot:ddshift}}  Therefore, studying asymptotic power requires us to establish approximation results for the \emph{non-centered} statistic~$S_n$, which are novel technical contributions, and which we obtain in Theorems~\ref{thm:SumGA} and~\ref{thm:SumGAnonloc}. These two theorems separately treat ``local'' and ``non-local'' deviations from the null. The proof strategy for obtaining Theorem~\ref{thm:sums} is discussed in more detail in Section~\ref{sec:pf1}.

One drawback of the above testing approach is that tests based on~$S_n$ are notoriously non-robust. Consider the \emph{adversarial contamination} setting, where an adversary first inspects the sample and hands over corrupted random vectors~$\tilde{X}_{1,n},\hdots,\tilde{X}_{n,n}$ in~$\R^d$ to the statistician satisfying
\begin{equation}\label{eq:contamfrac}
\envert[1]{\cbr[1]{i\in\cbr[0]{1,\hdots,n}:\tilde{X}_{i,n}\neq X_{i,n}}}
\leq
\overline{\eta}_n n,
\end{equation}
and where~$\overline{\eta}_n\in[0,1/2)$ is a non-random upper bound on the fraction of contaminated observations. Analyzing robustness properties of statistical procedures under adversarial contamination has recently attracted considerable interest, cf., e.g.,~\cite{lai2016agnostic}, \cite{cheng2019high}, \cite{diakonikolas2019robust}, \cite{hopkins2020robust}, \cite{LM21}, \cite{minsker2021robust}, \cite{bhatt2022minimax}, \cite{depersin2022robust}, \cite{dalalyan2022all}, \cite{minasyan2023statistically}, \cite{minsker2023efficient}, \cite{oliveira2025finite}, and \cite{Wins1}. It is easy to see that even for~$\overline{\eta}_n = 1/n$ (corresponding to only a single observation being contaminated) the test~$\varphi_n$ breaks down in this setting, cf.~Section~\ref{sec:breakdown}. Furthermore, the limiting results cited in the paragraph above~\eqref{eqn:centeredSN} essentially require the distributions of the~$X_{i,n}$ to possess \emph{exponentially decaying tails} to allow for an exponential growth of~$d$ in sample size, cf.~\cite{zhangwu2017} and~\cite{kock2024remark}. 

In the present article, we propose max-tests based on (normalized) \emph{quantile-winsorized} means $S_{n,W}$, cf.~Equation~\eqref{eq:S_nWS}. We establish their validity and study their power properties in the presence of adversarial contamination and heavy tails, i.e., only assuming the~$X_{i,n}$ to possess~$m > 2$ moments. Concerning their validity, Theorem~\ref{thm:Wins} shows that if $$\log(d)/{n}^{\frac{m-2}{5m-2}}\to 0 \quad \text{ and }\ \quad \sqrt{n\log(d)}\overline{\eta}_{n}^{1-\frac{1}{m}}\to 0,$$ then the test
\begin{equation*}
\varphi_{n,W}
=
\mathds{1}\cbr[1]{\enVert[0]{S_{n,W}}_\infty>c_{n,1-\alpha}},\qquad n\in\N,
\end{equation*}
has asymptotic size at most~$\alpha$. It would be natural to think that achieving this robustness to contamination and heavy tails through winsorization must entail some loss in power. This is, however, not the case: Theorem~\ref{thm:Wins} also shows that~$\varphi_{n,W}$ obeys an asymptotic characterization of its power function that is (first-order) \emph{identical} to that of~$\varphi_n$, which we obtain mainly as a benchmark in Theorem~\ref{thm:sums} (in a setting necessarily without contamination, and worse dimension dependence). That is, the robustness of~$\varphi_{n,W}$ is for free in terms of first-order asymptotics: whenever our size and power guarantees for~$\varphi_n$ are valid, i.e., when~$d/n^{m/2-1-\xi}\to 0$ for some~$\xi>0$ and~$\overline{\eta}_n=0$, then~$\varphi_{n,W}$ has exactly the same guarantees. Crucially, for any~$m > 2$, the guarantees for~$\varphi_{n,W}$ extend even to the case of 
\begin{itemize}
\item $\log(d)=o\del[1]{n^{(m-2)/(5m-2)}}$, i.e.~$d$ growing exponentially in~$n$, and
\item some adversarial contamination: $\overline{\eta}_n>0$ but~$\sqrt{n\log(d)}\overline{\eta}_n^{1-\frac1m}\to 0$.
\end{itemize}

As the main technical challenge to prove Theorem~\ref{thm:Wins}, we need to establish high-dimensional Gaussian approximations to the distribution of \emph{non-centered} normalized winsorized means, cf.~Theorems~\ref{thm:SumGAWins} and \ref{thm:SumGAWinsLarge}. These results extend Gaussian approximations for the distribution of \emph{centered} quantile-winsorized means in~\cite{kock2025high} in the same way as our approximation results for the non-centered statistic~$S_n$, which we discussed above, extend their centered versions. The limiting results in~\cite{kock2025high} would only deliver asymptotic size control of max-tests based on (normalized) winsorized means, but they are silent about asymptotic power, because the centered and non-centered statistics do not coincide under the alternative, due to a data-dependent shift (cf.~also Footnote~\ref{foot:ddshift}). 

Gaussian approximations for the distribution of statistics based on centered (non-normalized) trimmed means have been obtained in~\cite{resende2024robust}; practical disadvantages of the trimming-scheme used in the latter article have been discussed in~\cite{kock2025high}. \emph{Absent contamination}, asymptotic size control of tests based on certain winsorized means has recently been established in \cite{gao2021refinedcramertypemoderate} and \cite{liu2024robust}. Their winsorization methods are fundamentally different from ours, which is carefully constructed to ensure a minimal amount of winsorization (which is important to avoid power loss) while still guaranteeing validity. Neither \cite{gao2021refinedcramertypemoderate} nor \cite{liu2024robust} establish power properties of tests, which is the main focus of the present article.

Although~$\varphi_{n,W}$ has the same asymptotic power as~$\varphi_n$ based on the uncontaminated data when~$\sqrt{n\log(d)}\overline{\eta}_n^{1-1/m}\to 0$, this requires~$\overline{\eta}_n\to 0$. Theorem~\ref{thm:EtaNecessaryshort} shows that this is not a defect of~$\varphi_{n,W}$ but is necessary to achieve robustness for free at first order: If~$\sqrt{n\log(d)}\overline{\eta}_n^{1-1/m}\not\to 0$ then \emph{no} test can obey the Gaussian approximation in Theorem~\ref{thm:Wins} and \emph{every} test has asymptotic power no greater than its size against an alternative against which~$\varphi_{n,W}$ has non-trivial power when~$\sqrt{n\log(d)}\overline{\eta}_n^{1-1/m}\to 0$.

We also investigate the power properties of the test obtained after replacing the critical values~$c_{n,1-\alpha}$ in the definition of~$\varphi_{n, W}$ by data-driven (bootstrap-based) ones~$c_{n,1-\alpha}^B$, say. The latter take the correlation structure of the entries of the~$X_{i,n}$ into account. In Theorem~\ref{thm:BSWins} we show that the resulting test,~$\varphi_{n,W,B}$, say, has asymptotic rejection frequency \emph{equal to}~$\alpha$ uniformly over all distributions satisfying~$\mathsf{H}_{0,n}$ --- a stronger property than the size control provided for~$\varphi_n$ and~$\varphi_{n,W}$. We also characterize the asymptotic power of~$\varphi_{n,W,B}$, which is shown to be never smaller than that of~$\varphi_{n,W}$. However, Theorem~\ref{thm:comparison} reveals that the asymptotic power of~$\varphi_{n,W}$ and~$\varphi_{n,W,B}$ is often \emph{identical}. In particular,~$\varphi_{n,W}$ and~$\varphi_{n,W,B}$ are proven to have identical asymptotic power if the off-diagonal elements of the~$d\times d$ correlation matrix~$\Sigma_{0,n}$ of the~$X_{i,n}$ decay at a logarithmic rate as the distance between their indices increases. Thus, even for rather highly correlated~$X_{i,n}$ there is often no gain in asymptotic power from using bootstrap critical values, which may be somewhat surprising. On the other hand, we also exhibit very highly correlated settings and alternatives~$\mu_n$ against which~$\varphi_{n,W}$ has asymptotic power \emph{zero}, yet~$\varphi_{n,W,B}$ has asymptotic power \emph{one}. 

A recent literature in theoretical computer science studies robust distribution testing under an $\varepsilon$-fraction of corruptions; see, among others, \cite{DiakonikolasKaneStewart2017}, \cite{GeorgeCanonne2022}, and in particular~\cite{canonne2026full}. These works formulate nonasymptotic (Euclidean-norm) separated testing problems, and characterize the sample size required for reliable discrimination as a function of dimension, signal strength, and $\varepsilon$, often under the highly structured Gaussian location model $\mathsf{N}_d(\mu,I_d)$ or under additional sparsity assumptions. While these frameworks allow $\varepsilon$ to be non-vanishing, this does not come ``for free:'' the separation from the null required for reliable testing, the required sample size, or both generally depend on $\varepsilon$. Our objective is different: over a broad class requiring only finite $m$th moments and without structural restrictions on the correlation matrix $\Sigma_{0,n}$, we characterize the first-order power function uniformly over both null and alternative sequences and give conditions under which it coincides with the corresponding uncontaminated Gaussian benchmark. We show that this stronger no-first-order-loss objective is attainable when $\sqrt{n\log(d)}\overline{\eta}_n^{1-\frac1m}\to0$, and impossible for \emph{any} sequence of tests when this condition fails.

\section{Setting and notation}\label{sec:setting}

Consider a measurable space~$(\Omega,\mc{F})$ with the two arrays $X_{1,n},\hdots,X_{n,n}$ (uncontaminated observations) and~$\tilde{X}_{1,n},\hdots,\tilde{X}_{n,n}$ (potentially contaminated observations) of random vectors~$X_{i,n}=(X_{i,n,1},\hdots,X_{i,n,d_n})': \Omega \to \R^{d_n}$ and~$\tilde{X}_{i,n}=(\tilde{X}_{i,n,1},\hdots,\tilde{X}_{i,n,d_n})': \Omega \to \R^{d_n}$, and a sequence of random vectors~$Z_{d_n}: \Omega \to \R^{d_{n}}$. Throughout we assume that~$X_{1,n},\hdots,X_{n,n}$ and~$\tilde{X}_{1,n},\hdots,\tilde{X}_{n,n}$ satisfy~\eqref{eq:contamfrac}, and we often suppress the dependence of~$d=d_n$ on~$n$ in our notation. For any probability measure~$P$ on~$\mc{F}$, we denote by~$E_P$ the corresponding expectation operator.  Furthermore, we let (grant the expressions are defined)
\begin{itemize}
\item $\mu_n(P):=(\mu_{n,1}(P),\hdots,\mu_{n,d}(P))'=E_P X_{1,n}$;
\item $\Sigma_n(P)$ denote the covariance matrix of~$X_{1,n}$ under~$P$; and
\item for~$m\in[2,\infty)$ let~$\sigma_{m,n,j}^m(P):=E_P|X_{1,n,j}-\mu_{n,j}(P)|^m$.
\end{itemize} 
All of our results hold uniformly over sub-classes of distributions on~$(\Omega,\mc{F})$ imposing certain moment bounds, which we define next: for~$b_1\in(0,\infty)$,~$b_2\in(b_1,\infty)$, and~$m \geq 2$, we denote by~$\mc{P}(b_1,b_2,m)$ the sub-class of probability measures~$P$ on~$\mc{F}$ satisfying for every~$n\in\N$ that:\footnote{If~$b_1>b_2$ then~$\mc{P}(b_1,b_2,m)$ is empty by Lyapunov's inequality.}
\begin{enumerate}
\item The~$X_{1,n},\hdots,X_{n,n}$ are i.i.d.;
\item $E_P |X_{1,n,j}|^m<\infty$ for~$j=1,\hdots,d$;
\item $b_1 \leq \min_{j=1,\hdots,d}\sigma_{2,n,j}(P)$, and~$\max_{j=1,\hdots,d}\sigma_{m,n,j}(P)\leq b_2$; and
\item $Z_{d} \sim \mathsf{N}_d(0_{d},I_{d})$ and~$Z_{d}$ is independent of~$(X_{1,n},\hdots,X_{n,n},\tilde{X}_{1,n},\hdots,\tilde{X}_{n,n})$. 
\end{enumerate}
Furthermore, we denote the subset of~$\mc{P}(b_1,b_2,m)$ satisfying~$\mathsf{H}_{0,n}$ in~\eqref{eq:testingproblem} for every~$n\in\N$ by
\begin{equation*}
\mc{P}^0(b_1,b_2,m)=\cbr[1]{P\in \mc{P}(b_1,b_2,m): \mu_n(P)=0_d\text{ for every }n\in\N}.
\end{equation*}
We assume throughout that the setting, i.e.,~$(\Omega, \mathcal{F})$ and the above-defined arrays, are rich enough to ensure that~$ \mc{P}^0(b_1,b_2,m) \neq \emptyset$. 

We shall sometimes abbreviate~$Z_d$ by~$Z$ if the dimension is clear from the context. For~$\alpha\in(0,1)$, we let~$c_{n,1-\alpha}$ be the~$(1-\alpha)$-quantile of~$\enVert[0]{Z_{d}}_\infty$, that is~$P\del[1]{\enVert[0]{Z_{d}}_\infty>c_{n,1-\alpha}}= \alpha$, which in case~$d \to \infty$ implies (cf., e.g., Equations~(7) and~(8) in~\cite{kock2023consistency})
\begin{equation}\label{eq:CVDiag}
c_{n,1-\alpha}
=
\sqrt{2\log(d)}-\frac{\log\log(d)+\log(4\pi)}{2\sqrt{2\log(d)}}-\frac{\log\del[0]{-\log(1-\alpha)/2}+o(1)}{\sqrt{2\log(d)}}.
\end{equation}

\section{Tests based on normalized arithmetic means}\label{sec:arim}
To put the results on tests based on winsorized means in Section~\ref{sec:winsm} into perspective, we first state asymptotic size and power properties of the sequence of tests 
\begin{equation}\label{eqn:aritest}
\varphi_n=\mathds{1}\cbr[1]{\enVert[0]{S_n}_\infty>c_{n,1-\alpha}},\quad n\in\N.
\end{equation}
Recall from the introduction that 
\begin{equation}\label{eqn:Sndef}
S_{n,j}=\frac{1}{\sqrt{n}s_{n,j}}\sum_{i=1}^nX_{i,n, j},~~ \text{for} ~~ s_{n,j}^2=\frac{1}{n}\sum_{i=1}^n(X_{i,n,j}-\overline{X}_{n,j})^2,\quad\overline{X}_{n,j}=\frac{1}{n}\sum_{i=1}^nX_{i,n,j},
\end{equation}	
with~$S_{n,j}$ defined as~$0$ in case~$s_{n,j} = 0$ (all square roots are positive square roots if not mentioned otherwise). We stress that the critical values~$c_{n,1-\alpha}$ in~\eqref{eqn:aritest} are non-stochastic and  \emph{neither} depend on the covariance matrix~$\Sigma_n(P)$ nor on other unknown population quantities (such as~$b_1$ or~$b_2$). The tests in~\eqref{eqn:aritest} and the analysis provided in this section should be considered as a benchmark for similar results on more elaborate procedures discussed later on (e.g., tests based on winsorization and bootstrap critical values). 

Because~$S_n$ is based on~$X_{1, n}, \hdots, X_{n,n}$, the test statistic~$S_n$ is only computable in the non-contaminated setting, i.e., when $\overline{\eta}_n=0$. Hence, from an applied point of view, the results in this section all implicitly impose the condition that $\overline{\eta}_n=0$, even though the theorems do not state this assumption. One could in principle analyze properties of the related test $\tilde{\varphi}_n=\mathds{1}\cbr[1]{\enVert[0]{\tilde{S}_n}_\infty>c_{n,1-\alpha}}$ based on the test statistic~$\tilde{S}_n$ with coordinates
\begin{equation}\label{eq:notationcontamsum}
\tilde{S}_{n,j} := \frac{1}{\sqrt{n}\tilde{s}_{n,j}}\sum_{i=1}^n\tilde{X}_{i,n, j}, ~~ \text{for} ~~ \tilde{s}_{n,j}^2=\frac{1}{n}\sum_{i=1}^n(\tilde{X}_{i,n,j}-\overline{\tilde{X}}_{n,j})^2,\quad\overline{\tilde{X}}_{n,j}=\frac{1}{n}\sum_{i=1}^n\tilde{X}_{i,n,j},
\end{equation}
which is computable also in the contaminated setting. However, due to the non-robustness of the sample average, any non-trivial asymptotic result for~$\tilde{\varphi}_n$ would need to impose an uncontaminated sample at least for~$n$ large enough anyways, cf.~Theorem~\ref{thm:Breakdown} below. Hence, we formulate our first set of results in terms of~$S_n$.

Write~$D_n(P)=\text{diag}\del[1]{\sigma_{2,n,1}(P),\hdots,\sigma_{2,n,d}(P)}$ and denote the (common) correlation matrix of the uncontaminated observations~$X_{1,n},\hdots,X_{n,n}$ by (grant all variances are positive)~
\begin{equation}\label{eqn:sigma0def}
\Sigma_{0,n}(P)=D_n^{-1}(P)\Sigma_n(P) D_n^{-1}(P).
\end{equation} Furthermore, for any symmetric and positive semi-definite matrix~$A$, we denote by~$A^{1/2}$ its unique symmetric and positive semi-definite matrix square root. The following theorem recalls size control of the test in~\eqref{eqn:aritest} and fully characterizes its asymptotic power.
\begin{theorem}\label{thm:sums}
Let~$b_1\in(0,\infty)$,~$b_2\in(b_1,\infty)$, and~$m > 4$. Let~$\alpha\in(0,1)$, and suppose that~$d\geq 2$. If there exists a~$\xi\in(0,1)$ such that~$\frac{d}{n^{m/2-1-\xi}}\to 0$, then the following holds:
\begin{enumerate}[leftmargin=*]
\item Asymptotic size control:
\begin{equation}\label{eq:sizecherno}
\limsup_{n\to\infty} \sup_{P\in\mc{P}^0(b_1,b_2,m)} P\del[1]{\enVert[0]{S_n}_\infty >c_{n,1-\alpha}}\leq \alpha;
\end{equation}
and~$\lim_{n \to \infty} P_n(\enVert[0]{S_n}_\infty >c_{n,1-\alpha}) = \alpha$ for every sequence~$P_n \in \mc{P}^0(b_1,b_2,m)$ satisfying~$\Sigma_{0, n}(P_n) = I_d$ for~$n$ large enough. If~$d \to \infty$, then~$\lim_{n \to \infty} P_n(\enVert[0]{S_n}_\infty >c_{n,1-\alpha}) = \alpha$ holds for every sequence~$P_n$ in~$\mc{P}^0(b_1,b_2,m)$ contained in some~$\mathcal{P}(b_1, b_2, m, \bm{r})$ (defined in~\eqref{eqn:Pcorrbd} below) for a sequence~$\bm{r} = (r_l)$ such that~$\log(l)r_l \to 0$ as $l\to \infty.$
\item Gaussian approximation for the power function:
\begin{equation}
\begin{aligned}\label{eq:GaussApproxCherno}
\lim_{n\to\infty}\sup_{P\in\mc{P}(b_1,b_2,m)} \bigg|  &P \left( \enVert[0]{S_n}_\infty>c_{n,1-\alpha}\right) \\
& \hspace{1cm}-P\left(\enVert[0]{\Sigma_{0,n}^{1/2}(P)Z+\sqrt{n}D_n^{-1}(P)\mu_n(P)}_\infty>c_{n,1-\alpha}\right) \bigg|=0.
\end{aligned}
\end{equation}
\end{enumerate}
\end{theorem}
The first part of Theorem~\ref{thm:sums} shows that if~$d/n^{m/2-1-\xi}\to 0$ for some~$\xi>0$, then the critical values~$c_{n,1-\alpha}$ ensure that the asymptotic size of~$\varphi_n$ does not exceed~$\alpha$. If the setting is rich enough to admit a sequence~$P_n \in \mc{P}^0(b_1,b_2,m)$ satisfying (i)~$\Sigma_{0, n}(P_n) = I_d$ for~$n$ large enough, or (ii) that the degree of correlation is not too strong (cf.~\eqref{eqn:Pcorrbd} further below) \emph{and}~$d \to \infty$,  then the asymptotic size in~\eqref{eq:sizecherno} in fact equals~$\alpha$. Such sequences often exist and the asymptotic size hence often equals~$\alpha$. Together with~\eqref{eq:CVDiag} it then follows that~$c_{n,1-\alpha}$ are the smallest possible (hence power maximizing) \emph{data-independent} critical values up to terms of order~$o(1/\sqrt{\log(d)})$, ensuring that~$\varphi_n$ is of asymptotic size~$\alpha$. Tests based on data-driven (bootstrap) critical values are studied in Section~\ref{sec:BootstrapApprox} in the winsorized case, cf.~also Footnote~\ref{fn:boot}.  

Equation~\eqref{eq:GaussApproxCherno} in Theorem~\ref{thm:sums} uniformly characterizes the asymptotic power of~$\varphi_n$ as the power function of the supremum-norm based test in an appropriate Gaussian sequence model. Crucially, Theorem~\ref{thm:Wins} below will reveal that this asymptotic power is identical to that of tests based on winsorized means, but with the power results for the latter remaining valid even for (i)~$\log(d)=o\del[1]{n^{(m-2)/(5m-2)}}$, instead of the polynomial rate~$d/n^{m/2-1-\xi}\to 0$ for some~$\xi>0$ for~$\varphi_n$ and (ii) a moderate amount of adversarial contamination, instead of (implicitly) imposing~$\overline{\eta}_n=0$ for~$\varphi_n$.

\begin{remark}\label{rem:qiuchenshao}
By exploiting the self-normalized nature of~$S_n$, the asymptotic \emph{size} control in~\eqref{eq:sizecherno} continues to hold for (a test related to)~$\varphi_n$ by Proposition 3.1 in \cite{qiu2025self} even if~$\log(d)=o(n^{(m-2)/m})$ for~$m\in(2,3]$ over subsets of the null hypothesis restricting the correlation matrices~$\Sigma_{0,n}$. Specifically, that result, which is under~$\mathsf{H}_{0,n}$ only,  i) requires~$\max_{1\leq i<j\leq d}|\Sigma_{0,n,i,j}|$ to be bounded away from one, ii) restricts the sum of the absolute entries of~$\Sigma_{0,n}$, iii) requires the maximal squared row sum of~$\Sigma_{0,n}$ to be bounded, iv) requires~$\Sigma_{0,n}$ to be positive definite, and v) requires the inverse of~$\Sigma_{0,n}$ to be strictly diagonally dominant.\footnote{As a result, size control of a test based on that result can currently only be guaranteed for~$\log(d)=o(n^{(m-2)/m})$ over subsets of~$\mathsf{H}_{0,n}$ satisfying these conditions on~$\Sigma_{0,n}$.} These conditions on~$\Sigma_{0,n}$ are not imposed in Part 1.~of Theorem~\ref{thm:sums} which, however, imposes~$d/n^{m/2-1-\xi}\to 0$ for some~$\xi>0$. Although it is possible that future work will show that Gaussian approximations for~$S_n$ are valid for~$d$ growing exponentially in~$n$ under all of~$\mc{P}^0(b_1,b_2,m)$ (cf.~the \emph{conjecture} in Remark 2.4 of \cite{qiu2025self}), they will still need~$\overline{\eta}_n=0$ by Theorem~\ref{thm:Breakdown} below. We also stress that Proposition 3.1 in \cite{qiu2025self} is under~$\mathsf{H}_{0,n}$ only and thus silent about power. A full characterization of the power of tests based on~$S_n$ as in~\eqref{eq:GaussApproxCherno} (not only consistency against certain sparse alternatives) when~$d$ grows exponentially in~$n$ is likely to be challenging. A reason for this is that~$\max_{1\leq j\leq d}|s_{n,j}-\sigma_{2,n,j}(P)|$ need not be close to zero when~$m$ is small. The results for tests based on the winsorized means that we present in Section~\ref{sec:winsm} allow~$\overline{\eta}_n>0$, allow~$\log(d)=o\del[1]{n^{(m-2)/(5m-2)}}$ without further restrictions on~$\Sigma_{0,n}$, and provide a \emph{full characterization} of the power of the resulting test under~$\mathsf{H}_{1,n}$, i.e.,~are not obtained solely under~$\mathsf{H}_{0,n}$. 
\end{remark}

\subsection{Proof strategy of Theorem~\ref{thm:sums}}\label{sec:pf1}

As mentioned in the introduction, Gaussian approximations to the centered quantity~\eqref{eqn:centeredSN} can readily be obtained by extending existing Gaussian approximation results to~$S^{\dagger}_n$ (cf., e.g., Lemma~\ref{lem:GAStanfardizedSums}). Such extensions, however, only allow one to obtain asymptotic \emph{size} properties of the related tests, as~$S_n$ coincides with the centered quantity~\eqref{eqn:centeredSN} only under the null hypothesis~$\mu_n = 0$. Obtaining the asymptotic characterization of the power function given in Theorem~\ref{thm:sums} requires us to study Gaussian approximations for the non-centered quantity~$S_n$ in case~$\mu_n \neq 0$. To this end, we develop two results: we first establish Theorem~\ref{thm:SumGA}, which is informative as long as the deviation from the null is ``small'' (cf.~the first condition in~\eqref{eqn:locmu} for the formal condition). To establish this result, we combine Theorem~2.5 of~\cite{chernozhuokov2022improved} concerning centered and non-normalized statistics, the Gaussian anti-concentration-, Khatri-{\v{S}}id{\'a}k-, and Fuk-Nagaev-inequality. Then, in Theorem~\ref{thm:SumGAnonloc}, we establish a Gaussian approximation result informative in the remaining regime where the deviation from the null is not ``small.'' Combined with Theorem~\ref{thm:SumGA} this then delivers~\eqref{eq:GaussApproxCherno} in Theorem~\ref{thm:sums}. Note that~\eqref{eq:sizecherno} follows from~\eqref{eq:GaussApproxCherno} with~$\mu_n(P)=0$ and the Khatri-{\v{S}}id{\'a}k inequality, which shows that~$$P\del[1]{\enVert[0]{\Sigma_{0,n}^{1/2}(P)Z}_\infty>c_{n,1-\alpha}}\leq P\del[1]{\enVert[0]{Z}_\infty>c_{n,1-\alpha}}= \alpha.$$

\subsection{Breakdown of~$\varphi_n$ in the presence of contamination}\label{sec:breakdown}
The following result on the breakdown of~$\varphi_n$ when at least two observations can be adversarially contaminated is  trivial, but we decided to include it, as it helps to put the robustness of tests based on winsorized means in Section~\ref{sec:winsm} into perspective. We adopt the notational conventions around~\eqref{eq:notationcontamsum}.
\begin{theorem}\label{thm:Breakdown}
Fix~$n\geq 2$ and~$\alpha\in(0,1)$. Suppose the contaminated observations are obtained via~$\tilde{X}_{i,n}=X_{i,n}$ for~$i=1,\hdots,n-2$,~$\tilde{X}_{n-1,n}=(1,\hdots,1)'$, and $\tilde{X}_{n,n}=-\sum_{i=1}^{n-1}\tilde{X}_{i,n}$, such that~\eqref{eq:contamfrac} is satisfied with~$\overline{\eta}_n=2/n$. Then,~$\tilde{S}_n = 0_d$ and hence~$\tilde{\varphi}_{n}=0$, which implies, for any probability measure~$P$ on~$\mathcal{F}$, that
\begin{equation*}
P\del[1]{||\tilde{S}_n||_\infty>c_{n,1-\alpha}}=0;
\end{equation*}	
in particular the rejection frequency of~$\tilde{\varphi}_n$ does not depend on~$\mu_n(P)$, but equals zero under~$\mathsf{H}_{0,n}$ and~$\mathsf{H}_{1,n}$.
\end{theorem}
The contamination strategy in Theorem~\ref{thm:Breakdown} changes~$X_{n,n}$ to~$\tilde{X}_{n,n}=-\sum_{i=1}^{n-1}\tilde{X}_{i,n}$ to ensure that the vector of numerators~$\sum_{i=1}^n\tilde{X}_{i,n}$ of~$\tilde{S}_n$ equals zero. To avoid that any denominator~$\tilde{s}_{n,j}$ is zero (and the test statistic~$\tilde{S}_n$ would hence not be well-defined, or its conclusion would depend on how it is defined on exceptional sets), we set~$\tilde{X}_{n-1,n}=(1,\hdots,1)'$.\footnote{A denominator $\tilde{s}_{n,j}$ would equal zero without this safeguard (i.e., setting~$\tilde{X}_{i,n}=X_{i,n}$ for~$i=1,\hdots,n-1$ and~$\tilde{X}_{n,n}=-\sum_{i=1}^{n-1}\tilde{X}_{i,n}$) if and only if~$X_{1,n,j}=\hdots=X_{n-1,n,j}=0$.} Note that  for every sample size~$n\geq 2$ (and thus also asymptotically),~$\tilde{\varphi}_n$ has size and power equal to zero, irrespective of how large~$\mu_n(P)$ is. Thus,~$\tilde{\varphi}_{n}$ is completely uninformative about~$\mu_n$. The tests based on winsorized means discussed in Section~\ref{sec:winsm} do not suffer from this --- under the contamination in Theorem~\ref{thm:Breakdown} they have the same asymptotic size and power as if no contamination was present. 

The contamination strategy used in Theorem~\ref{thm:Breakdown} is convenient in that it breaks the test~$\tilde{\varphi}_n$ (i) for any~$\alpha \in (0, 1)$, and (ii) regardless of how the test is defined on that part of the sample space where the quotients in the test statistic are not well-defined. However, the contamination strategy requires \emph{two} observations to be contaminated, whereas the breakdown point of the sample average (before standardization) is~$1/n$. To see that the test~$\tilde{\varphi}_n$ can typically also be broken by only contaminating a single observation, one may proceed as follows. For every~$j=1,\hdots,d$ one can make, e.g., the first observation~$X_{1,n,j}$  so large as to ensure that~$|\tilde{S}_{n,j}|\approx \sqrt{n/(n-1)}$, which then leads to the analogous conclusion (of~$\tilde{\varphi}_n$ never rejecting) as in Theorem~\ref{thm:Breakdown} at least for all practically relevant~$\alpha$ (which satisfy~$c_{n,1-\alpha} > \sqrt{2}$). 

\section{Tests based on normalized winsorized means}\label{sec:winsm}
Although Theorem~\ref{thm:sums} establishes asymptotic size control and a complete characterization of the asymptotic power of the test~$\varphi_n$ based on~$S_n$, these guarantees require no adversarial contamination (i.e.,~$\overline{\eta}_n=0$) and~$d/n^{m/2-1-\xi}\to 0$ for some~$\xi>0$. If, for example,~$m=5$ this only allows dimension~$d$ to grow as fast as~$n^{3/2}$. We now show how replacing~$S_n$ by a vector of normalized \emph{winsorized} means can allow~$\log(d)$ to grow as fast as~$n^{\frac{m-2}{5m-2}}$ --- an exponential rate --- and simultaneously allow~$\overline{\eta}_n>0$. What is more, the resulting test~$\varphi_{n,W}$, say, obeys a Gaussian approximation to its power function that is \emph{identical} to that of~$\varphi_n$ whenever the more restrictive conditions on~$d$ and~$\overline{\eta}_n$ from Theorem~\ref{thm:sums} are satisfied. We stress that to implement~$\varphi_{n,W}$, one must know an upper bound~$\overline{\eta}_n$ on the fraction of contaminated observations.\footnote{When using~$\varphi_n$ one also implicitly takes a stance on~$\overline{\eta}_n$ as it must be zero for the test to be useful.} The true fraction need not be known.

\subsection{Winsorized means and tests based on them}
To define the test~$\varphi_{n,W}$, we first introduce the winsorized means that it is based on. Gaussian approximations to their distributions under~$\mathsf{H}_{0,n}$ were recently established in~\cite{kock2025high}. For real numbers~$x_1,\hdots,x_n$, denote by~$x_1^*\leq \hdots\leq x_n^*$ their non-decreasing rearrangement. For any pair of winsorization points~$-\infty<\alpha\leq\beta<\infty$ we define the function 
\begin{equation*}
\phi_{\alpha,\beta}(x)
=
\begin{cases}
\alpha\qquad \text{if }x<\alpha,\\
x\qquad \text{if }x\in[\alpha,\beta],\\
\beta\qquad \text{if }x>\beta.
\end{cases}
\end{equation*} 
Next, for an~$\eps_n\in(0,1/2)$ that will be specified in~\eqref{eq:epsfam} further below, let~$\hat{\alpha}_j=\tilde X_{\lceil \eps_n n \rceil, n,j}^*$ and~$\hat{\beta}_j=\tilde X_{\lfloor(1-\eps_n )n\rfloor + 1,n,j}^*$, and define~$S_{n,W}^{\dagger}\in\R^d$ as
\begin{equation}\label{eqn:winsmeandef}
S_{n,W, j}^{\dagger}
=
n^{-1/2}\sum_{i=1}^n\phi_{\hat\alpha_j,\hat\beta_j}(\tilde{X}_{i,n,j}),\quad j=1,\hdots,d.
\end{equation}
Thus, for each coordinate~$j$, the winsorization points~$\hat{\alpha}_j$ and~$\hat{\beta}_j$ are empirical quantiles (i.e., order statistics) of the contaminated random variables~$\tilde{X}_{1,n,j},\hdots,\tilde{X}_{n,n,j}$.

Under adversarial contamination~\eqref{eq:contamfrac} it is clear that even~$S_{n,W}^{\dagger}$ can perform arbitrarily badly unless at least the smallest and largest~$\overline{\eta}_n n$ observations are winsorized. Thus, one must choose~$\eps_n\geq \overline{\eta}_n$. In particular, we study~$\eps_n$ of the form
\begin{equation}\label{eq:epsfam}
\eps_n=\lambda_{1}\cdot \overline{\eta}_n +\lambda_{2}\cdot \frac{\log(dn)}{n},\qquad\lambda_1\in(1,\infty)\text{ and }\lambda_2\in(0,\infty).
\end{equation}
Imposing~$\lambda_1>1$ ensures that~$\eps_n\geq \overline{\eta}_n$ and we recommend selecting~$\lambda_1$ only slightly greater than one to avoid winsorizing unnecessarily. In the special case where~$\overline{\eta}_n=0$ there is no contamination, which is implicitly imposed in Section~\ref{sec:arim} for the test~$\varphi_n$ based on~$S_n$. 

To properly normalize the entries of~$S_{n,W}^{\dagger}$, and to define~$\varphi_{n,W}$, consider the following estimator of~$\Sigma_n(P)$ from~\cite{kock2025high}. Let
\begin{equation}\label{eq:epsprime}
\eps_n'= \lambda_1'\cdot\overline{\eta}_n+\lambda_2'\cdot\frac{\log(d^2n)}{n},\qquad \lambda_1'\in(1,\infty)\text{ and }\lambda_2'\in (0,\infty).
\end{equation}
Writing~$\hat{a}_j=\tilde{X}_{\lceil \eps_n'n\rceil,n,j}^*$,~$\hat{b}_j=\tilde{X}_{\lfloor(1-\eps_n' )n\rfloor + 1,n,j}^*$ and~$\tilde{\mu}_{n,j}=n^{-1}\sum_{i=1}^n\phi_{\hat{a}_j,\hat{b}_j}(\tilde{X}_{i,n,j})$ for~$j=1,\hdots,d$, define~$\tilde{\Sigma}_n$ as the matrix with entries 
\begin{equation}\label{eq:tildeSigma}
\tilde{\Sigma}_{n,j,k}
=
\frac{1}{n}\sum_{i=1}^n\sbr[1]{\phi_{\hat{a}_j,\hat{b}_j}(\tilde{X}_{i,n,j})-\tilde{\mu}_{n,j}}\sbr[1]{\phi_{\hat{a}_k,\hat{b}_k}(\tilde{X}_{i,n,k})-\tilde{\mu}_{n,k}},\qquad 1\leq j,k\leq d.
\end{equation}
The (robust) covariance matrix estimator~$\tilde{\Sigma}_n$ is positive semi-definite and symmetric by virtue of being a Gram matrix. Let~$\tilde{\sigma}_{n,j}=\tilde{\Sigma}_{n,j,j}^{1/2}$,~$j=1,\hdots,d$, denote the square roots of the diagonal elements of~$\tilde{\Sigma}_n$, and define the winsorized analogue~$S_{n,W}$ to~$S_n$ as 
\begin{equation}\label{eq:S_nWS}
S_{n,W,j}
=
\frac{S_{n,W, j}^{\dagger}}{\tilde{\sigma}_{n,j}}
=
\frac{1}{\sqrt{n}\tilde{\sigma}_{n,j}}\sum_{i=1}^n\phi_{\hat\alpha_j,\hat\beta_j}(\tilde{X}_{i,n,j}),\quad j=1,\hdots,d,
\end{equation}
which, for the sake of completeness, we define as~$0$ in case~$\tilde{\sigma}_{n,j} = 0$. Finally, recalling the definition of~$c_{n,1-\alpha}$ immediately prior to~\eqref{eq:CVDiag}, we define the test
\begin{equation*}
\varphi_{n,W}
=
\mathds{1}\cbr[1]{\enVert[0]{S_{n,W}}_\infty>c_{n,1-\alpha}},\quad n\in\N;
\end{equation*}
(tests based on bootstrap critical values will be considered in Section~\ref{sec:BootstrapApprox} below). We always implement~$S_{n,W}$ with~$\eps_{n}$ as in~\eqref{eq:epsfam} and~$\eps_{n}'$ as in~\eqref{eq:epsprime} without further mentioning.\footnote{The sequences~$\eps_{n}$ and~$\eps_{n}'$ both converge to zero if~$\sqrt{n\log(d)}\overline{\eta}_{n}^{1-\frac{1}{m}}\to 0$, and~$\log(d)/{n}^{\frac{m-2}{5m-2}}\to 0$, which is imposed in Theorems~\ref{thm:Wins},~\ref{thm:BSWins}, and~\ref{thm:comparison}. As a result, it eventually holds that~$\hat{\alpha}_j\leq \hat{\beta}_j$ and~$\hat{a}_j\leq \hat{b}_j$ for~$j=1,\hdots,d$.} Our simulations in Section~\ref{sec:Sims} suggest that~$\lambda_1=\lambda_1'=1.01$,~$\lambda_2=1$, and~$\lambda_2'=0.01$ work well in practice. Appendix~\ref{sec:RobCheck} studies the sensitivity to these choices.
\begin{theorem}\label{thm:Wins}
Let~$b_1\in(0,\infty)$,~$b_2\in(b_1,\infty)$,~$m>2$,~$\alpha\in(0,1)$ and suppose~$d\geq 2$. Fix~$\lambda_1,\lambda_1'\in(1,\infty)$ and~$\lambda_2,\lambda_2'\in(0,\infty)$. If~$$\sqrt{n\log(d)}\overline{\eta}_{n}^{1-\frac{1}{m}}\to 0 \quad \text{ and  } \quad \log(d)/{n}^{\frac{m-2}{5m-2}}\to 0,$$ then the following holds:
\begin{enumerate}[leftmargin=*]
\item Asymptotic size control:
\begin{equation}\label{eq:sizeWins}
\limsup_{n\to\infty} \sup_{P\in\mc{P}^0(b_1,b_2,m)} P\del[1]{\enVert[0]{S_{n,W}}_\infty >c_{n,1-\alpha}}\leq \alpha,
\end{equation}
and~$\lim_{n \to \infty} P_n(\enVert[0]{S_{n,W}}_\infty >c_{n,1-\alpha}) = \alpha$ for every sequence~$P_n \in \mc{P}^0(b_1,b_2,m)$ satisfying~$\Sigma_{0, n}(P_n) = I_d$ for~$n$ large enough. If~$d \to \infty$, then~$\lim_{n \to \infty} P_n(\enVert[0]{S_{n,W}}_\infty >c_{n,1-\alpha}) = \alpha$ holds for every sequence~$P_n$ in~$\mc{P}^0(b_1,b_2,m)$ contained in some~$\mathcal{P}(b_1, b_2, m, \bm{r})$ (defined in~\eqref{eqn:Pcorrbd} below) for a sequence~$\bm{r} = (r_l)$ such that~$\log(l)r_l \to 0$ as $l\to \infty.$

\item Gaussian approximation for the power function:
\begin{equation}
\begin{aligned}\label{eq:GaussApproxWins}
\lim_{n\to\infty}\sup_{P\in\mc{P}(b_1,b_2,m)} \bigg|  &P \left( \enVert[0]{S_{n,W}}_\infty>c_{n,1-\alpha}\right) \\
& \hspace{1cm}-P\left(\enVert[0]{\Sigma_{0,n}^{1/2}(P)Z+\sqrt{n}D_n^{-1}(P)\mu_n(P)}_\infty>c_{n,1-\alpha}\right) \bigg|=0.
\end{aligned}
\end{equation}
\end{enumerate}

\end{theorem}
Note that the conditions in Theorem~\ref{thm:sums} are stronger than the ones in Theorem~\ref{thm:Wins}. Comparing the statements established in those theorems, it follows that whenever the stronger conditions on~$d$ (and implicitly on~$\overline{\eta}_n$) from Theorem~\ref{thm:sums} are satisfied, then~$\varphi_n$ and~$\varphi_{n,W}$ both have asymptotic size controlled by~$\alpha$ and have \emph{identical} asymptotic power against \emph{all alternatives}. In addition, the guarantees for~$\varphi_{n,W}$ established in Theorem~\ref{thm:Wins} are valid even 
\begin{enumerate}
\item when~$d$ grows at an exponential rate in~$n$, i.e.~$\log(d)=o\del[1]{n^{(m-2)/(5m-2)}}$, far exceeding the polynomial rate~$d/n^{m/2-1-\xi}\to 0$ for some~$\xi>0$ imposed for~$\varphi_n$ in Theorem~\ref{thm:sums}; and
\item under a modest amount of adversarial contamination, i.e., as long as~$\overline{\eta}_n > 0$ satisfies the condition 
\begin{equation}\label{eqn:contdec}
\sqrt{n\log(d)}\overline{\eta}_n^{1-\frac{1}{m}}\to 0.
\end{equation}
Note that for the contamination~$\overline{\eta}_n=2/n$ studied in Theorem~\ref{thm:Breakdown}, condition~\eqref{eqn:contdec} is satisfied for~$m>2$ if $\log(d)/{n}^{\frac{m-2}{5m-2}}\to 0$. Thus, although the test~$\varphi_n$ breaks down under~$\overline{\eta}_n=2/n$  in that it has size and power equal to zero (cf.~Theorem~\ref{thm:Breakdown}), the robustified max-test~$\varphi_{n,W}$ obeys~\eqref{eq:sizeWins} and~\eqref{eq:GaussApproxWins}. In particular,~$\varphi_{n,W}$ is consistent if, e.g.,$$\sqrt{n}\enVert[0]{\mu_n(P_n)}_\infty/\sqrt{\log(d)}\to\infty,$$ cf.~\eqref{eq:NoGainMeanCond} of Theorem~\ref{thm:comparison}.
\end{enumerate}
Therefore, basing a test on the winsorized means~$S_{n,W}$ instead of~$S_n$ comes without loss of (first order) asymptotic power, with asymptotic validity for much larger~$d$, and allows for~$\overline{\eta}_n>0$. Note also that Theorem~\ref{thm:Wins} only imposes~$m>2$, whereas Theorem~\ref{thm:sums} imposes~$m > 4$. That the censoring of the tails entailed by winsorization does \emph{not} imply a loss of asymptotic power is perhaps surprising.

The proof strategy of Theorem~\ref{thm:Wins} is similar to that of Theorem~\ref{thm:sums} (cf.~Section~\ref{sec:pf1}), and entails establishing Gaussian approximations to the distribution of~$S_{n,W}$ that are informative for non-zero~$\mu_n$,~$\overline{\eta}_n>0$, and~$d$ growing at an exponential rate in~$n$, even when only~$m>2$ moments are imposed to exist (cf.~Theorems~\ref{thm:SumGAWins} and~\ref{thm:SumGAWinsLarge}). 

\subsection{Necessity of~\eqref{eqn:contdec} for (first order) free robustness}\label{sec:necess}

Theorem~\ref{thm:Wins} requires, among other things, the contamination decay condition~\eqref{eqn:contdec}. This requirement is not merely an artefact of the winsorized statistic or of the proof of Theorem~\ref{thm:Wins}. Indeed, Lemma 2.2 and its discussion in a recently updated\footnote{Tue., 28 Jul 2026.} version of~\cite{resende2024robust} imply that the uniform Gaussian approximation in~\eqref{eq:GaussApproxWins} can fail for any non-randomized statistic if (cf.~his Equation (2.4))
\begin{equation}\label{eqn:contdecf}
\limsup_{n\to\infty}\sqrt{n\log(d)}\overline{\eta}_{n}^{1-\frac{1}{m}}>0.
\end{equation}
Even if~\eqref{eq:GaussApproxWins} fails, this does not rule out the possibility that, under~\eqref{eqn:contdecf} and the remaining assumptions of Theorem~\ref{thm:Wins},~$\varphi_{n,W}$ (or a possibly randomized modification thereof) with asymptotic size~$\alpha \in (0, 1)$ has an asymptotic power function that is, under the alternative, asymptotically nowhere smaller, but somewhere larger than the approximation
\begin{equation}\label{eqn:GBench}
P \mapsto P\left(\enVert[0]{\Sigma_{0,n}^{1/2}(P)Z+\sqrt{n}D_n^{-1}(P)\mu_n(P)}_\infty>c_{n,1-\alpha}\right) =: \mc{G}_n(P)
\end{equation}
from~\eqref{eq:GaussApproxWins}. If this were true, although the Gaussian approximation in~\eqref{eq:GaussApproxWins} would fail, first order robustness would still be for free, as one would be able to (even) ``weakly dominate'' the power function~\eqref{eqn:GBench}! Using a similar argument as in the proof of Lemma 2.2 in the updated version of \cite{resende2024robust} mentioned above, we show that this is impossible. To be precise, we exhibit alternatives~$P_n^1$ with~$\mu_n(P_n^1) \neq 0$ and an admissible contamination scheme, such that under~\eqref{eqn:contdecf} and the remaining assumptions of Theorem~\ref{thm:Wins} \emph{every} test~$\psi_n$ of asymptotic size~$\alpha\in(0,1)$ also has asymptotic power at most~$\alpha<\lim_{n\to\infty}\mc{G}_n(P_n^1)$ against~$P_n^1$. 

We work with the following family of distributions: Let~$V_1, \hdots, V_d$ be independent random variables that are uniformly distributed on~$\{-1, 1\}$, let~$B$ be a Bernoulli distributed random variable with success probability~$p \in [0, 1]$, and suppose that~$V := (V_1, \hdots, V_d)'$ and~$B$ are independent. Define the distribution~$Q_{p,d}$, say, on~$\R^d$ via ($\mathcal{L}$ applied to a random quantity denotes its distribution)
\begin{equation}\label{eqn:defqpd}
Q_{p,d}:= \begin{cases}
\mc{L}(b_1V) & \text{ if } p = 0, \\
\mc{L}(b_1V+ap^{-1/m} B\bm{\iota}_d) & \text{ else},
\end{cases}
\end{equation}
where~$0 < b_1 < b_2 < \infty$,~$a := (b_2 - b_1)/2$, and~$\bm{\iota}_d:=(1,\hdots,1)'$.

Lemma~\ref{lem:HDDistConstruct} in Appendix~\ref{sec:AppLB} verifies that the requirements on the distribution of~$X_{i,n}$ imposed in Section~\ref{sec:setting}, and encoded through~$\mathcal{P}(b_1, b_2, m)$, are \emph{compatible} with the family~$Q_{p,d}$. The following assumption ensures that the required members of this family are attainable as laws of $X_{1,n}$.\footnote{That a setting, i.e., a measurable space~$(\Omega, \mathcal{F})$ and an array of observations~$X_{1, n}, \hdots, X_{n,n}$ such that Assumption~\ref{ass:lbd} holds \emph{exists} is easy to see, e.g., by working with a suitable infinite product space.} The interval specified for $c$ is nonempty precisely when \eqref{eqn:contdecf} holds.
\begin{assumption}\label{ass:lbd}
For a constant $$c \in \left(0, \limsup_{n\to\infty}
\sqrt{n\log(d)}\overline{\eta}_n^{1-1/m}\right),$$ and every~$n \in \N$ greater than~$M > 0$, say, there exist~$P^0_n \in \mc{P}^0(b_1,b_2,m)$ and~$P_n^1 \in \mc{P}(b_1,b_2,m)$, such that~$P^0_n \circ X^{-1}_{1,n} = Q_{0, d}$ and~$P^1_n \circ X^{-1}_{1,n} = Q_{p_n, d}$ , where
\begin{equation}\label{eqn:qkpdef}
p_n:=\kappa q_n,
\qquad
q_n:=(n\log(d \vee 2))^{-\frac{m}{2(m-1)}}, 	\qquad 
\kappa:= c^{\frac{m}{m-1}}/2.
\end{equation}
\end{assumption}
Under Assumption~\ref{ass:lbd}, eventually~$\mu_n(P_n^1)=\frac{a\kappa^{(m-1)/m}}{\sqrt{n\log(d \vee 2)}}\bm\iota_d\neq 0_d$ so that~$P^1_n \notin \mc{P}^0(b_1,b_2,m),$ i.e.,~$P^1_n$ is an element of the alternative hypothesis (cf.~Theorem~\ref{thm:EtaNecessary} in Appendix~\ref{sec:AppLB}). Furthermore, Theorem~\ref{thm:EtaNecessaryshort} only makes use of \emph{deterministic} contamination schemes: for~$\eta\in[0,1]$, let~$\mathfrak{C}_n(\eta)$ denote the collection of all measurable maps~$\mathsf{C}:(\R^d)^n\to(\R^d)^n$, satisfying
\begin{equation*}
\envert[1]{\cbr[1]{i\in\cbr[0]{1,\hdots,n}:\mathsf{C}_i(x)\neq x_i}}
\leq \eta n
\qquad\text{for every }x\in(\R^d)^n.
\end{equation*}
Thus, for any~$\mathsf{C}\in \mathfrak{C}_n(\eta)$ the contaminated data~$\mathsf{C}(X_{1,n},\hdots,X_{n,n}) = (\tilde{X}_{1,n},\hdots,\tilde{X}_{n,n})$ satisfies~$\envert[1]{\cbr[1]{i\in\cbr[0]{1,\hdots,n}:\tilde{X}_{i,n}\neq X_{i,n}}}\leq \eta n$.
\begin{theorem}\label{thm:EtaNecessaryshort}
Let~$b_1\in(0,\infty)$,~$b_2\in(b_1,\infty)$,~$m>2$, and~$\alpha\in(0,1)$. Suppose Assumption~\ref{ass:lbd} is satisfied,~$d\to\infty$, and~$\log(d)/n^{\frac{m-2}{5m-2}}\to0$. Then, there exists a sequence~$\mathsf{C}_n$ of contamination schemes satisfying~$\mathsf{C}_n \in \mathfrak{C}_n(\overline{\eta}_n)$ for every~$n \in \N$, such that for every sequence of tests~$\psi_n:(\R^d)^n\to[0,1]$, asymptotic size control at level~$\alpha$, i.e.,
\begin{equation}\label{eqn:asysizeshort}
\limsup_{n\to\infty}
\sup_{P\in\mc{P}^0(b_1,b_2,m)} 
E_P\sbr[1]{\psi_n(\mathsf{C}_n(X_{1,n}, \hdots, X_{n,n}))}
\leq\alpha,
\end{equation}
implies (no asymptotic power against the alternative~$P_{n}^1$ along a subsequence)
\begin{equation}\label{eqn:powshort}
\liminf_{n\to\infty}
E_{P_n^1}\sbr[1]{\psi_n(\mathsf{C}_n(X_{1,n}, \hdots, X_{n,n}))}
\leq\alpha < \lim_{n \to \infty} \mc{G}_n(P_n^1).
\end{equation}
\end{theorem}
Theorem~\ref{thm:EtaNecessaryshort} shows that when~\eqref{eqn:contdecf} holds, i.e.,~$$\limsup_{n\to\infty}\sqrt{n\log(d)}\overline{\eta}_{n}^{1-\frac{1}{m}}>0,$$ but the remaining conditions of Theorem~\ref{thm:Wins} are satisfied (as well as~$d\to \infty$), then \emph{any} test~$\psi_n$ of asymptotic size~$\alpha\in(0,1)$ has asymptotic power at most~$\alpha$ against the alternative~$P_n^1$ (along a subsequence). On the other hand, if instead~\eqref{eqn:contdec} holds, i.e., ~$$\lim_{n\to\infty}\sqrt{n\log(d)}\overline{\eta}_{n}^{1-\frac{1}{m}}=0,$$ then by~\eqref{eq:GaussApproxWins} of Theorem~\ref{thm:Wins}, the test~$\varphi_{n,W}$ (asymptotic size~$\alpha$) satisfies 
\begin{equation*}
\lim_{n\to\infty}E_{P_n^1}\varphi_{n,W}=\lim_{n\to\infty}\mc{G}_n(P_n^1)>\alpha.
\end{equation*}
Thus, if~\eqref{eqn:contdec} is satisfied then the winsorized test~$\varphi_{n,W}$ has asymptotic power exceeding its asymptotic size~$\alpha$ against~$P_n^1$, but if instead~\eqref{eqn:contdecf} is satisfied then \emph{no test} of asymptotic size~$\alpha$ can have asymptotic power exceeding~$\alpha$ against~$P_n^1$ (along a subsequence). Hence, the loss of power is caused by~\eqref{eqn:contdec} no longer holding and therefore no test can achieve robustness for free when this is the case.

The proof of Theorem~\ref{thm:EtaNecessaryshort} follows from  Theorem~\ref{thm:EtaNecessary} of Appendix~\ref{sec:AppLB}, which shows that for suitable~$\mathsf{C}_n\in\mathfrak{C}_n(\overline{\eta}_n)$ the total variation distance between~$P_{n}^0\circ (\mathsf{C}_n\circ (X_{1,n}, \hdots, X_{n,n}))^{-1}$ and~$P_{n}^1\circ (\mathsf{C}_n\circ (X_{1,n}, \hdots, X_{n,n}))^{-1}$ converges to zero along a subsequence from which it follows that~\eqref{eqn:asysizeshort} implies the first (i.e.~the weak) inequality of~\eqref{eqn:powshort}. The second inequality of~\eqref{eqn:powshort} follows by extreme value calculations for dependent Gaussians. Further details on the proof strategy can be found in Remark~\ref{rem:necproof} of Appendix~\ref{sec:AppLB}.  

\section{Data-driven (bootstrap) critical values}\label{sec:BootstrapApprox}

The test~$\varphi_{n,W}=\mathds{1}\cbr[0]{\enVert[0]{S_{n,W}}_\infty>c_{n,1-\alpha}}$ uses a critical value~$c_{n,1-\alpha}$ that does \emph{not} take into account the correlation structure of the entries of~$X_{1,n}$. This possibly leaves some room for improvement in terms of power. In this section, we show how one can replace~$c_{n,1-\alpha}$ by data-driven (bootstrap) critical values, taking into account the correlation structure of the entries of~$X_{1,n}$. We characterize the asymptotic power function of the resulting test~$\varphi_{n,W,B}$, say, for a large family of correlation structures. Finally, we investigate the question when --- and when not --- $\varphi_{n,W,B}$ has higher asymptotic power than~$\varphi_{n,W}$.\footnote{In case~$d/n^{m/2-1-\xi}\to 0$ for some~$\xi>0$ and~$\overline{\eta}_n=0$, one could also consider data-dependent critical values for tests based on the supremum of arithmetic means, i.e.,~$\|S_n\|_\infty$. However, since we are interested in situations wherein~$d$ increases at an exponential rate in~$n$ in the presence of~$m>2$ moments only and some contamination may be present, we refrain from studying this. We refer to~\cite{chang2017simulation} for results on size control and consistency against certain alternatives for tests based on~$\|S_n\|_\infty$ with data-driven critical values when~$d/n^{m/2-1-\xi}\to 0$ for some~$\xi>0$ and~$\overline{\eta}_n=0$.\label{fn:boot}} 

To define~$\varphi_{n,W,B}$, we define the correlation matrix~$\tilde{\Sigma}_{0,n}$ corresponding to~$\tilde{\Sigma}_{n}$ via
\begin{equation}\label{eqn:corrmatdef}
\tilde{\Sigma}_{0, n, i, j} := 
\begin{cases}
\tilde{\Sigma}_{n, i, j}/\sqrt{\tilde{\Sigma}_{n, i,i}\tilde{\Sigma}_{n,j,j}} & \text{ if } \tilde{\Sigma}_{n,i,i} \times \tilde{\Sigma}_{n,j,j} > 0, \\
0 & \text{ if } \tilde{\Sigma}_{n,i,i} \times \tilde{\Sigma}_{n,j,j} = 0 \text{ and } i \neq j, \\
1 & \text{ if } \tilde{\Sigma}_{n,i,i}= 0 \text{ and } i = j,
\end{cases}
\end{equation}	
noting that in the leading case where all diagonal elements of~$\tilde{D}_n:=\text{diag}\del[0]{\tilde{\sigma}_{n,1},\hdots,\tilde{\sigma}_{n,d}}$ are positive, we have~$\tilde{\Sigma}_{0,n}=\tilde{D}_n^{-1}\tilde{\Sigma}_n\tilde{D}_n^{-1}$. Mainly for completeness, the remaining cases in~\eqref{eqn:corrmatdef} take care of vanishing variance estimates.

Furthermore, given~$\alpha \in (0, 1)$, let~$c_{n,1-\alpha}^B$ be the~$(1-\alpha)$-quantile of the conditional distribution of~$\|\tilde{\Sigma}^{1/2}_{0,n}Z\|_\infty$ given~$\tilde{X}_{1,n},\hdots,\tilde{X}_{n,n}$, that is\footnote{Because~$\tilde{\Sigma}_{0,n}$ is, by definition, a correlation matrix, the distribution of~$\|\tilde{\Sigma}^{1/2}_{0,n}Z\|_\infty$ is (Lebesgue) absolutely continuous, and equality in~\eqref{eqn:bootcdef} can be achieved as required by suitably choosing the critical value.}

\begin{equation}\label{eqn:bootcdef}
P\left(\|\tilde{\Sigma}^{1/2}_{0,n}Z\|_\infty>c_{n,1-\alpha}^B\mid \tilde{X}_{1,n},\hdots,\tilde{X}_{n,n}\right)=\alpha.
\end{equation} 
We also denote by~$c_{n,1-\alpha}(\Sigma_{0,n}(P))$ the~$(1-\alpha)$-quantile of the distribution of~$\|\Sigma_{0,n}^{1/2}(P)Z\|_\infty$, that is 
\begin{equation*}
P\left(\|\Sigma_{0,n}^{1/2}(P)Z\|_\infty>c_{n,1-\alpha}(\Sigma_{0,n}(P))\right)=\alpha.
\end{equation*} Finally, set
\begin{equation*}
\varphi_{n,W,B}	
=
\mathds{1}\cbr[1]{\enVert[0]{S_{n,W}}_\infty>c_{n,1-\alpha}^B},\qquad n\in\N.
\end{equation*}

To state the asymptotic approximation of the power function of~$\varphi_{n,W,B}$ in Part 2~of Theorem~\ref{thm:BSWins} below, we introduce the following notation. For any sequence~$\bm{r} = (r_l)_{l \in \N}$ in~$[0, 1)$, set
\begin{equation}\label{eqn:Pcorrbd}
\mc{P}(b_1,b_2,m,\bm{r})
=
\cbr[2]{P\in \mc{P}(b_1,b_2,m):|\Sigma_{0,n,i,j}(P)|\leq r_{|i-j|}\text{ for all }1\leq i<j\leq d,\ n\in\N}.
\end{equation}
Therefore, a~$P\in \mc{P}(b_1,b_2,m)$ is in~$\mc{P}(b_1,b_2,m,\bm{r})$ if, for every~$l \in \N$, the~$l$-th off-diagonal elements of the correlation matrices~$\Sigma_{0,n}(P)$ are uniformly (in~$n$) bounded by~$r_l$. Hence, the sequence~$\bm{r}$ can be viewed as an upper bound on the decay of the correlations between the coordinates of~$X_{i,n}$ as the distance between their indices increases. 

For some of our results below, we have to assume that~$\bm{r}$ satisfies~$\log(l)r_l\to 0$ as~$l \to \infty$. Note that the resulting sets~$\mc{P}(b_1,b_2,m,\bm{r})$ allow for a wide, yet not exhaustive, range of correlated vectors~$X_{i,n}$. Decay conditions of that type were used already in~\cite{berman1964limit} to show that~$\max_{1\leq j\leq d}(\Sigma_{0,n}^{1/2}Z)_j$ has the same limiting distribution as~$\max_{1\leq j\leq d}Z_j$ when~$\Sigma_{0,n}$ is ``stationary'' in the sense of~$\Sigma_{0,n,i,j}=r_{|i-j|}$ for some sequence~$r_l$ with~$\log(l)r_l\to 0$; cf.~also Chapter 4 of \cite{leadbetter} for an overview of extreme value theory for correlated Gaussians, which we exploit in the proofs of Theorems~\ref{thm:BSWins} and~\ref{thm:comparison} below.

The following result establishes size guarantees for~$\varphi_{n,W,B}$ and an asymptotic approximation to its power function. 
\begin{theorem}\label{thm:BSWins}
Let~$b_1\in(0,\infty)$,~$b_2\in(b_1,\infty)$,~$m>2$,~$\alpha\in(0,1)$ and suppose~$d\geq 2$. Fix~$\lambda_1,\lambda_1'\in(1,\infty)$ and~$\lambda_2,\lambda_2'\in(0,\infty)$. If~$$\sqrt{n\log(d)}\overline{\eta}_{n}^{1-\frac{1}{m}}\to 0 \quad \text{ and  } \quad \log(d)/{n}^{\frac{m-2}{5m-2}}\to 0,$$ then the following holds:
\begin{enumerate}[leftmargin=*]
\item The rejection probability of~$\varphi_{n,W,B}$ converges to~$\alpha$ uniformly over~$\mc{P}^0(b_1,b_2,m)$: 
\begin{equation}\label{eq:BSsize}
\lim_{n\to\infty}\sup_{P\in\mc{P}^0(b_1,b_2,m)}\envert[2]{P\del[1]{\enVert[0]{S_{n,W}}_\infty>c_{n,1-\alpha}^B}-\alpha}
=
0.
\end{equation}
\item Suppose that~$d\to\infty$. For all sequences~$\bm{r} = (r_l)_{l \in \N}$ in~$[0, 1)$, such that~$\log(l)r_l\to 0$ as~$l \to \infty$ and~$\mc{P}(b_1,b_2,m, \bm{r}) \neq \emptyset$, it holds that 
\begin{equation}\label{eq:BSPower}
\begin{aligned}
\lim_{n\to\infty}\sup_{P\in\mc{P}(b_1,b_2,m, \bm{r})} \bigg|  &P \left( \enVert[0]{S_{n,W}}_\infty>c^B_{n,1-\alpha}\right) \\
& -P\left(\enVert[0]{\Sigma_{0,n}^{1/2}(P)Z+\sqrt{n}D_n^{-1}(P)\mu_n(P)}_\infty>c_{n,1-\alpha}(\Sigma_{0, n}(P))\right) \bigg|=0.
\end{aligned}
\end{equation}
\end{enumerate}
\end{theorem}
Clearly,~\eqref{eq:BSsize} of Theorem~\ref{thm:BSWins} implies that for all~$P\in \mc{P}^0(b_1,b_2,m)$, the rejection frequency of~$\varphi_{n,W,B}$ converges to~$\alpha$. For comparison, note that for the test~$\varphi_{n,W}$, which does not use data-driven critical values,~\eqref{eq:sizeWins} of Theorem~\ref{thm:Wins} does \emph{not} rule out that for some~$P\in \mc{P}^0(b_1,b_2,m)$ it is the case that~$\limsup_{n\to\infty}E_P\varphi_{n,W}$ is strictly smaller than~$\alpha$, as it only controls the \emph{maximal} rejection frequency over~$\mc{P}^0(b_1,b_2,m)$. In fact, it is not difficult to construct~$P\in\mc{P}^0(b_1,b_2,m)$ for which~$\lim_{n\to\infty}E_P\varphi_{n,W}=0$.

Next,~\eqref{eq:BSPower} establishes an asymptotic Gaussian approximation to the power function of~$\varphi_{n,W,B}$ similar to that for~$\varphi_{n,W}$ in~\eqref{eq:GaussApproxWins} but with~$c_{n,1-\alpha}(\Sigma_{0,n}(P))$ replacing the no smaller~$c_{n,1-\alpha}$ (cf. the opening discussion in the following section). The proof of Part~1.~of Theorem~\ref{thm:comparison} below reveals that~\eqref{eq:BSPower} remains valid even with~$c_{n,1-\alpha}$ in place of~$c_{n,1-\alpha}(\Sigma_{0,n}(P))$. Note also that the supremum in~\eqref{eq:BSPower} is taken over~$\mc{P}(b_1,b_2,m,\bm{r})\subseteq \mc{P}(b_1,b_2,m)$.

\subsection{Comparing the power of~$\varphi_{n,W}$ and~$\varphi_{n,W,B}$}

It is easy to see that the power function of the test~$\varphi_{n,W,B}$ uniformly dominates the power function of~$\varphi_{n,W}$. To this end, just note that by the Khatri-{\v{S}}id{\'a}k inequality as well as independence of~$Z$ and~$\tilde{X}_{1,n},\hdots,\tilde{X}_{n,n}$ for~$P\in \mc{P}(b_1,b_2,m)$,
\begin{equation*}
P\left(\|\tilde{\Sigma}^{1/2}_{0,n}Z\|_\infty\leq c_{n,1-\alpha}\mid \tilde{X}_{1,n},\hdots,\tilde{X}_{n,n}\right)
\geq
P(\|Z\|_\infty\leq c_{n,1-\alpha})
=1-\alpha.	
\end{equation*}
Thus,~$c_{n,1-\alpha}^B\leq c_{n,1-\alpha}$, so that for every probability measure~$P$ on~$(\Omega,\mc{F})$
\begin{equation}\label{eqn:powermono}
E_P\varphi_{n,W} \leq E_P\varphi_{n,W,B}; 
\end{equation}
that is for every~$n\in\N$ the rejection frequency under the null and alternative is \emph{no smaller} for~$\varphi_{n,W,B}$ than it is for~$\varphi_{n,W}$. Let us emphasize that the relation in~\eqref{eqn:powermono} alone does not imply that $\varphi_{n,W,B}$ would be ``better'' than~$\varphi_{n,W}$, as not only the power but also the size of~$\varphi_{n,W,B}$ is typically larger than that of~$\varphi_{n,W}$. The relation becomes more informative from an asymptotic point of view: according to Theorems~\ref{thm:Wins} and~\ref{thm:BSWins} both tests~$\varphi_{n,W}$ and~$\varphi_{n,W,B}$ control asymptotic size ($\varphi_{n,W}$ can be asymptotically conservative, while~$\varphi_{n,W,B}$ provides exact asymptotic size control). Relation~\eqref{eqn:powermono} now implies that the (first order) asymptotic power of~$\varphi_{n,W,B}$ is never lower than that of~$\varphi_{n,W}$. Of course, this statement does not rule
out that the asymptotic power of these two tests is identical for many distributions, which turns out to be the case and which is established in the following result.
\begin{theorem}\label{thm:comparison}
Let~$b_1\in(0,\infty)$,~$b_2\in(b_1,\infty)$,~$m>2$,~$\alpha\in(0,1)$ and~$d\geq 2$. Fix~$\lambda_1,\lambda_1'\in(1,\infty)$ and~$\lambda_2,\lambda_2'\in(0,\infty)$. If~$$\sqrt{n\log(d)}\overline{\eta}_{n}^{1-\frac{1}{m}}\to 0 \quad \text{ and  } \quad \log(d)/{n}^{\frac{m-2}{5m-2}}\to 0,$$ then the following holds:
\begin{enumerate}[leftmargin=*]
\item Suppose that~$d\to\infty$. For all sequences~$\bm{r} = (r_l)_{l \in \N}$ in~$[0, 1)$, such that~$\log(l)r_l\to 0$ and~$\mc{P}(b_1,b_2,m, \bm{r}) \neq \emptyset$, it holds that 
\begin{equation}\label{eq:NoGain}
\lim_{n\to\infty}\sup_{P\in\mc{P}(b_1,b_2,m,\bm{r})}\envert[2]{P\del[1]{\enVert[0]{S_{n,W}}_\infty > c_{n,1-\alpha}}-P\del[1]{\enVert[0]{S_{n,W}}_\infty > c^B_{n,1-\alpha}}}
=  
0.	
\end{equation}

\item For~$P_n\in\mc{P}(b_1,b_2,m)$ such that~$\sqrt{n}\enVert[0]{\mu_n(P_n)}_\infty/\sqrt{\log(d)}\to\infty$, it holds that
\begin{equation}\label{eq:NoGainMeanCond}
\lim_{n\to\infty}P_n\del[1]{\enVert[0]{S_{n,W}}_\infty > c_{n,1-\alpha}}= 1\quad\text{and}\quad \lim_{n\to\infty}P_n\del[1]{\enVert[0]{S_{n,W}}_\infty > c^B_{n,1-\alpha}}=1.
\end{equation}
\item Suppose~$d\to\infty$. For~$P_n\in\mc{P}(b_1,b_2,m)$ such that~$\mu_n(P_n) = \mathfrak{m}_n \bm{\iota}_d$ and~$\Sigma_{n}(P_n) = b_1^2 \bm{\iota}_d \bm{\iota}_d'$, where~$\bm{\iota}_d = (1, \hdots, 1)' \in \R^d$ and~$\mathfrak{m}_n=\log^{1/4}(d)/\sqrt{n}$ for~$n$ large enough, it holds that
\begin{equation}\label{eq:Gain}
\lim_{n\to\infty}P_n\del[1]{\enVert[0]{S_{n,W}}_\infty > c_{n,1-\alpha}}= 0 
\quad\text{and}\quad 
\lim_{n\to\infty}P_n\del[1]{\enVert[0]{S_{n,W}}_\infty > c^B_{n,1-\alpha}}= 1.
\end{equation}
\end{enumerate}	
\end{theorem}

Parts 1.~and 2.~of Theorem~\ref{thm:comparison} provide sufficient conditions on the correlation matrix~$\Sigma_{0,n}$ and deviation from the null~$\mu_n$, respectively, under which the difference in the asymptotic power of~$\varphi_{n,W}$ and~$\varphi_{n,W,B}$ converges to zero. The condition in Part 1.~requires the correlation to decay to~$0$ at a certain (slow) rate. The condition in Part~2.~requires the deviation from the null~$\mu_n$ to be non-negligible in a certain way. Under each of these conditions, there is \emph{no asymptotic power gain} from using the data-driven critical values~$c_{n,1-\alpha}^B$ instead of~$c_{n,1-\alpha}$. Specifically, Part~1.~shows that even when the entries of~$X_{1,n}$ exhibit substantial correlation (in the sense of the off-diagonal elements of~$\Sigma_{0,n}$ decaying at a logarithmic rate only as the distance between their indices increases),~$\varphi_{n,W}$ and~$\varphi_{n,W,B}$ have the same asymptotic power. This essentially follows from Parts~2.~of Theorems~\ref{thm:Wins} and~\ref{thm:BSWins},~$c_{n,1-\alpha}(\Sigma_{0,n})=c_{n,1-\alpha}+o(1/\sqrt{\log(d)})$ for such~$\Sigma_{0,n}$ (cf.~Lemma~\ref{lem:extreme}), along with Gaussian anti-concentration. Part~2.~of Theorem~\ref{thm:comparison} shows that for sufficiently large deviations~$\mu_n$ from the null, it does not matter which critical values one uses since the tests~$\varphi_{n,W}$ and~$\varphi_{n,W,B}$ will both be consistent.

Finally, (the proof of) Part 3.~guarantees the existence of very highly correlated~$X_{1,n}$ under which~$c_{n,1-\alpha}(\Sigma_{0,n})$ is bounded (because the strong correlation essentially implies the effective dimension of the data to be one-dimensional). Recalling that eventually~$c_{n,1-\alpha}\geq \sqrt{\log(d)}\to\infty$, this allows one to construct alternatives such that~\eqref{eq:Gain} is satisfied, i.e.,~such that the asymptotic power of~$\varphi_{n,W}$ is zero, yet~$\varphi_{n,W,B}$ is consistent.

\begin{remark}\label{rem:sumspowercomp}
Since~$\varphi_n$ and~$\varphi_{n,W}$ have identical asymptotic power whenever the stronger conditions on~$d$ and~$\overline{\eta}_n$ underlying the former in Theorem~\ref{thm:sums} are satisfied, the  comparison between~$\varphi_{n,W}$ and~$\varphi_{n,W,B}$ in Theorem~\ref{thm:comparison} implies analogous comparison statements between~$\varphi_{n}$ and~$\varphi_{n,W,B}$ under these stronger conditions.
\end{remark}

\section{Numerical evidence}\label{sec:Sims}
In this section we numerically investigate the finite-sample performance of the tests~$\varphi_{n}$, $\varphi_{n,W}$, and~$\varphi_{n,W,B}$ studied in Sections~\ref{sec:arim},~\ref{sec:winsm}, and~\ref{sec:BootstrapApprox}, respectively.\footnote{With a slight abuse of notation, we stress that~$\varphi_n$ is implemented on the potentially contaminated data. That is~$\varphi_n$ coincides with~$\tilde{\varphi}_n$ defined before~\eqref{eq:notationcontamsum}.} Furthermore, let~$\check{\Sigma}_{0,n}$ denote the usual sample correlation matrix of the vectors~$\tilde{X}_{1,n},\hdots,\tilde{X}_{n,n}$ (with an analogous convention for dealing with vanishing variance estimates as in~\eqref{eqn:corrmatdef}). Given~$\alpha \in (0, 1)$, let~$c_{n,1-\alpha}^{B,\text{sum}}$ be the~$(1-\alpha)$-quantile of the conditional distribution of~$\|\check{\Sigma}^{1/2}_{0,n}Z\|_\infty$ given~$\tilde{X}_{1,n},\hdots,\tilde{X}_{n,n}$, that is
\begin{equation*}
P\left(\|\check{\Sigma}^{1/2}_{0,n}Z\|_\infty>c_{n,1-\alpha}^{B,\text{sum}}\mid \tilde{X}_{1,n},\hdots,\tilde{X}_{n,n}\right)=\alpha.
\end{equation*}
We also implement the test
\begin{equation*}
\varphi_{n,B}=\mathds{1}\cbr[1]{\enVert[0]{S_n}_\infty>c_{n,1-\alpha}^{B,\text{sum}}},\quad n\in\N,
\end{equation*} 
which is the bootstrap version of~$\varphi_n$ studied in~\cite{chang2017simulation} (absent contamination). Finally, we invert the confidence sets for high-dimensional population means developed in \cite{liu2024robust} and \cite{gao2021refinedcramertypemoderate} to obtain further benchmarks~$\varphi_{n,LL}$ and~$\varphi_{n,GSS}$, respectively. We stress that the theoretical size guarantees for $\varphi_{n,LL}$ and~$\varphi_{n,GSS}$ were developed in the setting of no contamination ($\overline{\eta}_n=0$), which explains why they deteriorate in the presence of contamination. As pointed out in the conclusion of~\cite{liu2024robust}, it would be of interest to theoretically study~$\varphi_{n,LL}$ in a setting with contaminated data.

In all simulations below, the uncontaminated~$X_{i,n}$ are i.i.d.~multivariate~$t$-distributed with scale matrix~$\Omega_{n}$ with~$\Omega_{n,i,j}=\rho^{|i-j|}$ for~$\rho\in\cbr[0]{0,0.9,1}$ and~$\nu=4.01$ degrees of freedom, implying that the covariance matrix is~$\Sigma_n=\frac{\nu}{\nu-2}\cdot\Omega_n$ and the entries of the~$X_{i,n}$ having finite~$m$th moment for~$m<\nu=4.01$. The case of~$\rho=1$ reflects a setting wherein the ``effective dimension'' is one since the entries of each~$X_{i,n}$ have a correlation of one. This case of extreme dependence is included to illustrate Part 3.~of Theorem~\ref{thm:comparison}. We use~$\alpha=0.05$ and all results are based on~$2{,}000$ Monte Carlo replications. Within each replication,~$c_{n,1-\alpha}^B$ and~$c_{n,1-\alpha}^{B,\text{sum}}$ are chosen as the empirical~$(1-\alpha)$-quantiles based on $500$ draws from~$\|\tilde{\Sigma}^{1/2}_{0,n}Z\|_\infty$ and~$\|\check{\Sigma}^{1/2}_{0,n}Z\|_\infty$ (obtained via multiplier bootstrap), respectively.\footnote{We experimented with basing the calculation of~$c_{n,1-\alpha}^B$ and~$c_{n,1-\alpha}^{B,\text{sum}}$ on fewer or more than~$500$ draws from ~$\|\tilde{\Sigma}^{1/2}_{0,n}Z\|_\infty$ and~$\|\check{\Sigma}^{1/2}_{0,n}Z\|_\infty$, respectively. However, this did not make any meaningful difference as this simulation error ``averaged out'' over the $2{,}000$ Monte Carlo replications.}

Throughout, we implement~$\varphi_{n,W}$ and~$\varphi_{n,W,B}$ with~$\lambda_1=\lambda_1'=1.01$,~$\lambda_2=1$, and~$\lambda_2'=0.01$. We choose~$\lambda_2'$ smaller than~$\lambda_2$ since smaller~$\lambda_2'$ leads to smaller~$\eps_n'$ and hence winsorizing fewer observations in the construction of~$\tilde{\sigma}_{n,j}$. This, in turn, leads to a larger denominator of the test statistic~$S_{n,W,j}=\frac{S^{\dagger}_{n, W,j}}{\tilde{\sigma}_{n,j}}$ [cf.~the definitions around \eqref{eq:S_nWS}]. Thus, for the test~$\varphi_{n,W}$, whose critical value~$c_{n,1-\alpha}$ is constant, lower~$\lambda_2'$ mechanically leads to a lower rejection frequency and hence a more conservative test, which is typically preferred over a liberal one. For~$\varphi_{n,W,B}$, its bootstrap critical value~$c_{n,1-\alpha}^B$ also changes with~$\lambda_2'$ through the estimated correlation matrix, so no analogous monotonicity is implied. We stress that there is nothing special about the exact value of~$\lambda_2'=0.01$ as $\varphi_{n,W}$ and~$\varphi_{n,W,B}$ only depend on~$\lambda_2'$ via the induced winsorization indices. For example, in the simulation settings studied, all~$\lambda_2'\in(0,0.04]$ lead to numerically \emph{identical} implementations of~$\varphi_{n,W}$ and~$\varphi_{n,W,B}$ as they all induce the same number of winsorized observations. Tables~\ref{tab:lambda_q0} and~\ref{tab:lambda_q3} in Section \ref{sec:RobCheck} of the supplementary appendix study the effect of larger~$\lambda_2'$ on null rejection probabilities, which are increasing in~$\lambda_2'$ and thus justify the current ``small'' choice. 

In the settings \emph{with} contamination that we study, the adversary replaces $q=3$ randomly chosen~$X_{i,n}$ by~$\tilde{X}_{i,n}=(500,\hdots,500)'$ but leaves the remaining~$X_{i,n}$ untouched.~$\varphi_{n,W}$ and~$\varphi_{n,W,B}$ are implemented with~$\overline{\eta}_n=q/n$. In Section \ref{sec:RobCheck} of the supplementary appendix, cf.~Table~\ref{tab:q_mis} there, we study the sensitivity of $\varphi_{n,W}$ and~$\varphi_{n,W,B}$ to misspecifying~$\overline{\eta}_n$ (or, equivalently,~$q$).

\subsection{Null rejection frequencies}
First, Table~\ref{tab:1} considers empirical null rejection frequencies absent contamination for all combinations of~$n\in\cbr[0]{100, 200, 500,1{,}000}$ and~$d\in\cbr[0]{200,500,1{,}000,2{,}000}$. The table reveals that apart from occasional over-rejections for~$n=100$ and~$n=200$, all tests appear to control size reasonably well.  Furthermore, by construction (cf.~the Khatri-{\v{S}}id{\'a}k inequality) the null rejection frequencies of the \emph{non}-bootstrap based tests~$\varphi_{n}$ and~$\varphi_{n,W}$ using the conservative critical values~$c_{n,1-\alpha}$ are lower than those of the bootstrap based counterparts $\varphi_{n,B}$ and~$\varphi_{n,W,B}$, respectively (in particular when~$\rho\in\cbr[0]{0.9,1}$). As~$\rho$ increases, the rejection frequency of $\varphi_{n}$ and~$\varphi_{n,W}$ decrease since the critical values~$c_{n,1-\alpha}$ become more conservative as they are based on the ``worst-case'' setting of~$\Sigma_{0,n}=I_d$, corresponding to~$\rho=0$. For~$\rho=1$, the rejection frequency of $\varphi_{n}$ and~$\varphi_{n,W}$ is essentially zero. Since the critical values of the bootstrap based tests $\varphi_{n,B}$ and~$\varphi_{n,W,B}$ ``adapt'' to~$\rho$, the rejection frequencies of these tests are often closer to~$\alpha=0.05$ (in particular for~$\rho=1$). The non-bootstrap based test~$\varphi_{n,GSS}$ is conservative for large~$\rho$ in particular, whereas the bootstrap-based~$\varphi_{n,LL}$ has nominal rejection frequencies closer to $\alpha=0.05$. Finally, we note that the good control of the null rejection frequencies of~$\varphi_n$ in this setting \emph{absent} contamination can be explained by the self-normalized nature of~$S_n$ and the results in~\cite{qiu2025self}, cf.~also Remark~\ref{rem:qiuchenshao}.

Table~\ref{tab:2} considers the setting of contaminated data. As expected, the tests~$\varphi_n$ and $\varphi_{n,B}$ based on~$S_n$ perform poorly: the presence of $q=3$ large ``outliers'' in each coordinate means that~$|S_{n,j}|\approx \sqrt{\frac{qn}{n-q}}\approx \sqrt{q}\approx 1.73$ for~$j=1,\hdots, d$, explaining the observed null rejection frequency of essentially zero for these tests as the critical values are larger than this.\footnote{Even when~$\rho=1$ such that~$\Sigma_{0,n}=\bm{\iota}_d \bm{\iota}_d'$ one has that~$c_{n,1-\alpha}(\Sigma_{0,n}(P))=1.96>1.73$. Thus, also the bootstrap critical values typically exceed~$1.73$, explaining the empirical rejection frequencies of essentially zero of~$\varphi_{n}$ and~$\varphi_{n,B}$. Note also that since~$\max_{1\leq j\leq d}|S_{n,j}|\approx \sqrt{q}$ increases in~$q$ one can expect rejection frequencies of close to one for~$\varphi_{n,B}$ for larger~$q$ for~$\rho = 1$. Indeed unreported simulations yield rejection frequencies of essentially one for~$\varphi_{n,B}$ already for~$q=5$ (without changing the value of the contaminated observations).} On the other hand, although the tests based on winsorized means have inflated null rejection frequencies when~$n=100$, they are much less affected by the contamination (and, as we shall see, have much higher power). The fact that~$\varphi_{n,W}$ has a null rejection frequency of zero when~$\rho=1$ is not due to contamination but a result of the critical value~$c_{n,1-\alpha}$ being too conservative as already observed in Table~\ref{tab:1} absent contamination. Finally,~$\varphi_{n,LL}$ is oversized even for large~$n$ (in particular for large~$\rho$). The opposite is the case for~$\varphi_{n,GSS}$, but recall that these latter two tests are not developed for a setting with adversarial contamination but instead designed to handle heavy tails. Overall,~$\varphi_{n,W,B}$ has empirical rejection closest to~$\alpha=0.05$.

Tables~\ref{tab:nu301_q0} and~\ref{tab:nu301_q3} of Section~\ref{sec:RobCheck} of the supplementary appendix repeat the size investigation of Tables~\ref{tab:1} and~\ref{tab:2}, but for the more heavy-tailed setting of~$\nu=3.01$ instead of~$\nu=4.01$. The findings are qualitatively identical and further comments are provided in the appendix.   
\begin{landscape}
\begin{table}[p]
\centering
\setlength{\tabcolsep}{4pt}
\renewcommand{\arraystretch}{0.88}
\begin{minipage}[t]{0.40\linewidth}
\centering
{\footnotesize 
\makebox[\linewidth][c]{$\rho=0$}\par\vspace{2pt}

\resizebox{\linewidth}{!}{%
\begin{tabular}{cccccccc}
\toprule
$n$&$d$& $\varphi_n$ & $\varphi_{n,B}$ & $\varphi_{n,W}$ & $\varphi_{n,W,B}$ & $\varphi_{n,LL}$ & $\varphi_{n,GSS}$ \\ 
\midrule
\multirow{4}{*}{$100$}
&$200$ & 0.058 & 0.063 & 0.031 & 0.034 & 0.058 & 0.004 \\ 
&$500$ & 0.066 & 0.080 & 0.037 & 0.043 & 0.066 & 0.002 \\ 
&$1000$ & 0.059 & 0.077 & 0.035 & 0.043 & 0.063 & 0.001 \\ 
&$2000$ & 0.069 & 0.085 & 0.037 & 0.042 & 0.066 & 0.000 \\
\midrule 
\multirow{4}{*}{$200$}
&$200$ & 0.046 & 0.051 & 0.025 & 0.025 & 0.046 & 0.017 \\ 
&$500$ & 0.059 & 0.066 & 0.034 & 0.036 & 0.058 & 0.014 \\ 
&$1000$ & 0.058 & 0.065 & 0.029 & 0.034 & 0.057 & 0.011 \\ 
&$2000$ & 0.050 & 0.059 & 0.025 & 0.029 & 0.048 & 0.009 \\
\midrule  
\multirow{4}{*}{$500$}
& $200$ & 0.046 & 0.049 & 0.029 & 0.029 & 0.044 & 0.028 \\ 
& $500$ & 0.044 & 0.051 & 0.026 & 0.032 & 0.039 & 0.024 \\ 
& $1000$ & 0.052 & 0.055 & 0.028 & 0.033 & 0.037 & 0.029 \\ 
& $2000$ & 0.051 & 0.059 & 0.033 & 0.037 & 0.044 & 0.027 \\ 
\midrule 
\multirow{4}{*}{$1000$}
& $200$ & 0.051 & 0.054 & 0.035 & 0.038 & 0.045 & 0.041 \\ 
& $500$ & 0.047 & 0.052 & 0.025 & 0.026 & 0.046 & 0.037 \\ 
& $1000$ & 0.038 & 0.042 & 0.024 & 0.025 & 0.047 & 0.029 \\ 
&$2000$ & 0.043 & 0.049 & 0.027 & 0.031 & 0.041 & 0.029 \\ 
\bottomrule
\end{tabular}}
}
\end{minipage}
\hspace{0.04\linewidth}
\begin{minipage}[t]{0.40\linewidth}
\centering
{\footnotesize 
\makebox[\linewidth][c]{$\rho=0.9$}\par\vspace{2pt}

\resizebox{\linewidth}{!}{%
\begin{tabular}{cccccccc}
\toprule
$n$&$d$& $\varphi_n$ & $\varphi_{n,B}$ & $\varphi_{n,W}$ & $\varphi_{n,W,B}$ & $\varphi_{n,LL}$ & $\varphi_{n,GSS}$ \\
\midrule
\multirow{4}{*}{$100$}
&$200$ & 0.038 & 0.071 & 0.021 & 0.042 & 0.068 & 0.002 \\ 
&$500$ & 0.040 & 0.070 & 0.024 & 0.043 & 0.069 & 0.002 \\ 
&$1000$ & 0.034 & 0.074 & 0.018 & 0.042 & 0.056 & 0.001 \\ 
&$2000$ & 0.052 & 0.096 & 0.029 & 0.047 & 0.068 & 0.000 \\
\midrule 
\multirow{4}{*}{$200$}
&$200$ & 0.038 & 0.069 & 0.025 & 0.043 & 0.049 & 0.013 \\ 
&$500$ & 0.035 & 0.061 & 0.021 & 0.038 & 0.068 & 0.008 \\ 
&$1000$ & 0.036 & 0.061 & 0.023 & 0.041 & 0.043 & 0.005 \\ 
&$2000$ & 0.031 & 0.059 & 0.016 & 0.028 & 0.046 & 0.004 \\
\midrule  
\multirow{4}{*}{$500$}
& $200$ & 0.024 & 0.054 & 0.017 & 0.035 & 0.046 & 0.017 \\ 
& $500$ & 0.029 & 0.059 & 0.018 & 0.037 & 0.050 & 0.018 \\ 
& $1000$ & 0.028 & 0.048 & 0.017 & 0.032 & 0.049 & 0.011 \\ 
& $2000$ & 0.031 & 0.056 & 0.017 & 0.030 & 0.050 & 0.014 \\ 
\midrule 
\multirow{4}{*}{$1000$}
& $200$ & 0.030 & 0.056 & 0.017 & 0.035 & 0.051 & 0.023 \\ 
& $500$ & 0.033 & 0.053 & 0.021 & 0.040 & 0.052 & 0.025 \\ 
& $1000$ & 0.028 & 0.051 & 0.021 & 0.037 & 0.043 & 0.025 \\ 
&$2000$ & 0.022 & 0.041 & 0.014 & 0.025 & 0.042 & 0.017 \\ 
\bottomrule
\end{tabular}}
}
\end{minipage}
\par\vspace{0.5cm}
\begin{minipage}[t]{0.40\linewidth}
\centering
{\footnotesize
\makebox[\linewidth][c]{$\rho=1$}\par\vspace{2pt}

\resizebox{\linewidth}{!}{%
\begin{tabular}{cccccccc}
\toprule
$n$&$d$& $\varphi_n$ & $\varphi_{n,B}$ & $\varphi_{n,W}$ & $\varphi_{n,W,B}$ & $\varphi_{n,LL}$ & $\varphi_{n,GSS}$ \\
\midrule
\multirow{4}{*}{$100$}
& $200$  & 0.001 & 0.056 & 0.001 & 0.036 & 0.056 & 0.000 \\
& $500$  & 0.001 & 0.047 & 0.001 & 0.036 & 0.051 & 0.000 \\
& $1000$ & 0.001 & 0.046 & 0.000 & 0.035 & 0.056 & 0.000 \\
& $2000$ & 0.000 & 0.053 & 0.000 & 0.038 & 0.052 & 0.000 \\
\midrule
\multirow{4}{*}{$200$}
& $200$  & 0.001 & 0.054 & 0.001 & 0.048 & 0.061 & 0.001 \\
& $500$  & 0.001 & 0.057 & 0.000 & 0.048 & 0.060 & 0.000 \\
& $1000$ & 0.001 & 0.057 & 0.000 & 0.037 & 0.048 & 0.000 \\
& $2000$ & 0.000 & 0.049 & 0.000 & 0.035 & 0.056 & 0.000 \\
\midrule
\multirow{4}{*}{$500$}
& $200$  & 0.001 & 0.050 & 0.001 & 0.041 & 0.053 & 0.001 \\
& $500$  & 0.000 & 0.045 & 0.000 & 0.034 & 0.055 & 0.000 \\
& $1000$ & 0.000 & 0.054 & 0.000 & 0.041 & 0.058 & 0.000 \\
& $2000$ & 0.000 & 0.052 & 0.000 & 0.039 & 0.053 & 0.000 \\
\midrule
\multirow{4}{*}{$1000$}
& $200$  & 0.000 & 0.057 & 0.000 & 0.047 & 0.051 & 0.000 \\
& $500$  & 0.000 & 0.051 & 0.000 & 0.043 & 0.050 & 0.000 \\
& $1000$ & 0.000 & 0.053 & 0.000 & 0.043 & 0.050 & 0.000 \\
& $2000$ & 0.000 & 0.044 & 0.000 & 0.035 & 0.052 & 0.000 \\
\bottomrule
\end{tabular}}
}
\end{minipage} 
\caption{\footnotesize  Empirical rejection frequencies under the null absent contamination. Each tabular has~$\Sigma_{n,i,j}=\frac{\nu}{\nu-2}\rho^{|i-j|}$ for~$\nu=4.01$ and~$\rho$ as indicated.}
\label{tab:1}
\end{table}
\end{landscape}

\begin{landscape}
\begin{table}[p]
\centering
\setlength{\tabcolsep}{4pt}
\renewcommand{\arraystretch}{0.88}
\begin{minipage}[t]{0.40\linewidth}
\centering
{\footnotesize 
\makebox[\linewidth][c]{$\rho=0$}\par\vspace{2pt}

\resizebox{\linewidth}{!}{%
\begin{tabular}{cccccccc}
\toprule
$n$&$d$& $\varphi_n$ & $\varphi_{n,B}$ & $\varphi_{n,W}$ & $\varphi_{n,W,B}$ & $\varphi_{n,LL}$ & $\varphi_{n,GSS}$ \\ 
\midrule
\multirow{4}{*}{$100$}
&$200$ & 0.000 & 0.009 & 0.080 & 0.105 & 0.076 & 0.000 \\ 
&$500$ & 0.000 & 0.007 & 0.080 & 0.113 & 0.092 & 0.000 \\ 
&$1000$ & 0.000 & 0.011 & 0.075 & 0.119 & 0.100 & 0.000 \\ 
&$2000$ & 0.000 & 0.007 & 0.073 & 0.123 & 0.090 & 0.000 \\
\midrule 
\multirow{4}{*}{$200$}
&$200$ & 0.000 & 0.005 & 0.051 & 0.059 & 0.101 & 0.000 \\ 
&$500$ & 0.000 & 0.005 & 0.039 & 0.051 & 0.103 & 0.000 \\ 
&$1000$ & 0.000 & 0.005 & 0.043 & 0.055 & 0.090 & 0.000 \\ 
&$2000$ & 0.000 & 0.006 & 0.042 & 0.058 & 0.107 & 0.000 \\
\midrule  
\multirow{4}{*}{$500$}
& $200$ & 0.000 & 0.003 & 0.051 & 0.057 & 0.129 & 0.000 \\ 
& $500$ & 0.000 & 0.003 & 0.041 & 0.046 & 0.147 & 0.000 \\ 
& $1000$ & 0.000 & 0.003 & 0.029 & 0.036 & 0.124 & 0.000 \\ 
& $2000$ & 0.000 & 0.002 & 0.033 & 0.041 & 0.123 & 0.000 \\ 
\midrule 
\multirow{4}{*}{$1000$}
& $200$ & 0.000 & 0.004 & 0.043 & 0.048 & 0.122 & 0.000 \\ 
& $500$ & 0.000 & 0.002 & 0.030 & 0.030 & 0.122 & 0.000 \\ 
& $1000$ & 0.000 & 0.005 & 0.033 & 0.041 & 0.129 & 0.000 \\ 
&$2000$ & 0.000 & 0.003 & 0.035 & 0.041 & 0.121 & 0.000 \\ 
\bottomrule
\end{tabular}}
}
\end{minipage}
\hspace{0.04\linewidth}
\begin{minipage}[t]{0.40\linewidth}
\centering
{\footnotesize 
\makebox[\linewidth][c]{$\rho=0.9$}\par\vspace{2pt}

\resizebox{\linewidth}{!}{%
\begin{tabular}{cccccccc}
\toprule
$n$&$d$& $\varphi_n$ & $\varphi_{n,B}$ & $\varphi_{n,W}$ & $\varphi_{n,W,B}$ & $\varphi_{n,LL}$ & $\varphi_{n,GSS}$ \\ 
\midrule
\multirow{4}{*}{$100$}
&$200$ & 0.000 & 0.009 & 0.048 & 0.105 & 0.090 & 0.000 \\ 
&$500$ & 0.000 & 0.007 & 0.047 & 0.120 & 0.083 & 0.000 \\ 
&$1000$ & 0.000 & 0.011 & 0.049 & 0.127 & 0.098 & 0.000 \\ 
&$2000$ & 0.000 & 0.009 & 0.069 & 0.149 & 0.097 & 0.000 \\
\midrule 
\multirow{4}{*}{$200$}
&$200$ & 0.000 & 0.004 & 0.028 & 0.069 & 0.117 & 0.000 \\ 
&$500$ & 0.000 & 0.006 & 0.030 & 0.074 & 0.100 & 0.000 \\ 
&$1000$ & 0.000 & 0.006 & 0.030 & 0.065 & 0.084 & 0.000 \\ 
&$2000$ & 0.000 & 0.005 & 0.028 & 0.058 & 0.093 & 0.000 \\
\midrule  
\multirow{4}{*}{$500$}
& $200$ & 0.000 & 0.004 & 0.021 & 0.050 & 0.126 & 0.000 \\ 
& $500$ & 0.000 & 0.006 & 0.024 & 0.054 & 0.131 & 0.000 \\ 
& $1000$ & 0.000 & 0.004 & 0.024 & 0.041 & 0.118 & 0.000 \\ 
& $2000$ & 0.000 & 0.001 & 0.018 & 0.036 & 0.118 & 0.000 \\ 
\midrule 
\multirow{4}{*}{$1000$}
& $200$ & 0.000 & 0.004 & 0.028 & 0.049 & 0.133 & 0.000 \\ 
& $500$ & 0.000 & 0.004 & 0.024 & 0.042 & 0.139 & 0.000 \\ 
& $1000$ & 0.000 & 0.004 & 0.026 & 0.048 & 0.139 & 0.000 \\ 
&$2000$ & 0.000 & 0.003 & 0.025 & 0.044 & 0.122 & 0.000 \\ 
\bottomrule
\end{tabular}}
}
\end{minipage}
\par\vspace{0.5cm}
\begin{minipage}[t]{0.40\linewidth}
\centering
{\footnotesize
\makebox[\linewidth][c]{$\rho=1$}\par\vspace{2pt}

\resizebox{\linewidth}{!}{%
\begin{tabular}{cccccccc}
\toprule
$n$&$d$& $\varphi_n$ & $\varphi_{n,B}$ & $\varphi_{n,W}$ & $\varphi_{n,W,B}$ & $\varphi_{n,LL}$ & $\varphi_{n,GSS}$ \\
\midrule
\multirow{4}{*}{$100$}
& $200$  & 0.000 & 0.011 & 0.000 & 0.053 & 0.144 & 0.000 \\
& $500$  & 0.000 & 0.008 & 0.000 & 0.058 & 0.150 & 0.000 \\
& $1000$ & 0.000 & 0.009 & 0.001 & 0.043 & 0.129 & 0.000 \\
& $2000$ & 0.000 & 0.009 & 0.000 & 0.045 & 0.139 & 0.000 \\
\midrule
\multirow{4}{*}{$200$}
& $200$  & 0.000 & 0.005 & 0.001 & 0.045 & 0.157 & 0.000 \\
& $500$  & 0.000 & 0.004 & 0.000 & 0.050 & 0.167 & 0.000 \\
& $1000$ & 0.000 & 0.009 & 0.000 & 0.042 & 0.163 & 0.000 \\
& $2000$ & 0.000 & 0.004 & 0.000 & 0.039 & 0.160 & 0.000 \\
\midrule
\multirow{4}{*}{$500$}
& $200$  & 0.000 & 0.005 & 0.001 & 0.051 & 0.165 & 0.000 \\
& $500$  & 0.000 & 0.005 & 0.000 & 0.049 & 0.163 & 0.000 \\
& $1000$ & 0.000 & 0.008 & 0.000 & 0.043 & 0.144 & 0.000 \\
& $2000$ & 0.000 & 0.006 & 0.000 & 0.041 & 0.159 & 0.000 \\
\midrule
\multirow{4}{*}{$1000$}
& $200$  & 0.000 & 0.011 & 0.000 & 0.042 & 0.165 & 0.000 \\
& $500$  & 0.000 & 0.018 & 0.000 & 0.048 & 0.170 & 0.000 \\
& $1000$ & 0.000 & 0.013 & 0.000 & 0.050 & 0.174 & 0.000 \\
& $2000$ & 0.000 & 0.011 & 0.000 & 0.048 & 0.156 & 0.000 \\
\bottomrule
\end{tabular}}
}
\end{minipage}

\caption{\footnotesize Empirical rejection frequencies under the null in the presence of contamination. The adversary replaces~$q=3$ randomly chosen~$X_{i,n}$ by~$\tilde{X}_{i,n}=(500,\hdots,500)'$ and leaves the remaining observations unchanged. Each tabular has~$\Sigma_{n,i,j}=\frac{\nu}{\nu-2}\rho^{|i-j|}$ for~$\nu=4.01$ and~$\rho$ as indicated.}
\label{tab:2}
\end{table}
\end{landscape}

\subsection{Power}
\begin{figure}[H]
\centering

\begin{minipage}{6cm}
\centering
\footnotesize No Contamination, $\rho=0$\\[-15pt]
\includegraphics[width=\linewidth]{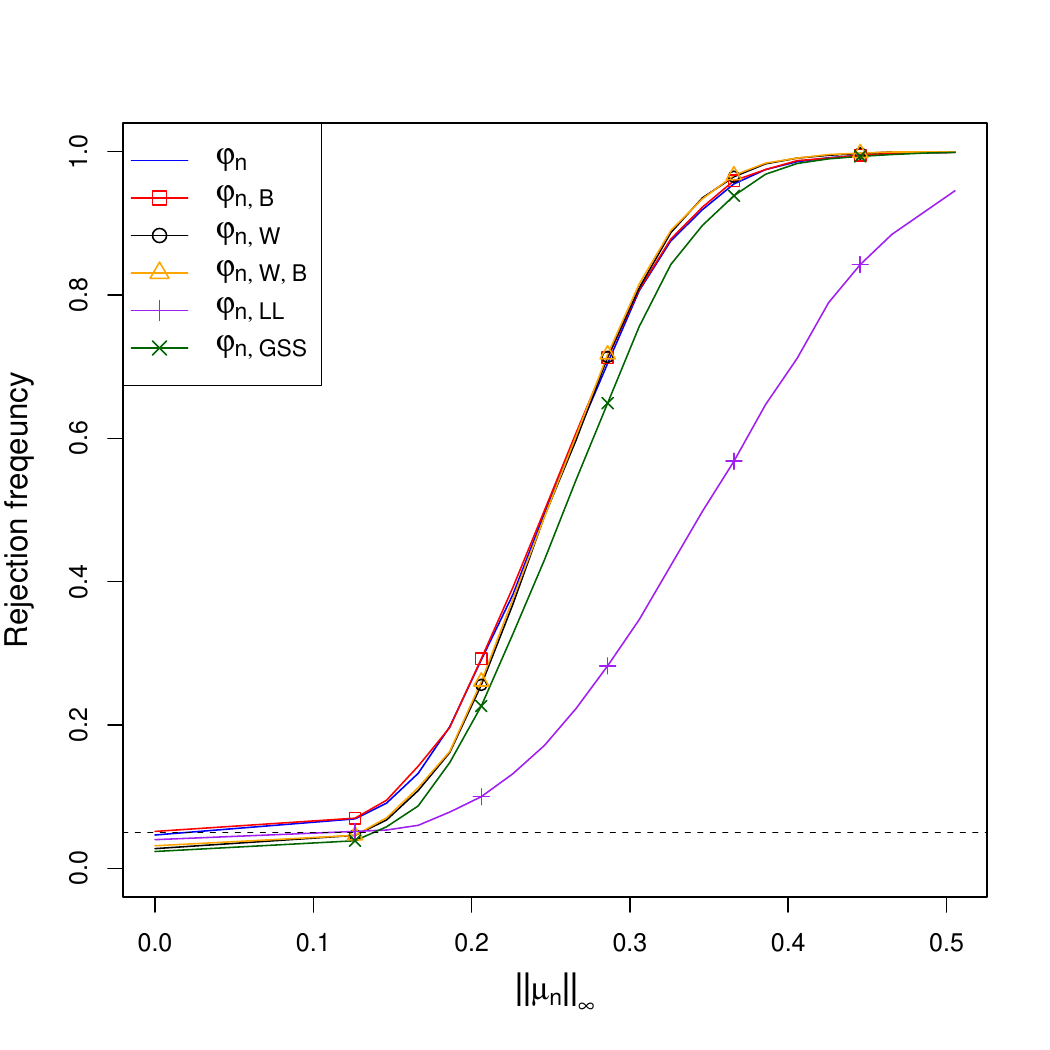}
\end{minipage}
\hspace{0.8cm}
\begin{minipage}{6cm}
\centering
\footnotesize Contamination, $\rho=0$\\[-15pt]
\includegraphics[width=\linewidth]{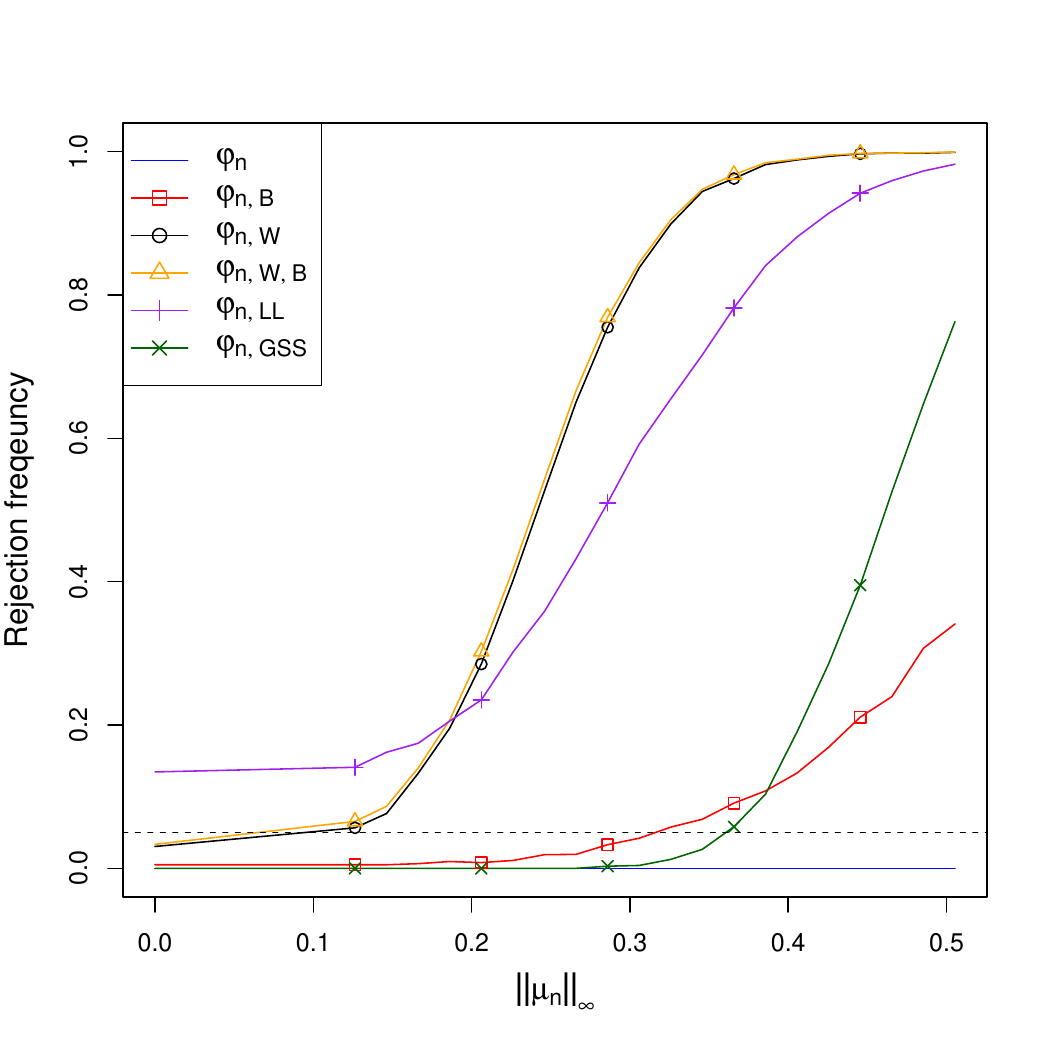}
\end{minipage}

\vspace{0.3cm}

\begin{minipage}{6cm}
\centering
\footnotesize No Contamination, $\rho=0.9$\\[-15pt]
\includegraphics[width=\linewidth]{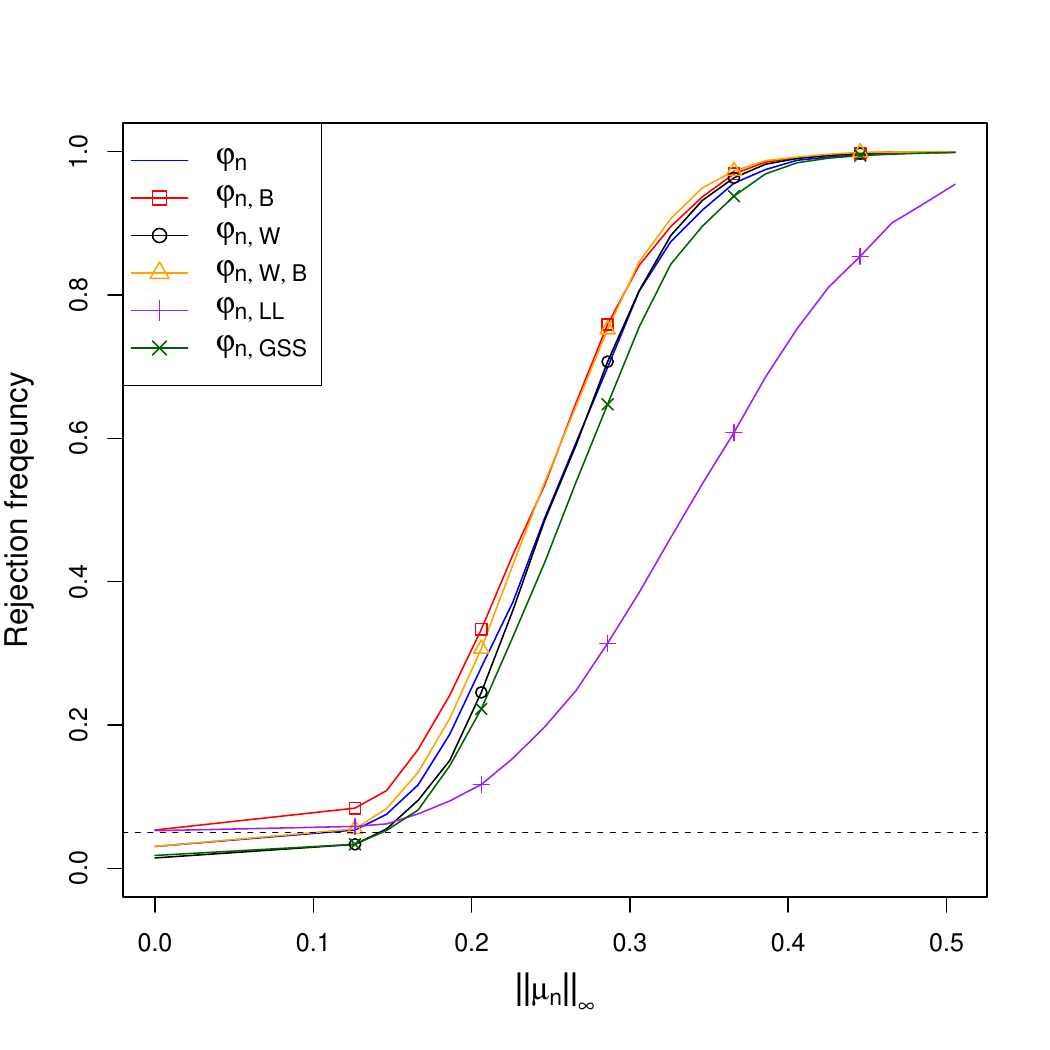}
\end{minipage}
\hspace{0.8cm}
\begin{minipage}{6cm}
\centering
\footnotesize Contamination, $\rho=0.9$\\[-15pt]
\includegraphics[width=\linewidth]{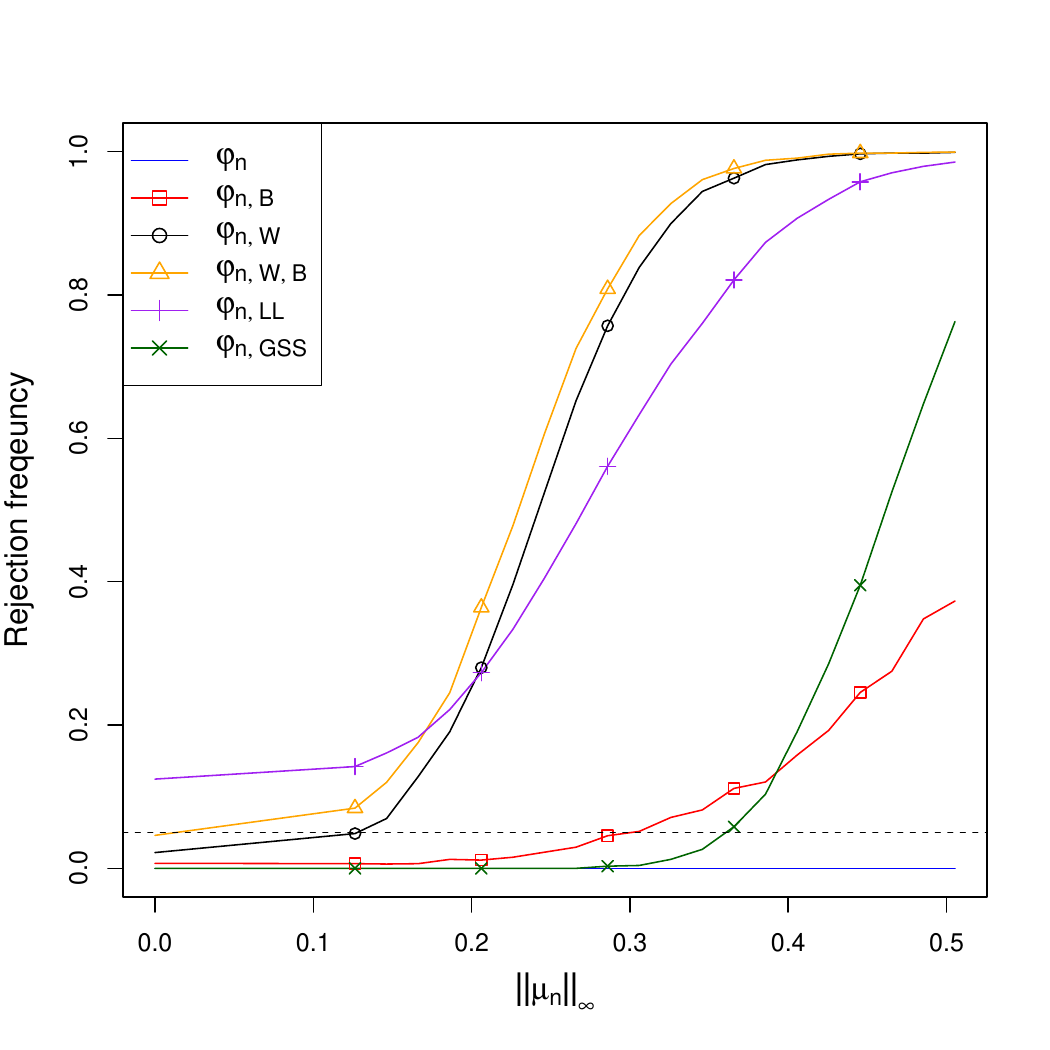}
\end{minipage}

\vspace{0.3cm}

\begin{minipage}{6cm}
\centering
\footnotesize No Contamination, $\rho=1$\\[-15pt]
\includegraphics[width=\linewidth]{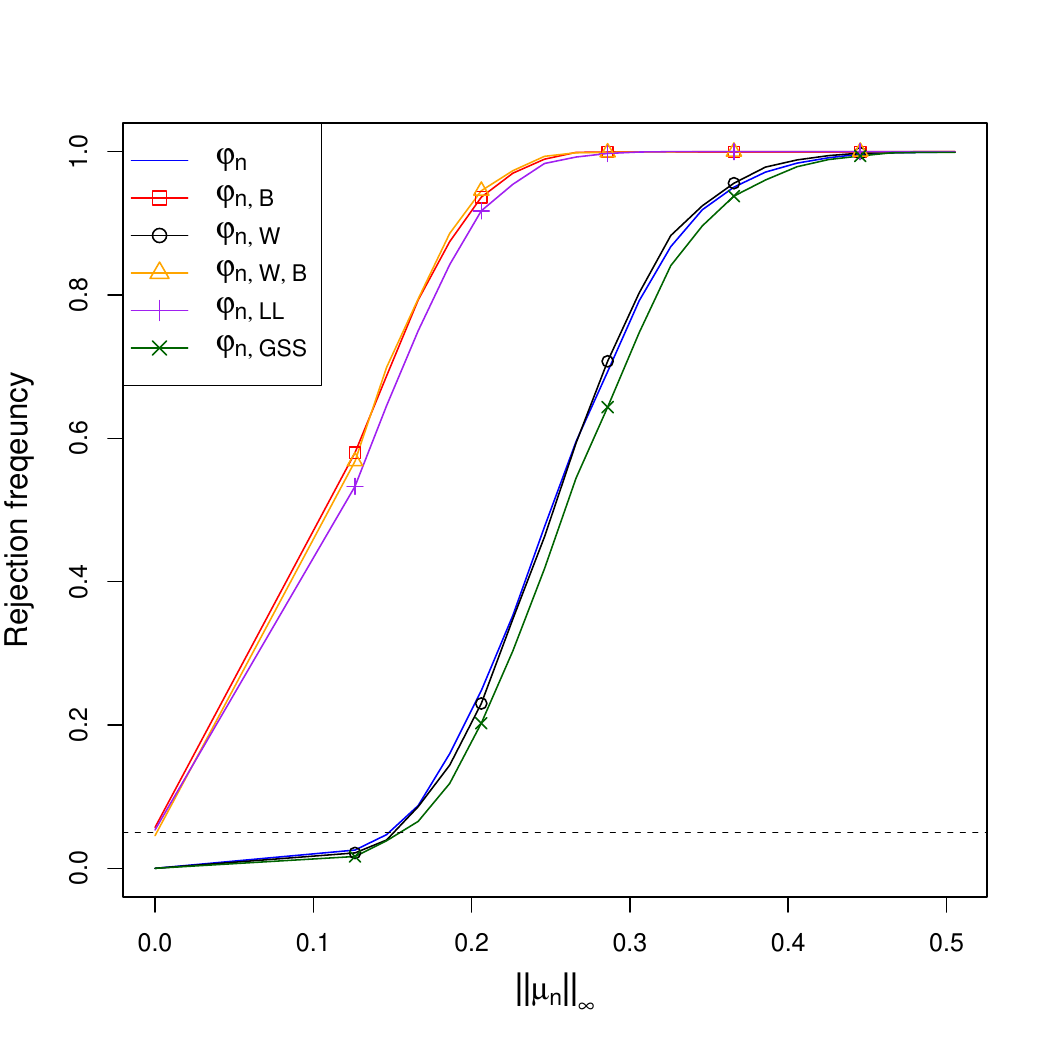}
\end{minipage}
\hspace{0.8cm}
\begin{minipage}{6cm}
\centering
\footnotesize Contamination, $\rho=1$\\[-15pt]
\includegraphics[width=\linewidth]{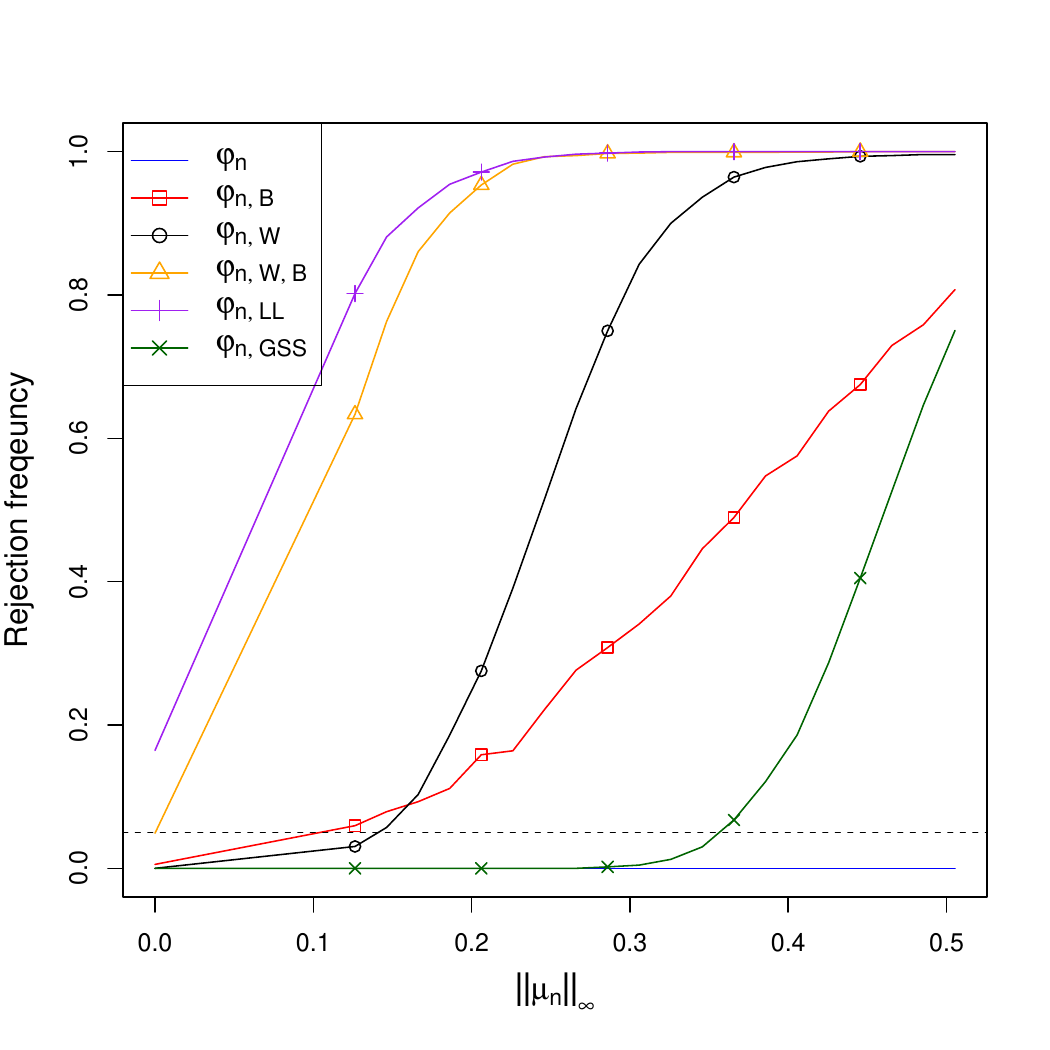}
\end{minipage}
\caption{\footnotesize Rejection frequencies as a function of~$||\mu_n||_\infty$.~$n=500$ and~$d=1{,}000$. The left three plots consider the case of no contamination and the right three plots the case of contamination. Each plot has~$\Sigma_{n,i,j}=\frac{\nu}{\nu-2}\rho^{|i-j|}$ for~$\nu=4.01$ and~$\rho$ as indicated.}
\label{fig:1}
\end{figure}

We study the power of the tests in the setting of~$n=500$ and~$d=1{,}000$ and consider a sparse deviation~$\mu_n=(a,0,\hdots,0)'$ from $\mathsf{H}_{0,n}:\mu_n=0_d$ for a range of~$a\geq 0$. Figure~\ref{fig:1} plots the power of the tests as a function of~$||\mu_n||_\infty=a$. 

The left column of Figure~\ref{fig:1} considers the case without contamination ($q=0$). Since the conditions of Part 1.~of Theorem~\ref{thm:comparison} are satisfied for~$\rho\in\cbr[0]{0,0.9}$, it is not surprising that the power gains of~$\varphi_{n,W,B}$ over~$\varphi_{n,W}$ are small in these cases.\footnote{Theorem~\ref{thm:comparison} compares the power of~$\varphi_{n,W}$ and $\varphi_{n,W,B}$ but \emph{absent} contamination Figure~\ref{fig:1} indicates that a similar comparison holds between~$\varphi_{n}$ and~$\varphi_{n,B}$, cf.~also Remark~\ref{rem:sumspowercomp}.} For these values of~$\rho$ the test~$\varphi_{n,LL}$ has the lowest power. For~$\rho=1$, Part 3.~of Theorem~\ref{thm:comparison} helps to explain why the power of the bootstrap based test~$\varphi_{n,W,B}$ can be much larger than that of~$\varphi_{n,W}$ for small to moderate values of~$||\mu_n||_\infty$. Part 2.~of that theorem explains why the power differences vanish for sufficiently large~$||\mu_n||_\infty$ since then both~$\varphi_{n,W}$ and~$\varphi_{n,W,B}$ are consistent. It is noteworthy that in this extreme correlation setting ($\rho = 1$), the three bootstrap-based tests taking into account this dependence ($\varphi_{n,B}$,~$\varphi_{n,W,B}$ and~$\varphi_{n,LL}$) have much higher power than the three tests using deterministic ``worst-case'' critical values ($\varphi_{n}$,~$\varphi_{n,W}$ and~$\varphi_{n,GSS}$).

The right column of Figure~\ref{fig:1} considers the case of contaminated data ($q=3$). The power functions of the tests~$\varphi_{n,W}$ and~$\varphi_{n,W,B}$, which are based on the robust winsorized means, are hardly affected by the contamination (compare the left and right panels in each row). Thus, Parts 1.--3. of Theorem~\ref{thm:comparison} again explain the relative power of $\varphi_{n,W}$ and~$\varphi_{n,W,B}$ across all values of~$\rho$. The remaining four tests~$\varphi_n$,~$\varphi_{n,B}$,~$\varphi_{n,LL}$, and~$\varphi_{n,GSS}$ are not robust to contamination and the power curves should be interpreted with much caution. For example,~$\varphi_{n,LL}$ exhibits power similar to~$\varphi_{n,W,B}$, but also has much larger null rejection probabilities, cf.~Table~\ref{tab:2} with~$n=500$ and~$d=1{,}000$.

\subsection{Further numerical results} Section \ref{sec:RobCheck} in the supplementary appendix provides further numerical checks, which have already briefly been discussed in the present section. In particular, it investigates the effect of misspecifying~$q$ in the implementation of~$\varphi_{n,W}$ and~$\varphi_{n,W,B}$ and the sensitivity to the value of~$\lambda_2'$. The robustness checks also study the case of~$\nu=3.01$ instead of~$\nu=4.01$, corresponding to a more heavy tailed setting.

\section{Conclusion}

We studied quantile-winsorized max-tests and established their uniform validity and power properties in regimes where~$\log(d)=o\del[1]{n^{(m-2)/(5m-2)}}$  and only~$m > 2$ moments exist. These tests are robust to adversarial contamination, which standard max-tests are not, without sacrificing asymptotic power, i.e., their robustness comes ``for free'' at first order. The condition~$\sqrt{n\log(d)}\overline{\eta}_n^{1-1/m}\to 0$ imposed on~$\overline{\eta}_n$ is proven to be necessary for this. Whether one bases winsorized max-tests on fixed or bootstrap critical values does not change their asymptotic properties as long as the correlation between coordinates decays sufficiently quickly. As that decay is unknown in practice and since power gains are possible under slow decay of correlations, pairing our winsorized max-tests with bootstrap critical values is our recommended choice for practical use.

\bibliographystyle{ecta} 
\bibliography{ref}		

\newpage
\begin{appendix}
\numberwithin{equation}{section}	

\section{Proof of Theorem~\ref{thm:sums}}
\subsection{Auxiliary results concerning variance estimation}
The following lemma is probably well-known, but we were unable to pinpoint an exact reference. Similar results can be found in, e.g., Lemma~D.5 in~\cite{chernozhukov2019inference}.
\begin{lemma}\label{lem:variance}
Let~$b>0$ and~$m > 4$. Then, there exists a (positive) constant~$C$, say, depending only on~$b$ and~$m$, such that the following holds: for every~$n \in \N$, every~$z>1$, and for~$Y_1,\hdots,Y_n$ i.i.d.~random variables with mean~$\mu$, variance~$\sigma^2$, and such that~$\E|Y_1-\mu|^m\leq b^m$, we have
\begin{equation}\label{eqn:fnbd}
\P\del[2]{\envert[0]{s_n^2-\sigma^2}
\geq
C\sbr[1]{n^{-1/2}\sqrt{\log (z)}+n^{2/m-1}z^{2/m}+n^{-1}\log (z)}}
\leq
\frac{8}{z},
\end{equation} 
where~$s_n^2:=n^{-1}\sum_{i=1}^n\del[1]{Y_i-\overline{Y}_n}^2$ and~$\overline{Y}_n:=n^{-1}\sum_{i=1}^nY_i$.
\end{lemma}
\begin{proof}
Fix~$b, m, n, z$ and~$Y_1,\hdots,Y_n$ as in the present lemma's statement, and write
\begin{equation}\label{eq:varest1}
s_{n}^2-\sigma^2
=
\frac{1}{n}\sum_{i=1}^n\sbr[2]{\del[1]{Y_{i}-\mu}^2-\sigma^2}-\del[1]{\overline{Y}_{n}-\mu}^2.
\end{equation}
The~$(Y_{i}-\mu)^2-\sigma^2$ are i.i.d.~with expectation zero, and~$\E|\del[1]{Y_1-\mu}^2-\sigma^2|^{m/2}
\leq
2^{m/2}b^m$.
By the Fuk-Nagaev inequality, as stated in Equation~(1.8) of~\cite{rio2017constants} and applied with
\begin{itemize}
\item $X_i$ there being~$(Y_i-\mu)^2-\sigma^2$,
\item $\sigma^2$ there being~$n\E\sbr[1]{\del[1]{Y_{1}-\mu}^2-\sigma^2}^2\leq C_1 n$ for~$C_1$ a constant depending only on~$b$,
\item $q$ there being~$m/2 > 2$,
\item $C_q(X)$ there being upper bounded by~$	\left(n\E\envert[1]{\del[1]{Y_1-\mu}^2-\sigma^2}^{m/2}\right)^{2/m}
\leq
\del[1]{n2^{m/2}b^m}^{2/m}$,
\end{itemize}
it follows that, for a constant~$C_2$ depending only on~$b$ and~$m$, we have
\begin{equation*}
\P\del[3]{\frac{1}{n}\sum_{i=1}^n\sbr[1]{\del[1]{Y_{i}-\mu}^2-\sigma^2}\geq C_2\del[1]{n^{-1/2}\sqrt{\log z}+n^{2/m-1}z^{2/m}}}\leq \frac{2}{z}.
\end{equation*}
Combining this bound with an identical argument for~$-\frac{1}{n}\sum_{i=1}^n\sbr[1]{\del[1]{Y_{i}-\mu}^2-\sigma^2}$ yields
\begin{equation}\label{eq:varest2}
\P\del[3]{\envert[2]{\frac{1}{n}\sum_{i=1}^n\sbr[1]{\del[1]{Y_{i}-\mu}^2-\sigma^2}}\geq C_2\del[1]{n^{-1/2}\sqrt{\log z}+n^{2/m-1}z^{2/m}}}\leq \frac{4}{z}.
\end{equation}
Writing~$\overline{Y}_n-\mu
=
\frac{1}{n}\sum_{i=1}^n(Y_i-\mu)$,
the Fuk-Nagaev inequality as in Equation~(1.8) in~\cite{rio2017constants} with~$X_i$ there being~$Y_i-\mu$,~$q$ there being~$m$, and similar arguments to above, yields
\begin{equation}\label{eq:varest3}
\P\del[2]{\envert[0]{\overline{Y}_n-\mu}\geq C_3\del[1]{n^{-1/2}\sqrt{\log z}+n^{1/m-1}z^{1/m}}}
\leq \frac{4}{z},
\end{equation}
for a constant~$C_3$ depending only on~$b$ and~$m$. Thus, for a constant~$C_4$ depending only on~$b$ and~$m$, it holds that (we also used that~$n^{2/m-2} \leq n^{2/m-1}$)
\begin{equation}\label{eq:varest3*}
\P\del[2]{\envert[0]{\overline{Y}_n-\mu}^2\geq C_4\del[1]{n^{-1}\log z+n^{2/m-1}z^{2/m}}}
\leq \frac{4}{z}.
\end{equation}
Combining the bounds in~\eqref{eq:varest2} and~\eqref{eq:varest3*} with~\eqref{eq:varest1} yields~\eqref{eqn:fnbd}.
\end{proof}

Recall that we denote~$s_{n,j}^2=n^{-1}\sum_{i=1}^n\del[1]{X_{i,n,j}-\overline{X}_{n,j}}^2$, where~$\overline{X}_{n,j}=n^{-1}\sum_{i=1}^nX_{i,n,j}$. 

\begin{corollary}\label{cor:varestim}
Let~$b_1\in(0,\infty)$,~$b_2\in(b_1,\infty)$, and~$m > 4$. Then, there exists a (positive) constant~$C$, say, depending only on~$b_2$ and~$m$, such that for every~$n\geq 3$ and every~$P\in \mc{P}(b_1,b_2,m)$, we have
\begin{equation}\label{eqn:varestim*}
P\del[4]{\max_{j=1,\hdots,d}\envert[1]{s_{n,j}^2-\sigma_{2,n,j}^2(P)}> C \mathfrak{a}(n, d, m)}
\leq \frac{8}{\log(n)},
\end{equation}
where we abbreviated
\begin{equation}\label{eqn:fadef}
\mathfrak{a}(n, d, m) :=  \frac{\log(dn)}{n}+\sqrt{\frac{\log(dn)}{n}}+\frac{[d\log(n)]^{2/m}}{n^{1-2/m}};
\end{equation}
if~$\mathfrak{a}(n, d, m) <\frac{b_1^2}{2C}$, then, setting~$C' := \sqrt{2}C/b_1^3$, it furthermore holds that
\begin{equation}\label{eqn:varestim**}
P\left(\min_{j = 1, \hdots, d} s_{n,j} \geq \frac{b_1}{\sqrt{2}} \text{ and } \max_{j=1,\hdots,d}\envert[1]{s_{n,j}^{-1}-\sigma_{2,n,j}^{-1}(P)} \leq C'\mathfrak{a}(n, d, m)   \right)
\geq 1- \frac{8}{\log(n)}.
\end{equation}
\end{corollary}

\begin{proof}
To obtain~\eqref{eqn:varestim*} we use a union bound over~$j=1,\hdots,d$, apply Lemma~\ref{lem:variance} with~$z=d\log(n) \geq \log(3) > 1$, and use~$\log(d\log(n))\leq \log(dn)$. 

If~$\mathfrak{a}(n, d, m) <\frac{b_1^2}{2C}$, then~\eqref{eqn:varestim*} shows that with probability at least~$1-8/\log(n)$
\begin{equation}\label{eqn:ints^2}
\max_{j=1,\hdots,d}\envert[1]{s_{n,j}^2-\sigma_{2,n,j}^2(P)} \leq C \mathfrak{a}(n, d, m) < \frac{b_1^2}{2}.
\end{equation}
From~\eqref{eqn:ints^2}, together with~$\min_{j = 1, \hdots, d} \sigma_{2,n,j}(P) \geq b_1 > 0$ as a consequence of~$P\in \mc{P}(b_1,b_2,m)$, we conclude that~$s_{n,j} \geq b_1/\sqrt{2} > 0$ for every~$j = 1, \hdots, d$, and
\begin{equation*}
\max_{j=1,\hdots,d}\envert[1]{s_{n,j}^{-1}-\sigma_{2,n,j}^{-1}(P)} \leq \max_{j=1,\hdots,d}\frac{|s_{n,j}-\sigma_{2,n,j}(P)|}{s_{n,j}\sigma_{2,n,j}(P)}
\leq 
\frac{\sqrt{2}}{b_1^2} \max_{j=1,\hdots,d}\envert[1]{s_{n,j}-\sigma_{2,n,j}(P)}.
\end{equation*}
Note that, because~$1 \leq \sigma_{2,n,j}(P)/b_1$ and~$s_{n,j} \geq 0$,
\begin{equation*}
|s_{n,j}-\sigma_{2,n,j}(P)| \leq |s_{n,j}-\sigma_{2,n,j}(P)| |s_{n,j}+\sigma_{2,n,j}(P)|/b_1 = |s_{n,j}^2-\sigma^2_{2,n,j}(P)|/b_1.
\end{equation*}
Therefore, whenever~\eqref{eqn:ints^2} holds,~$\min_{j =1, \hdots, d}s_{n,j} \geq b_1/\sqrt{2} > 0$, and $$\max_{j=1,\hdots,d}\envert[1]{s_{n,j}^{-1}-\sigma_{2,n,j}^{-1}(P)} \leq \frac{\sqrt{2}}{b_1^3} \max_{j=1,\hdots,d}|s_{n,j}^2-\sigma^2_{2,n,j}(P)| \leq C' \mathfrak{a}(n, d, m).$$ Because~\eqref{eqn:ints^2} holds with probability at least~$1-8/\log(n)$ this establishes~\eqref{eqn:varestim**}.
\end{proof}

\subsection{Gaussian approximation for centered standardized sums}\label{sec:GA_ST_sums}
Denote by~$\mc{H}=\mc{H}_n$ the class of (generalized) hyperrectangles in~$\R^d$, that is the class of all sets of the form
\begin{equation*}
H=\cbr[1]{x\in\R^d:a_j\leq x_j\leq b_j\text{ for all }j=1,\hdots, d},
\end{equation*}
where~$-\infty\leq a_j\leq b_j\leq \infty$ for all~$j = 1, \hdots, d$. 
The following lemma, a Gaussian approximation theorem for centered and unstandardized averages, is a special case of Theorem~2.5 in~\cite{chernozhuokov2022improved} (together with the discussion immediately after Theorem~2.1 in the same reference); note that the latter theorem imposes the existence of a fourth moment via their ``Condition M.''
\begin{lemma}\label{lem:cherno}
Let~$b_1\in(0,\infty)$,~$b_2\in(b_1,\infty)$, and~$m \geq 4$. Then, there exists a (positive) constant~$C$, say, depending only on~$b_1$,~$b_2$ and~$m$, such that for every~$P\in\mc{P}(b_1,b_2,m)$, if~$n \geq 3$ and~$d \geq 2$, it holds that
\begin{align*}
&\sup_{H\in\mc{H}}\envert[3]{P\del[2]{\frac{1}{\sqrt{n}}\sum_{i=1}^n[X_{i,n}-\mu_n(P)]\in H}-P\del[1]{\Sigma_n^{1/2}(P)Z_d\in H}}\\
& \leq 
C\del[4]{\sbr[3]{\frac{d^{2/m}\log^{5}(d)}{n}}^{1/4}+\sbr[3]{\frac{d^{2/m}\log^{3-2/m}(d)}{n^{1-2/m}}}^{1/2}} =: C \cdot \mathfrak{b}(n, d, m).
\end{align*}
\end{lemma}
\begin{proof}
We set up for an application of Theorem 2.5 in \cite{chernozhuokov2022improved} with~$q$ there being our~$m$, and~$X_i$ there being~$[X_{i,n}-\mu_n(P)]$. Note that 
\begin{equation*}
E_P\enVert[0]{X_{1,n}-\mu_n(P)}_\infty^m
\leq
\del[3]{\sum_{j=1}^dE_P\envert[1]{X_{1,n,j}-\mu_{n,j}(P)}^m}
\leq
b_2^md,
\end{equation*}
such that~$B_n$ in Theorem 2.5 of~\cite{chernozhuokov2022improved} can be chosen as~$\max(1,b_2)d^{1/m}$. ``Condition M'' of~\cite{chernozhuokov2022improved} is also satisfied, because for all~$j=1,\hdots,d$
\begin{equation*}
E_P\del[1]{X_{i,n,j}-\mu_{n,j}(P)}^4
\leq
b_2^4
\leq
b_2^2B_n^2.
\end{equation*}
Now, the claimed bound follows from Theorem 2.5 of~\cite{chernozhuokov2022improved} (in which~$n \geq 3$ and~$d \geq 2$ are maintained throughout, cf.~the beginning of their Section~2.2); the discussion immediately after Theorem~2.1 in the same reference explains how to convert the bound in their Theorem 2.5 into one valid uniformly over~$\mc{H}$.  
\end{proof}

The following auxiliary result is a consequence of the Gaussian anti-concentration inequality stated in Theorem 1 of~\cite{chernozhukov2017detailed}, cf.~also Lemma~A.1 in~\cite{chernozhukov2017central}. For any non-empty set~$A\subseteq \R^d$ and~$\zeta >0$, let $$A^{\zeta,\infty}=\cbr[1]{x\in\R^d:\inf_{y\in A}\|x-y\|_\infty\leq \zeta},$$ and $$A^{-\zeta,\infty}=\cbr[1]{x\in\R^d:\mc{B}_\infty(x,\zeta) \subseteq A} \quad \text{ where } \quad \mc{B}_{\infty}(x,\zeta)=\cbr[1]{y\in\R^d:\|y-x\|_{\infty} \leq \zeta}.$$
For any non-empty set~$A\subseteq \R^d$ we set~$A^{0,\infty}=A$. We also set~$\emptyset^{\zeta,\infty}=\emptyset^{-\zeta,\infty}=\emptyset$ for~$\zeta\geq 0$.
\begin{lemma}\label{lem:auxlem}
Let~$U$ and~$V$ be~$d$-variate random vectors, and let~$\Xi$ be a~$d\times d$-dimensional covariance matrix having the square root of its diagonal elements bounded from below by~$\xi > 0$. Let~$a\geq0$, $b\geq 0$, $c\geq 0$, and~$H\in\mc{H}$ be such that
\begin{equation}\label{eqn:auxGapp1}
\max_{v \in \{-1, 1\}}\left| \P(U \in H^{vb, \infty}) - \P(\Xi^{1/2} Z_d \in H^{vb, \infty}) \right| \leq a ~~ \text{and} ~~ \P\left(\|U - V\|_{\infty} > b\right) \leq c.
\end{equation}
Then,
\begin{equation}\label{eqn:auxGapp2}
\left| \P(V \in H) - \P(\Xi^{1/2} Z_d \in H) \right| \leq a + \left(\sqrt{2 \log(d)} + 4\right)b \xi^{-1} + c.
\end{equation}
\end{lemma}
\begin{proof}
If~$H$ is empty,~\eqref{eqn:auxGapp2} trivially holds. Hence, we assume that~$H$ is non-empty. The Gaussian anti-concentration inequality in Theorem~1 of~\cite{chernozhukov2017detailed}, cf.~its ``two-sided'' version in Lemma~E.1 of \cite{kock2025high}, which we re-state in Lemma~\ref{lem:anticonc} below for convenience, delivers 
\begin{equation*}
\max_{v \in \{-1, 1\}} \left| \P\left(\Xi^{1/2} Z_d \in H^{v b, \infty}\right) - \P\left(\Xi^{1/2} Z_d \in H\right) \right| \leq \left(\sqrt{2 \log(d)} + 4\right) b \xi^{-1}.
\end{equation*}
In combination with the first condition in~\eqref{eqn:auxGapp1} we obtain 
\begin{equation}\label{eqn:gace2}
\max_{v \in \{-1, 1\}} \left| \P\left(U \in H^{v b, \infty}\right) - \P\left(\Xi^{1/2} Z_d \in H\right) \right| \leq a + \left(\sqrt{2 \log(d)} + 4\right) b \xi^{-1}.
\end{equation}
Furthermore, 
\begin{equation*}
\begin{aligned}
\{V \in H\} &\subseteq \{U \in H^{b, \infty}\} &&\cup \{\|U - V\|_{\infty} > b\}, \text{ and }\\
\{U \in H^{-b, \infty}\} &\subseteq \{V \in H\} &&\cup \{\|U - V\|_{\infty} > b\},
\end{aligned}
\end{equation*}
so that the second condition in~\eqref{eqn:auxGapp1} yields
\begin{equation*}
- c + \P(U \in H^{-b, \infty}) \leq \P(V \in H) \leq c + \P(U \in H^{b, \infty}),
\end{equation*}
which, combined with~\eqref{eqn:gace2}, delivers~\eqref{eqn:auxGapp2}.
\end{proof}

\begin{lemma}[Lemma~E.1 of \cite{kock2025high}]\label{lem:anticonc}
Let~$Z$ in~$\R^d$ with~$d\geq 1$ be such that~$Z\sim~\mathsf{N}_d(0,\Sigma)$ with~$\Sigma_{j,j} \geq \underline{\sigma}^2$ for all~$j=1,\hdots, d$ and some~$\underline{\sigma}^2>0$. Then, for all real numbers~$\delta_1$ and~$\delta_2$ and all~$a$ and~$b$ in~$\overline{\R}^d$, it holds that
\begin{equation}\label{eq:twogauss1}
\P\del[1]{Z\in[a+\delta_1,b+\delta_2]}
\leq
\P\del[1]{Z \in[a,b]}+\frac{\overline{\delta}}{\underline{\sigma}}\del[1]{\sqrt{2\log(d)}+4},
\end{equation}
where~$\overline{\delta}=|\delta_1|\vee |\delta_2|$, and
\begin{equation}\label{eq:twogauss2}
\P\del[1]{Z\in[a+\delta_1,b+\delta_2]}
\geq
\P\del[1]{Z \in[a,b]}-\frac{\overline{\delta}}{\underline{\sigma}}\del[1]{\sqrt{2\log(d)}+4}.
\end{equation}
\end{lemma}

For ease of reference, we recall the following instance of the Khatri-{\v{S}}id{\'a}k inequality (\cite{khatri1967certain} and \cite{vsidak1967rectangular}), e.g., from Corollary 2.4.6 in~\cite{gine2016mathematical}.

\begin{theorem}[Khatri-{\v{S}}id{\'a}k inequality]\label{thm:ksi}
Let the mean zero random variables~$z_1, \hdots, z_l$ be jointly normal. Then, for every~$x \geq 0$, it holds that
\begin{equation}
\mathbb{P}\left(\max_{i=1, \hdots, l} |z_i| > x\right) \leq 1 - \prod_{i = 1}^l \mathbb{P}(|z_i| \leq x).
\end{equation}
\end{theorem}

The following lemma combines the Gaussian approximation for centered but \emph{non}-standardized sums over~$\mc{H}$ in Lemma~\ref{lem:cherno} with the bound in~\eqref{eqn:varestim**} from Corollary~\ref{cor:varestim} to establish a Gaussian approximation for centered and standardized sums. Its proof is similar to that of Theorem 4.1 in~\cite{kock2025high}. Set~$$\hat{D}_n:=\text{diag}\left(s_{n,1},\hdots,s_{n,d}\right),$$ and recall~$D_n(P)=\text{diag}(\sigma_{2,n,1}(P),\hdots,\sigma_{2,n,d}(P))$ and~$	\Sigma_{0,n}(P)=D_n^{-1}(P)\Sigma_n(P) D_n^{-1}(P)$ from~\eqref{eqn:sigma0def}. Furthermore, denote by~$\hat{D}^{-1}_n$ the diagonal matrix with diagonal entries~$s^{-1}_{n, j}$ for those indices~$j$ for which~$s_{n,j} > 0$ and~$0$ else.

\begin{lemma}\label{lem:GAStanfardizedSums}
Let~$b_1\in(0,\infty)$,~$b_2\in(b_1,\infty)$, and~$m > 4$. Then, there exists a (positive) constant~$C$, say, depending only on~$b_1$,~$b_2$ and~$m$, such that for every~$n\geq 3$,~$d\geq 2$, and~$P\in\mc{P}(b_1,b_2,m)$
\begin{equation}\label{eq:stsumsbound}
\begin{aligned}
&\sup_{H\in\mc{H}}\envert[3]{P\del[3]{\hat{D}_n^{-1}\frac{1}{\sqrt{n}}\sum_{i=1}^n\sbr[1]{X_{i,n}-\mu_{n}(P)}\in H}-P\del[1]{\Sigma_{0,n}^{1/2}(P)Z_d\in H}} \\
& \quad \quad \leq C \cdot \mathfrak{c}(n,d,m) + \frac{8}{\log(n)},
\end{aligned}
\end{equation}	
where~$$\mathfrak{c}(n,d,m) := \sqrt{\log(d)\log(dn)} \mathfrak{a}(n, d, m) + \mathfrak{b}(n, d, m).$$ If~$\frac{d}{n^{m/2-1-\xi}}\to 0$ for some~$\xi\in(0,1)$, then~$\mathfrak{c}(n,d,m) \to 0$.

\end{lemma}
\begin{proof}
Fix~$b_1, b_2, m, n, d$ and $P$ as in the present lemma's statement, and suppress the dependence of~$\sigma_{2, n,j}(P)$ and~$\mu_{n}(P)$ on~$P$ in the sequel. Let~$C_1$ denote the constant ``$C$'' from Corollary~\ref{cor:varestim}, which depends only on~$b_2$ and~$m$. We first note that if~$\mathfrak{a}(n, d, m) \geq\frac{b_1^2}{2C_1}$ holds, then the upper bound in~\eqref{eq:stsumsbound} exceeds~$1$ for any constant~$C$ satisfying~$C \geq 2C_1/b_1^2$. Because~$C_1/b_1^2$ only depends on~$b_1, b_2$ and~$m$, the condition~$C \geq 2C_1/b_1^2$ can be enforced (by enlarging the constant~$C$ that will be derived subsequently, if necessary). We shall therefore assume without loss of generality that~$\mathfrak{a}(n, d, m)< b_1^2/(2C_1)$. Set 
\begin{equation}
U_n := D_n^{-1}\frac{1}{\sqrt{n}}\sum_{i=1}^n\sbr[1]{X_{i,n}-\mu_{n}} \quad \text{and} \quad V_n := \hat{D}_n^{-1}\frac{1}{\sqrt{n}}\sum_{i=1}^n\sbr[1]{X_{i,n}-\mu_{n}}.
\end{equation}
By Corollary~\ref{cor:varestim}, for~$C_2 := \sqrt{2}C_1/b_1^3$ depending only on~$b_1,b_2$, and~$m$, with probability at least~$1-8/\log(n)$ it holds that
\begin{equation}
\min_{j = 1, \hdots, d} s_{n,j} \geq b_1/\sqrt{2} \quad \text{ and } \quad \max_{j=1,\hdots,d}|s^{-1}_{n,j}-\sigma^{-1}_{2, n,j}|
\leq
C_2\mathfrak{a}(n, d, m)\label{eq:An2},
\end{equation}
from which it follows that
\begin{equation}
\|U_n - V_n\|_\infty 
\leq 
C_2\mathfrak{a}(n, d, m) \enVert[3]{\frac{1}{\sqrt{n}}\sum_{i=1}^n[X_{i,n}-\mu_{n}]}_{\infty}
\label{eq:An1}.
\end{equation}
Using the Khatri-{\v{S}}id{\'a}k inequality (as stated in Theorem~\ref{thm:ksi}), we obtain
\begin{align*}
P\del{\max_{j=1,\hdots,d}|(\Sigma_n^{1/2}Z_d)_j|> b_2\sqrt{2\log(2dn)}} &  \leq  1 - \prod_{j = 1}^d P\del{|\sigma_{2, n,j} (Z_d)_j| \leq b_2\sqrt{2\log(2dn)}} \\
& = P\del{\max_{j=1,\hdots,d}|\sigma_{2, n,j} (Z_d)_j|> b_2\sqrt{2\log(2dn)}} \\&\leq 
\frac{1}{n},
\end{align*}
the last inequality following from a union bound,~$\max_{j=1,\hdots,d}\sigma_{2, n,j}\leq b_2$ as~$P\in\mc{P}(b_1,b_2,m)$, together with~$\P\del[0]{|\bm{z}|> t}\leq 2\exp(-t^2/2)$ for~$t\geq 0$ and~$\bm{z} \sim \mathsf{N}_1(0,1)$.
Thus, by Lemma \ref{lem:cherno}, and for a constant~$K_1$ depending only on~$b_1,b_2$, and~$m$, we have (absorbing~$n^{-1}$ into the dominating~$\mathfrak{b}(n, d, m)$)
\begin{equation}
P\del[3]{\big\|\frac{1}{\sqrt{n}}\sum_{i=1}^n[X_{i,n}-\mu_{n}]\big\|_{\infty} \leq b_2\sqrt{2\log(2dn)}}
\geq 1- 
K_1\mathfrak{b}(n, d, m)\label{eq:An3}.
\end{equation} 
Hence, by~\eqref{eq:An2}--\eqref{eq:An3} and for a constant~$C_3$ depending only on~$b_1, b_2$, and~$m$, we have
\begin{equation}\label{eq:An4}
P\del[3]{\|U_n - V_n\|_\infty \leq C_3\sqrt{\log(dn)}\mathfrak{a}(n, d, m)}
\geq 
1 - \frac{8}{\log(n)}- K_1 \mathfrak{b}(n, d, m).
\end{equation}
For any~$u \in \R^d$ and~$H' = \cbr[1]{x\in\R^d:a_j\leq x_j\leq b_j\text{ for all }j=1,\hdots, d}\in \mathcal{H}$, we have~$D_n^{-1} u \in H'$ if and only if $u \in \cbr[1]{x\in\R^d:a_j\sigma_{2, n,j}\leq x_j\leq b_j\sigma_{2, n,j}\text{ for all }j=1,\hdots, d} = D_n H' \in \mathcal{H}$. Furthermore, the (so-defined) map~$D_n$ acts bijectively on~$\mathcal{H}$ because~$\sigma_{2, n,j} \geq b_1 > 0$. Therefore,
\begin{equation}\label{eqn:CLTU}
\begin{aligned}
&\sup_{H' \in \mathcal{H}}\envert[3]{P\left(U_n \in H'\right)-P\del[1]{\Sigma_{0,n}^{1/2}Z_d\in H'}} \\& = \sup_{H' \in \mathcal{H}}
\envert[3]{P\bigg(\frac{1}{\sqrt{n}}\sum_{i=1}^n\sbr[1]{X_{i,n}-\mu_{n}} \in H'\bigg)-P\del[1]{\Sigma_{n}^{1/2}Z_d\in H'}} \leq
K_1\mathfrak{b}(n, d, m),
\end{aligned}
\end{equation}
the inequality following from Lemma~\ref{lem:cherno}. Equations~\eqref{eq:An4} and~\eqref{eqn:CLTU} show that Lemma \ref{lem:auxlem} is applicable with~$U = U_n$, $V = V_n$,~$\Xi = \Sigma_{0,n}$,~$\xi = 1$, $a = K_1\mathfrak{b}(n, d, m)$, $b = C_3\sqrt{\log(dn)}\mathfrak{a}(n, d, m)$,~$c = \frac{8}{\log(n)}+K_1 \mathfrak{b}(n, d, m)$, and for every~$H \in \mathcal{H}$. We conclude that for a constant~$C$ depending only on~$b_1$, $b_2$, and~$m$, the supremum in~\eqref{eq:stsumsbound} is bounded from above by
\begin{align*}
&a + (\sqrt{2 \log(d)} + 4)b \xi^{-1} + c \\
&= 2 K_1\mathfrak{b}(n, d, m) + (\sqrt{2 \log(d)} + 4) C_3\sqrt{\log(dn)}\mathfrak{a}(n, d, m) + \frac{8}{\log(n)}  \\
&\leq C\left(\sqrt{\log(d)\log(dn)} \mathfrak{a}(n, d, m) + \mathfrak{b}(n, d, m) \right) + 8/\log(n).
\end{align*}
We have therefore established~\eqref{eq:stsumsbound}. To verify the remaining statement, first recall that
\begin{equation}
\begin{aligned}\label{eqn:summaryabc}
\mathfrak{a}(n, d, m) &= \frac{\log(dn)}{n}+\sqrt{\frac{\log(dn)}{n}}+\frac{[d\log(n)]^{2/m}}{n^{1-2/m}},\\
\mathfrak{b}(n, d, m) &= \sbr[3]{\frac{d^{2/m}\log^{5}(d)}{n}}^{1/4}+\sbr[3]{\frac{d^{2/m}\log^{3-2/m}(d)}{n^{1-2/m}}}^{1/2},\\
\mathfrak{c}(n,d,m) &= \sqrt{\log(d)\log(dn)} \mathfrak{a}(n, d, m) + \mathfrak{b}(n, d, m).
\end{aligned}
\end{equation}
Now, if there exists a~$\xi\in(0,1)$ such that~$\frac{d}{n^{m/2-1-\xi}}\to 0$, then obviously (recall that~$m > 4$)~$$\sqrt{\log(d)\log(dn)} \mathfrak{a}(n, d, m) \to 0 \quad \text{ and } \quad \mathfrak{b}(n, d, m) \to 0,$$
so that~$\mathfrak{c}(n,d,m) \to 0$ follows.
\end{proof}

\subsection{Gaussian approximation for~$S_n$ when~$\enVert[0]{\mu_n(P)}_\infty$ is ``small''}

\begin{theorem}\label{thm:SumGA}
Let~$b_1\in(0,\infty)$,~$b_2\in(b_1,\infty)$, and~$m > 4$. Then, there exists a (positive) constant~$C$, say, depending only on~$b_1$,~$b_2$ and~$m$, such that, for every~$n\geq 3$,~$d\geq 2$, and~$P\in\mc{P}(b_1,b_2,m)$,
\begin{equation}\label{eqn:SUMGA}
\begin{aligned}
&\sup_{H\in\mc{H}}\envert[1]{P\del[1]{S_{n}\in H}-P\del[1]{\Sigma_{0,n}^{1/2}(P)Z_{d}+\sqrt{n}D_{n}^{-1}\mu_{n}(P)\in H}} \\
&\leq 
C \cdot \left[\mathfrak{c}(n, d, m) + \mathfrak{a}(n, d, m)\sqrt{n\log(d)}\|\mu_n\|_{\infty}\right] + \frac{16}{\log(n)}.
\end{aligned}
\end{equation}
Let a sequence~$P_{n} \in \mc{P}(b_1,b_2,m)$ and a subsequence~$n'$ of~$n$ be given, such that~$d' \geq 2$, and such that for some~$\xi\in(0,1)$ it holds that
\begin{equation}\label{eqn:locmu}
\mathfrak{a}(n', d', m)\sqrt{n'\log(d')}\|\mu_{n'}\|_{\infty} \to 0 \quad \text{and} \quad \frac{d'}{{n'}^{m/2-1-\xi}}\to 0.
\end{equation} 
Then, the upper bound in~\eqref{eqn:SUMGA} converges to~$0$ along~$n'$.
\end{theorem}
\begin{proof}
We first establish the upper bound in~\eqref{eqn:SUMGA}. Fix~$b_1, b_2, m, n, d$ and $P$ as in the present theorem's statement. For simplicity, we suppress the dependence of all quantities (in particular $\Sigma_{0,n}$ and~$\mu_{n}$) on~$P$. Let~$C_1$ denote the constant ``$C$'' from Corollary~\ref{cor:varestim}, which depends only on~$b_2$ and~$m$. Arguing as in the beginning of the proof of Lemma~\ref{lem:GAStanfardizedSums}, we can focus on the case when~$\mathfrak{a}(n, d, m)< b_1^2/(2C_1)$. Fix~$H \in \mathcal{H}$. For $T_n:=\hat{D}_n^{-1}\frac{1}{\sqrt{n}}\sum_{i=1}^n\sbr[1]{X_{i,n}-\mu_n}+\sqrt{n}D_n^{-1}\mu_n$ (recall the definition of~$\hat{D}_n^{-1}$ right before Lemma~\ref{lem:GAStanfardizedSums}) and~$(H-\sqrt{n}D_n^{-1}\mu_n) =: H^* \in \mc{H}$ we have
\begin{equation*}
P\del[1]{T_n\in H}
= P\del[1]{T_n - \sqrt{n}D_n^{-1}\mu_n\in H^*} = 
P\del[2]{\hat{D}_n^{-1}\frac{1}{\sqrt{n}}\sum_{i=1}^n\sbr[1]{X_{i,n}-\mu_n}\in H^*}.
\end{equation*}
Hence, for~$a_n := C \cdot  \mathfrak{c}(n, d, m) + 8/\log(n)$ (the constant~$C$ depending only on~$b_1$,~$b_2$ and~$m$) Lemma~\ref{lem:GAStanfardizedSums} yields~$|P(T_n - \sqrt{n}D_n^{-1}\mu_n\in H^*)-P(\Sigma_{0,n}^{1/2}Z_d\in H^*)|
\leq a_n.$ By the same argument
\begin{equation}\label{eq:boundTn}
\sup_{H' \in \mathcal{H}}	\envert[2]{P\del[1]{T_n - \sqrt{n}D_n^{-1}\mu_n\in H'}-P\del[1]{\Sigma_{0,n}^{1/2}Z_d\in H'}}
\leq a_n.
\end{equation}
Recall that~$C_1$ (depending only on~$b_1,b_2$, and~$m$) denotes the constant ``$C$'' from Corollary~\ref{cor:varestim},~$\mathfrak{a}(n, d, m)< b_1^2/(2C_1)$, and denote~$C_1':= \sqrt{2}C_1/b_1^3$. Hence, by Corollary~\ref{cor:varestim}, with probability at least~$1 - 8/\log(n)$ it holds that
\begin{equation*}
\min_{j = 1, \hdots, d} s_{n,j} \geq \frac{b_1}{\sqrt{2}} \quad \text{ and } \quad \max_{j=1,\hdots,d}|s^{-1}_{n,j}-\sigma^{-1}_{2, n,j}|
\leq 
C_1' \cdot \mathfrak{a}(n, d, m),
\end{equation*}
from which it follows that
\begin{equation*}
\enVert[1]{S_n-T_n}_{\infty}
=
\enVert[1]{\sqrt{n}(\hat{D}_n^{-1}-D_n^{-1})\mu_n}_{\infty}
\leq
\sqrt{n} C_1' \cdot \mathfrak{a}(n, d, m) \|\mu_n\|_\infty =: \overline{A}_n,
\end{equation*}
where we used that~$S_n
=
T_n+\sqrt{n}(\hat{D}_n^{-1}-D_n^{-1})\mu_n$ (recall Equation~\eqref{eqn:Sndef}, including our convention when variance estimates vanish). Thus,
it holds that
\begin{equation}\label{eq:boundSn}
P\del[1]{\enVert[1]{S_n-T_n}_{\infty}> \overline{A}_n}\leq 8/\log(n).
\end{equation}
Using~\eqref{eq:boundTn} and~\eqref{eq:boundSn}, we can now apply Lemma~\ref{lem:auxlem} with~$H = H^*$,~$U = T_n - \sqrt{n}D_n^{-1} \mu_n$, $V = S_n - \sqrt{n}D_n^{-1} \mu_n$, ~$\Xi = \Sigma_{0, n}$,~$\xi = 1$,~$a = a_n$, $b = \overline{A}_n$ and~$c = 8/\log(n)$ to conclude 
\begin{align*}
a_n + \left(\sqrt{2 \log(d)} + 4\right)\overline{A}_n + 8/\log(n)  &\geq \left| P(S_n - \sqrt{n}D_n^{-1} \mu_n \in H^*) - P(\Sigma_{0, n}^{1/2} Z_d \in H^*) \right| \\
&= \left| P(S_n \in H) - P(\Sigma_{0, n}^{1/2} Z_d + \sqrt{n}D_n^{-1}\mu_n \in H) \right|.
\end{align*}
There exists a constant~$C$, say, depending only on~$b_1,b_2$, and~$m$, such that the upper bound~$a_n +(\sqrt{2 \log(d)} + 4)\overline{A}_n + 8/\log(n)$ is dominated by
\begin{equation}
C \cdot \left[\mathfrak{c}(n, d, m) + \mathfrak{a}(n, d, m)\sqrt{n\log(d)}\|\mu_n\|_{\infty}\right] + 16/\log(n),
\end{equation}
establishing~\eqref{eqn:SUMGA}. To prove the remaining statement, recall (from the proof of  Lemma~\ref{lem:GAStanfardizedSums} around~\eqref{eqn:summaryabc}, but now along the subsequence~$n'$) that under the second condition in~\eqref{eqn:locmu} it holds that~$\mathfrak{c}(n', d', m) \to 0$. The first condition imposes~$\mathfrak{a}(n, d, m)\sqrt{n\log(d)}\|\mu_n\|_{\infty} \to 0$ along the subsequence~$n'$, so that we can conclude.
\end{proof}

\subsection{Gaussian approximation for~$S_n$ when~$\enVert[0]{\mu_n(P)}_\infty$ is not ``small''}

We consider the asymptotic behavior of~$S_{n}$ when the first condition in~\eqref{eqn:locmu} is not satisfied.

\begin{theorem}\label{thm:SumGAnonloc}
Let~$b_1\in(0,\infty)$,~$b_2\in(b_1,\infty)$, and~$m > 4$. Let a sequence~$P_{n} \in \mc{P}(b_1,b_2,m)$ and a subsequence~$n'$ of~$n$ be given, such that~$d' \geq 2$, and such that for some~$\xi\in(0,1)$ it holds that
\begin{equation}\label{eqn:locmunot}
\mathfrak{a}(n', d', m)\sqrt{n'\log(d')}\|\mu_{n'}\|_{\infty} \to c \in (0, \infty] \quad \text{and} \quad \frac{d'}{{n'}^{m/2-1-\xi}}\to 0.
\end{equation} 
Then, for every sequence~$\mathcal{H}^*_n$ satisfying $\emptyset \neq \mathcal{H}^*_n \subseteq \mathcal{H}_n$ along~$n'$ and such that $$r_{n'} := \sup \left \{\|x\|_{\infty}: x \in H,~ H \in \mathcal{H}^*_{n'}\right \}= O(\sqrt{\log(d')}),$$ it holds that
\begin{equation}\label{eqn:SUMGA2}
\sup_{H\in\mc{H}^*_{n'} }\envert[1]{P_{n'}\del[1]{S_{n'}\in H}-P_{n'}\del[1]{\Sigma_{0,n'}^{1/2}(P_{n'})Z_{d'}+\sqrt{n'}D_{n'}^{-1}\mu_{n'}(P_{n'})\in H}} \to 0;
\end{equation}
more specifically, for every sequence~$H_n$ such that~$H_n \in \mathcal{H}^*_n$ along~$n'$, we then have 
\begin{equation}\label{eqn:SUMGA3}
P_{n'}\del[1]{S_{n'}\in H_{n'}} \to 0 \quad \text{  and } \quad P_{n'}\del[1]{\Sigma_{0,n'}^{1/2}(P_{n'})Z_{d'}+\sqrt{n'}D_{n'}^{-1}\mu_{n'}(P_{n'})\in H_{n'}} \to 0.
\end{equation}
\end{theorem}

\begin{proof}
Fix~$b_1, b_2, m, n', d'$ and $P_n$ as in the present theorem's statement. We establish~\eqref{eqn:SUMGA3}, write~$n$ instead of~$n'$, and drop the dependence of several quantities on~$P_n$ when possible. Write $W_n=\hat{D}_n^{-1}\frac{1}{\sqrt{n}}\sum_{i=1}^n\sbr[1]{X_{i,n}-\mu_n}$ (recall the definition of~$\hat{D}_n^{-1}$ right before Lemma~\ref{lem:GAStanfardizedSums}) and lower bound, using~$H_n \subseteq \{x \in \R^d: -r_n \leq x_j \leq r_n \text{ for all } j = 1, \hdots, d\}$ (recall that~$H_n\in\mc{H}_n^*$, and Equation~\eqref{eqn:Sndef}, including our convention when variance estimates vanish),
\begin{equation*}
P_n\del[1]{S_n \notin H_n} \geq P_n\del[1]{\|S_n\|_{\infty} > r_n} \\
= P_n\left( \|W_n +\sqrt{n}D_n^{-1}\mu_n+\sqrt{n}(\hat{D}_n^{-1}-D_n^{-1})\mu_n\|_\infty > r_n\right),
\end{equation*}
which is no smaller than
\begin{equation}\label{eq:LBsums}
P_n\left(\|\sqrt{n}D_n^{-1}\mu_n+\sqrt{n}(\hat{D}_n^{-1}-D_n^{-1})\mu_n\|_\infty > r_n + \|W_n\|_\infty\right).
\end{equation}
Lemma~\ref{lem:GAStanfardizedSums} delivers
\begin{equation*}
\sup_{H\in\mc{H}_n} \left| P_n\left(W_n \in H\right)-P_n\left(\Sigma_{0,n}^{1/2}Z_d\in H\right) \right|
\to 0.
\end{equation*}	
By the Khatri-{\v{S}}id{\'a}k inequality (and a union bound)~$\|\Sigma_{0,n}^{1/2}Z_d\|_{\infty} = O_{P_n}(\sqrt{\log(d)})$, so that
\begin{equation}\label{eqn:rhslogn}
r_n + \|W_n\|_\infty = O_{P_n}\del[1]{\sqrt{\log(d)}} \quad \text{ and } \quad r_n + \|\Sigma_{0,n}^{1/2}Z_d\|_{\infty} = O_{P_n}\del[1]{\sqrt{\log(d)}}.
\end{equation}
That~$d/n^{m/2-1-\xi}\to 0$ for some~$\xi\in(0,1)$ also implies (recall that~$m > 4$ holds)
\begin{equation}\label{eqn:arate}
\mathfrak{a}(n, d, m) = \frac{\log(dn)}{n}+\sqrt{\frac{\log(dn)}{n}}+\frac{[d\log(n)]^{2/m}}{n^{1-2/m}} = o\left(\left[\frac{\log(n)}{n^{\xi}}\right]^ {2/m}\right).
\end{equation}
Therefore, without loss of generality~$\mathfrak{a}(n,d,m) < b_1^2/(2C_1)$ for~$C_1$ as in Corollary~\ref{cor:varestim}. The latter corollary implies that for a constant~$C'$ depending only on~$b_1, b_2$ and~$m$, it holds that (without loss of generality, we assume that~$n \geq 3$)
\begin{equation}\label{eqn:intermedins}
\min_{j = 1, \hdots, d} s_{n,j} \geq b_1/\sqrt{2} \text{ and } \enVert[0]{\hat{D}_n^{-1}-D_n^{-1}}_\infty \leq C'\mathfrak{a}(n, d, m)
\end{equation}
with probability at least~$1-8/\log(n)$.
Abbreviate~$\sigma_{n,j}=\sigma_{2,n,j}$. From~$\sigma_{n, j} \leq b_2$ for~$j = 1, \hdots, d$, it follows that
\begin{equation}\label{eqn:mulow*}
\enVert[1]{D_n^{-1}\mu_n}_\infty = \max_{j = 1, \hdots, d} |\mu_{n,j}|/\sigma_{n,j} \geq  \|\mu_{n}\|_{\infty} / b_2,
\end{equation} so that
\begin{align*}
\enVert[2]{\sqrt{n}D_n^{-1}\mu_n+\sqrt{n}(\hat{D}_n^{-1}-D_n^{-1})\mu_n}_\infty
&\geq
\enVert[1]{\sqrt{n}D_n^{-1}\mu_n}_\infty-\enVert[1]{\sqrt{n}(\hat{D}_n^{-1}-D_n^{-1})\mu_n}_\infty\\
&\geq
\sqrt{n}\enVert[0]{\mu_n}_\infty\del[1]{ b^{-1}_2-\enVert[0]{\hat{D}_n^{-1}-D_n^{-1}}_\infty}.
\end{align*}
Because~$\mathfrak{a}(n, d, m)\sqrt{n\log(d)}\|\mu_{n}\|_{\infty} \geq C$, for some positive real number~$C$, say, at least for~$n$ large enough (which we assume without loss of generality), it follows from the statement around~\eqref{eqn:intermedins} that with probability at least~$1-8/\log(n)$ the far left-hand side of the previous display is lower bounded by
\begin{equation}\label{eqn:calogbd}
\frac{C}{\mathfrak{a}(n, d, m) \sqrt{\log(d)}}\del[1]{b_2^{-1}-C'\mathfrak{a}(n, d, m)} \geq \frac{C}{\mathfrak{a}(n, d, m) \sqrt{\log(d)}},
\end{equation}
the inequality (cf.~also~\eqref{eqn:arate}) holding for~$n$ sufficiently large (upon adjusting~$C$ downward). Summarizing, with probability converging to~$1$  $$\enVert[2]{\sqrt{n}D_n^{-1}\mu_n+\sqrt{n}(\hat{D}_n^{-1}-D_n^{-1})\mu_n}_\infty  \geq \frac{C}{\mathfrak{a}(n, d, m) \sqrt{\log(d)}}.$$ Hence, recalling~\eqref{eq:LBsums}, we have
\begin{align*}
P_n\del[1]{S_n \notin H_n} &\geq P_n\left(\|\sqrt{n}D_n^{-1}\mu_n+\sqrt{n}(\hat{D}_n^{-1}-D_n^{-1})\mu_n\|_\infty > r_n + \|W_n\|_\infty\right) \\
& \geq
P_n\left(C > \mathfrak{a}(n, d, m) \sqrt{\log(d)}[r_n + \|W_n\|_\infty ]\right) + o(1) \to 1,
\end{align*}
the convergence following from~\eqref{eqn:rhslogn} and~\eqref{eqn:arate}. 

To prove the second statement in~\eqref{eqn:SUMGA3}, note that
\begin{align*}
P_{n}\del[1]{\Sigma_{0,n}^{1/2}Z_{d}+\sqrt{n}D_{n}^{-1}\mu_{n}\notin H_{n}} &\geq  P_n\del[1]{\enVert[0]{\Sigma_{0,n}^{1/2}Z_d+\sqrt{n}D_{n}^{-1}\mu_{n}}_\infty> r_n}\\
&\geq
P_n\del[1]{\sqrt{n} \enVert[0]{D_{n}^{-1}\mu_{n}}_\infty>r_n + \enVert[0]{\Sigma_{0,n}^{1/2}Z_d}_\infty} \\
&\geq P_n(C \geq b_2\mathfrak{a}(n, d, m)\sqrt{\log(d)}[r_n + \enVert[0]{\Sigma_{0,n}^{1/2}Z_d}_\infty]) \to 1,
\end{align*}
where we used~\eqref{eqn:mulow*} and $\mathfrak{a}(n, d, m)\sqrt{n\log(d)}\|\mu_{n}\|_{\infty} \geq C > 0$ in the third inequality, and the convergence follows (as above) from~\eqref{eqn:rhslogn} and~\eqref{eqn:arate}. 
\end{proof}

\subsection{Proof of Theorem~\ref{thm:sums}}

Before we present the proof of Theorem~\ref{thm:sums}, we recall the following elementary fact. 
\begin{lemma}\label{lem:unifseq}
Let~$A$ be a non-empty set, and~$b_n:A \to [0, \infty)$ be a sequence of functions. Then, the following statements are equivalent:
\begin{enumerate}
\item $\lim_{n \to \infty} \sup_{a \in A} b_n(a) = 0$.
\item For every sequence~$a_n \in A$, we have~$b_n(a_n) \to 0$.
\item For every sequence~$a_n \in A$ and every subsequence~$n'$ of~$n$, there exists a subsequence~$n''$ of~$n'$, such that~$b_{n''}(a_{n''})  \to 0$.
\end{enumerate}
\end{lemma}

\begin{proof}[Proof of Theorem~\ref{thm:sums}:]
To establish Part 1, let~$P_n$ be a sequence in~$\mc{P}^0(b_1,b_2,m)$, and suppress the dependence of all quantities on~$P_n$. Theorem~\ref{thm:SumGA} gives (note that~$\mu_n = 0$ implies~\eqref{eqn:locmu})
\begin{equation*}
\envert[1]{P_{n}\del[1]{\enVert[0]{S_{n}}_\infty>c_{n,1-\alpha}}-P_n\del[1]{\enVert[0]{\Sigma_{0,n}^{1/2}Z_d}_\infty>c_{n,1-\alpha}}}\to 0.	
\end{equation*}
By the Khatri-{\v{S}}id{\'a}k inequality (the inequality being an equality in case~$\Sigma_{0,n} = I_d$)
\begin{equation*}
P_n\del[1]{\enVert[0]{\Sigma_{0,n}^{1/2}Z_d}_\infty>c_{n,1-\alpha}}
\leq 
P_n\del[1]{\enVert[0]{Z_d}_\infty>c_{n,1-\alpha}} = \alpha,
\end{equation*}
which establishes the first two statements in Part~1 (respectively choosing~$P_n$ to approximate the supremum in~\eqref{eq:sizecherno}, or~$P_n$ satisfying~$\Sigma_{0, n}(P_n) = I_d$ for~$n$ large enough). The last statement follows from Lemma~\ref{lem:extreme}, which shows that if~$d \to \infty$, then for every sequence~$P_n$ in~$\mathcal{P}(b_1, b_2, m, \bm{r})$ we have~$P_n(\|\Sigma_{0,n}^{1/2}Z_d\|_\infty>c_{n,1-\alpha}) \to \alpha$ under the assumption that~$\log(l)r_l \to 0$ as $l\to \infty$.

To establish~\eqref{eq:GaussApproxCherno}, we use the equivalence between Parts 1 and 3 in Lemma~\ref{lem:unifseq}. Let~$P_n$ be a sequence in~$\mc{P}(b_1,b_2,m)$. Let~$n'$ be a subsequence of~$n$. If~$\mathfrak{a}(n', d', m)\sqrt{n'\log(d')}\|\mu_{n'}\|_{\infty} \to 0$, then Theorem~\ref{thm:SumGA} shows that
\begin{equation}\label{eqn:locconv}
\envert[2]{P_{n'}\del[1]{\enVert[0]{S_{n'}}_\infty>c_{n',1-\alpha}}-P_{n'}\del[1]{\enVert[0]{\Sigma_{0,n'}^{1/2}Z+\sqrt{n'}D_{n'}^{-1}\mu_{n'}}_\infty>c_{n',1-\alpha}}}
\to 0.
\end{equation}
Otherwise,~$n'$ possesses a subsequence~$n''$ along which~$\mathfrak{a}(n'', d'', m)\sqrt{n''\log(d'')}\|\mu_{n''}\|_{\infty} \to c \in (0, \infty]$. Theorem~\ref{thm:SumGAnonloc} together with~$c_{n,1-\alpha} = O(\sqrt{\log(d)})$ then shows~\eqref{eqn:locconv} along~$n''$.
\end{proof}

\subsection{Proof of Theorem \ref{thm:Breakdown}}
Fix~$n\geq 2$ and~$\alpha\in(0,1)$. Note that~$\tilde{X}_{n,n}=-\sum_{i=1}^{n-1}\tilde{X}_{i,n}$ implies that~$\sum_{i=1}^n \tilde{X}_{i,n}=0_d$. Thus, also~$\overline{\tilde{X}}_n=0_d$ and hence (interpreting a sum over an empty index set as~$0$, if necessary)
\begin{equation*}
n \tilde{s}_{n,j}^2=\sum_{i=1}^{n-2}\tilde{X}_{i,n,j}^2+1+\del[3]{-\sum_{i=1}^{n-1} \tilde{X}_{i,n,j}}^2>0,\qquad \text{for }j=1,\hdots,d.
\end{equation*}
Thus,
\begin{align*}
\tilde{S}_{n,j}
=
\frac{\sum_{i=1}^n \tilde{X}_{i,n,j}}{\sqrt{n} \tilde{s}_{n,j}}
=0,\qquad\text{for }j=1,\hdots,d,
\end{align*}
and also~$\tilde{\varphi}_n=\mathds{1}\cbr[0]{||\tilde{S}_n||_\infty>c_{n,1-\alpha}}=0$ because~$c_{n,1-\alpha}>0$.

\section{Proof of Theorem~\ref{thm:Wins}}
\subsection{Auxiliary results concerning variance estimation}
Recall that for~$j=1,\hdots,d$ we denote~$\tilde{\sigma}_{n,j}^2=\tilde{\Sigma}_{n,j,j}$, with~$\tilde{\Sigma}_{n}$ defined in~\eqref{eq:tildeSigma}.
For the purpose of using the finite-sample bounds on the precision of~$\tilde{\Sigma}_n$ from~\cite{kock2025high}, we introduce the condition
\begin{align}\label{eq:epscond'}
2\eps_n' +\frac{\log(d^2n)}{n}+\sqrt{\del[2]{\frac{\log(d^2n)}{n}}^2+4\frac{\log(d^2n)}{n}\eps_n'}<1.	
\end{align}
We stress that \eqref{eq:epscond'} is \emph{not} needed in 
\begin{itemize}
\item[i)] any of our results on asymptotic size and power of the tests~$\varphi_{n,W}$ and~$\varphi_{n,W,B}$, as the conditions imposed in Theorems~\ref{thm:Wins},~\ref{thm:BSWins}, and~\ref{thm:comparison} imply that~\eqref{eq:epscond'} is eventually satisfied.
\item[ii)] the finite-sample Gaussian approximation in \eqref{eqn:SUMGASW} of Theorem~\ref{thm:SumGAWins} below, as this approximation is argued to be trivially valid when~\eqref{eq:epscond'} does not hold (cf.~the argument in the beginning of the proof of Theorem~\ref{thm:SumGAWins}).
\end{itemize}

We emphasize that the statement in~\eqref{eq:covestimGram} below is Theorem 3.1 in \cite{kock2025high} translated to our setting laid out in Section~\ref{sec:setting}.

\begin{lemma}\label{lem:varwinsestim}
Let~$b_1\in(0,\infty)$,~$b_2\in(b_1,\infty)$, and~$m > 2$. Assume that~$n > 3$ and~$d \geq 2$. Fix~$\lambda_1'\in(1,\infty)$,~$\lambda_2'\in(0,\infty)$ and let
\begin{equation}\label{eqn:ddef}
\mathfrak{d}(n,d, m, \overline{\eta}_n) := \overline{\eta}_n^{1-\frac{2}{m}}+\del[2]{\frac{\log(dn)}{n}}^{1-\frac{1}{(m/2)\wedge 2}}.
\end{equation}
If~$\varepsilon_n'$ is as in~\eqref{eq:epsprime} and satisfies~\eqref{eq:epscond'}, then, for a (positive) constant~$C$, say, depending only on~$b_2,\lambda_1',\lambda_2',m$ and for all~$P\in\mc{P}(b_1,b_2,m)$, we have
\begin{equation}\label{eq:covestimGram}
P\del[4]{\max_{1\leq j,k\leq d}\envert[1]{\tilde{\Sigma}_{n,j,k}-\Sigma_{n,j,k}(P)}> C \cdot \mathfrak{d}(n,d, m, \overline{\eta}_n)}
\leq \frac{24}{n};
\end{equation}	
if furthermore~$\mathfrak{d}(n, d, m, \overline{\eta}_n) <\frac{b_1^2}{2C}$, then, setting~$C' := \sqrt{2}C/b_1^3$, it holds that
\begin{equation}\label{eq:winsvarinv}
P\del[4]{\min_{j = 1, \hdots, d} \tilde{\sigma}_{n,j} \geq \frac{b_1}{\sqrt{2}} \text{ and } \max_{j=1,\hdots,d}\envert[1]{\tilde{\sigma}_{n,j}^{-1}-\sigma_{2,n,j}^{-1}(P)}\leq C'\cdot \mathfrak{d}(n,d, m, \overline{\eta}_n)}\geq1 - \frac{24}{n},
\end{equation}
and there exists a (positive) constant~$C''$, say, depending only on~$b_1, b_2, \lambda_1',\lambda_2', m$, such that 
\begin{equation}\label{eqn:correst}
P\del[4]{\min_{j = 1, \hdots, d} \tilde{\sigma}_{n,j} \geq \frac{b_1}{\sqrt{2}} \text{ and } \max_{1\leq j,k\leq d}\envert[1]{\tilde{\Sigma}_{0,n,j,k}-\Sigma_{0,n, j,k}(P)}
\leq
C'' \cdot\mathfrak{d}(n,d, m, \overline{\eta}_n)}
\geq 1 -\frac{24}{n}.
\end{equation}
\end{lemma}
\begin{proof}
The statement in~\eqref{eq:covestimGram} is Theorem 3.1 in~\cite{kock2025high} (where~$n > 3$ and~$d \geq 2$ is maintained throughout, and which holds by assumption). To establish~\eqref{eq:winsvarinv} we only need to slightly adapt the argument used to derive~\eqref{eqn:varestim**} from~\eqref{eqn:varestim*} in the proof of Corollary~\ref{cor:varestim} [specifically: replace~$s_{n,j}$ by~$\tilde{\sigma}_{n,j}$, $\mathfrak{a}(n,d,m)$ by~$\mathfrak{d}(n,d,m,\overline{\eta}_n)$, and the references to~\eqref{eqn:varestim*} and~\eqref{eqn:varestim**} by references to~\eqref{eq:covestimGram} and~\eqref{eq:winsvarinv}, respectively].

It remains to verify~\eqref{eqn:correst}. We suppress the dependence of several quantities on~$P$ in what follows. By~\eqref{eq:covestimGram} there exists an event~$E_n$ of probability at least~$1-\frac{24}{n}$ on which
\begin{equation}\label{eqn:alrsho}
\max_{1\leq j,k\leq d}\envert[1]{\tilde{\Sigma}_{n,j,k}-\Sigma_{n,j,k}}
\leq	
C\mathfrak{d}(n,d, m, \overline{\eta}_n) \leq \frac{b_1^2}{2}.
\end{equation}
Because by assumption~$\min_{j=1,\hdots,d}\sigma_{2,n,j}^2\geq b_1^2$, we have~$\min_{j=1,\hdots,d}\tilde{\sigma}_{n,j}^2\geq b_1^2/2$ on~$E_n$. Abbreviating~$\sigma_{n,j}=\sigma_{2,n,j}$, note that for all~$1\leq j,k\leq d$ and on~$E_n$
\begin{align*}
\tilde{\Sigma}_{0,n,j,k}-\Sigma_{0,n,j,k}
&=
\frac{\tilde{\Sigma}_{n,j,k}}{\tilde{\sigma}_{n,j}\tilde{\sigma}_{n,k}}-\frac{\Sigma_{n,j,k}}{\sigma_{n,j}\sigma_{n,k}}\\
&=
\frac{(\tilde{\Sigma}_{n,j,k}-\Sigma_{n,j,k})\sigma_{n,j}\sigma_{n,k}+\Sigma_{n,j,k}(\sigma_{n,j}\sigma_{n,k}-\tilde{\sigma}_{n,j}\tilde{\sigma}_{n,k})}{\tilde{\sigma}_{n,j}\tilde{\sigma}_{n,k}\sigma_{n,j}\sigma_{n,k}},
\end{align*}
so that, recalling~$\sigma_{n,j}\leq b_2$ for all~$j=1,\hdots,d$, for some constant~$C''$ depending only on~$b_1,b_2,\lambda_1',\lambda_2'$, and~$m$, we have on~$E_n$ (where~\eqref{eqn:alrsho} implies~$\tilde{\sigma}_{n,j}\leq \sqrt{b_2^2 + b_1^2/2}$ for all~$j=1,\hdots,d$) that
\begin{equation}
\max_{1\leq j,k\leq d}\envert[1]{\tilde{\Sigma}_{0,n,j,k}-\Sigma_{0,n, j,k}} \leq C'' \max_{1\leq j,k\leq d}\envert[1]{\tilde{\Sigma}_{n,j,k}-\Sigma_{n,j,k}}.
\end{equation}
We hence conclude with~\eqref{eqn:alrsho} by adjusting~$C''$.
\end{proof}

\subsection{Gaussian approximation for~$S_{n,W}$ when~$\enVert[0]{\mu_n(P)}_\infty$ is ``small''}
The following statement is an analogue to Theorem~\ref{thm:SumGA}, but now for~$S_{n,W}$ instead of~$S_n$. 
\begin{theorem}\label{thm:SumGAWins}
Let~$b_1\in(0,\infty)$,~$b_2\in(b_1,\infty)$, and~$m > 2$. Fix~$\lambda_1,\lambda_1'\in(1,\infty)$, and~$\lambda_2,\lambda_2'\in(0,\infty)$. Then, there exists a (positive) constant~$C$, say, depending only on~$b_1$,~$b_2$,~$\lambda_1$,~$\lambda_2$,~$\lambda_1'$,~$\lambda_2'$, and~$m$, such that, for every~$n > 3$,~$d\geq 2$, and~$P\in\mc{P}(b_1,b_2,m)$, if~$\eps_n,\eps_n'\in(0,1/2)$, with~$\eps_n$ as in~\eqref{eq:epsfam} and~$\eps_n'$ as in~\eqref{eq:epsprime}, it holds that
\begin{equation}\label{eqn:SUMGASW}
\begin{aligned}
&\sup_{H\in\mc{H}}\envert[1]{P\del[1]{S_{n, W}\in H}-P\del[1]{\Sigma_{0,n}^{1/2}(P)Z_{d}+\sqrt{n}D_{n}^{-1}\mu_{n}(P)\in H}} \\
&\leq 
C \cdot \left[\mathfrak{e}(n, d, m, \overline{\eta}_n) + \mathfrak{d}(n, d, m, \overline{\eta}_n)\sqrt{n\log(d)}\|\mu_n\|_{\infty}\right] + \frac{24}{n},
\end{aligned}
\end{equation}
where (cf.~\eqref{eqn:ddef} for the definition of~$\mathfrak{d}(n, d, m, \overline{\eta}_n)$)
\begin{align*} 
\mathfrak{e}(n, d, m, \overline{\eta}_n) :=~& \mathfrak{f}(n, d, m, \overline{\eta}_n)+\sqrt{\log(d)\log(dn)}\mathfrak{d}(n,d, m, \overline{\eta}_n), \\
\mathfrak{f}(n, d, m, \overline{\eta}_n) :=~& \del[4]{\sbr[3]{\frac{\log^{5-\frac{2}{m}}(dn)}{n^{1-\frac{2}{m}}}}^{\frac{1}{4}}
+
\sbr[3]{\overline{\eta}_n^{1-\frac{1}{m}}+\sbr[2]{\frac{\log(dn)}{n}}^{1-\frac{1}{m}}}\sqrt{n\log(d)}} \notag \\
&+
\del[4]{\log^2(d)\sbr[3]{\overline{\eta}_n^{1-\frac{2}{m}}+\sbr[2]{\frac{\log(dn)}{n}}^{1-\frac{2}{m}}}}^{1/2}, 
\end{align*}
Let a sequence~$P_{n} \in \mc{P}(b_1,b_2,m)$ and a subsequence~$n'$ of~$n$ be given, such that~$d' \geq 2$, 
\begin{equation}\label{eqn:divndeta}
\sqrt{n'\log(d')}\overline{\eta}_{n'}^{1-\frac{1}{m}}\to 0 \quad \text{ and } \quad \log(d')/{n'}^{\frac{m-2}{5m-2}}\to 0,
\end{equation}
and
\begin{equation}\label{eq:extrasslocgauss}
\sqrt{n'\log(d')} \mathfrak{d}(n',d',m,\overline{\eta}_{n'})\|\mu_{n'}(P_{n'})\|_\infty	
\to
0.
\end{equation}
Then the upper bound in~\eqref{eqn:SUMGASW} converges to~$0$ along~$n'$.
\end{theorem}

Before we present a proof of Theorem~\ref{thm:SumGAWins}, let us note that, similar to the proof of Theorem~\ref{thm:SumGA} which relies on Lemma~\ref{lem:GAStanfardizedSums}, the proof of Theorem~\ref{thm:SumGAWins} makes use of the following Gaussian approximation result for centered and standardized sums, which is a version of Theorem 4.1 in~\cite{kock2025high}, and which we re-state for the convenience of the reader (using some notation already introduced in stating Theorem~\ref{thm:SumGAWins}). To state that result, 
let~$S_{n,W,c}\in\R^d$ be defined via its coordinates (with the convention of setting it~$0$ in case~$\tilde{\sigma}_{n,j} = 0$)
\begin{equation}\label{eqn:snwcdef}
S_{n,W,c,j}
:=
\frac{1}{\sqrt{n}\tilde{\sigma}_{n,j}}\sum_{i=1}^n\sbr[1]{\phi_{\hat\alpha_j,\hat\beta_j}(\tilde{X}_{i,n,j})-\mu_{n,j}},\qquad j=1,\hdots,d.
\end{equation}
Note that~\cite{kock2025high} denote the expression in~\eqref{eqn:snwcdef} by~``$S_{n, W, S}$,'' which they leave undefined in case a variance estimate vanishes. In terms of our notation, and with our convention, their Theorem 4.1 delivers the following.
\begin{theorem}[a version of Theorem 4.1 in~\cite{kock2025high}]\label{thm:41kp}
Let~$b_1\in(0,\infty)$,~$b_2\in(b_1,\infty)$, and~$m > 2$. Fix~$\lambda_1,\lambda_1'\in(1,\infty)$, and~$\lambda_2,\lambda_2'\in(0,\infty)$. Then, there exists a (positive) constant~$C^*$, say, depending only on~$b_1$,~$b_2$,~$\lambda_1$,~$\lambda_2$,~$\lambda_1'$,~$\lambda_2'$, and~$m$, such that, for every~$n > 3$,~$d\geq 2$, and~$P\in\mc{P}(b_1,b_2,m)$, if~$\eps_n,\eps_n'\in(0,1/2)$, with~$\eps_n$ as in~\eqref{eq:epsfam} and~$\eps_n'$ as in~\eqref{eq:epsprime}, it holds that
\begin{equation}\label{eqn:SUMGASWcent}
\begin{aligned}
\sup_{H\in\mc{H}}\envert[1]{P\del[1]{S_{n, W, c}\in H}-P\del[1]{\Sigma_{0,n}^{1/2}(P)Z_{d} \in H}} 
\leq 
C^* \cdot \mathfrak{e}(n, d, m, \overline{\eta}_n),
\end{aligned}
\end{equation}
The upper bound in~\eqref{eqn:SUMGASWcent} converges to~$0$ along any subsequence~$n'$ of~$n$ for which~\eqref{eqn:divndeta} holds.
\end{theorem}
For notational convenience, let~$\tilde{D}_n=\text{diag}(\tilde{\sigma}_{n,1},\hdots,\tilde{\sigma}_{n,d})$, and denote by~$\tilde{D}_n^{-1}$ the diagonal matrix with diagonal entries~$\tilde{\sigma}^{-1}_{n,j}$ for those indices~$j$ for which~$\tilde{\sigma}_{n,j} > 0$ and~$0$ else.

\begin{proof}[Proof of Theorem~\ref{thm:SumGAWins}:]
We argue analogously to the proof of Theorem~\ref{thm:SumGA}, but instead of basing our argument on Corollary~\ref{cor:varestim} and Lemma~\ref{lem:GAStanfardizedSums}, we use Lemma~\ref{lem:varwinsestim} and Theorem 4.1 in~\cite{kock2025high} (in the version stated in Theorem~\ref{thm:41kp} above), respectively. We first establish the upper bound in~\eqref{eqn:SUMGASW}. Fix~$b_1, b_2, m, n, \lambda_1,\lambda_2,\lambda_1',\lambda_2', d$ and $P$ as in the present theorem's statement. Suppose that~$\eps_n,\eps_n'\in(0,1/2)$. For simplicity, we suppress the dependence of all quantities (in particular $\Sigma_{0,n}$ and~$\mu_{n}$) on~$P$.

Let~$C_1$ denote the constant ``$C$'' from Lemma~\ref{lem:varwinsestim}, which depends only on~$b_2$,~$\lambda_1'$,~$\lambda_2'$, and~$m$. We first note that if~$\mathfrak{d}(n, d, m, \overline{\eta}_n) \geq\frac{b_1^2}{2C_1}$ holds, then the upper bound in~\eqref{eqn:SUMGASW} exceeds~$1$ for any constant~$C$ satisfying~$C \geq 2C_1/b_1^2$. Because~$C_1/b_1^2$ only depends on~$b_1, b_2$, $\lambda_1'$,~$\lambda_2'$, and~$m$, the condition~$C \geq 2C_1/b_1^2$ can be enforced. We can therefore assume without loss of generality that~$\mathfrak{d}(n, d, m, \overline{\eta}_n) < b_1^2/(2C_1)$. Similarly, it follows from~$\frac{\log(d^2n)}{n}\leq {\lambda_2'}^{-1}{\eps_n'}$ that for a constant~$C_2$ depending only on~$\lambda_2'$ the left-hand side of~\eqref{eq:epscond'} is bounded from above by
\begin{align*}
2\eps_n'+{\lambda_2'}^{-1}{\eps_n'}+\sqrt{{\lambda_2'}^{-2}{\eps_n'}^2+4{\lambda_2'}^{-1}{\eps_n'}^2}
\leq
C_2\eps_n'.	
\end{align*}
Note that if~$C_2\eps_n'\geq 1$ then by the definition of~$\eps_n'$ in~\eqref{eq:epsprime}
\begin{align*}
\overline{\eta}_n\geq \frac{1}{2C_2\lambda_1'}\qquad\text{or}\qquad \frac{\log(dn)}{n}
=
\frac{\log([dn]^2)/2}{n}
\geq
\frac{\log(d^2n)}{2n}
\geq
\frac{1}{4C_2\lambda_2'}.
\end{align*} 
In either case, the upper bound in~\eqref{eqn:SUMGASW} is valid after suitably adjusting~$C$ (if needed). Thus, we can also assume that~$C_2\eps_n'< 1$, implying that~\eqref{eq:epscond'} is satisfied. In particular~\eqref{eq:winsvarinv} holds.

Let~$H \in \mathcal{H}$, and recall~$S_{n,W,c}\in\R^d$ from~\eqref{eqn:snwcdef} (which we set~$0$ in case~$\tilde{\sigma}_{n,j} = 0$). For~$T_{n,W}:=S_{n,W,c}+\sqrt{n}D_n^{-1}\mu_n$ and~$(H-\sqrt{n}D_n^{-1}\mu_n) =: H^* \in \mc{H}$ we have~$P\del[1]{T_{n,W}\in H}
= P(S_{n,W,c} \in H^*)$. Hence, for~$a_n := C^* \cdot \mathfrak{e}(n, d, m, \overline{\eta}_n)$ (the constant~$C^*$ depending only on~$b_1$,~$b_2$,~$\lambda_1,\lambda_2,\lambda_1',\lambda_2'$, and~$m$), Theorem 4.1 in~\cite{kock2025high} (in the version stated in Theorem~\ref{thm:41kp} above) yields~$$|P(T_{n,W} - \sqrt{n}D_n^{-1}\mu_n\in H^*)-P(\Sigma_{0,n}^{1/2}Z_d\in H^*)|
\leq a_n.$$ By the same argument,
\begin{equation}\label{eq:boundTnW}
\sup_{H' \in \mathcal{H}}	\envert[2]{P\del[1]{T_{n,W} - \sqrt{n}D_n^{-1}\mu_n\in H'}-P\del[1]{\Sigma_{0,n}^{1/2}Z_d\in H'}}
\leq a_n.
\end{equation}
Recall that~$C_1$ (depending only on~$b_1,b_2$,~$\lambda_1',\lambda_2'$, and~$m$) denotes the constant ``$C$'' from Lemma~\ref{lem:varwinsestim}, that~$\mathfrak{d}(n, d, m, \overline{\eta}_n) < b_1^2/(2C_1)$, and denote $C'_1 := \sqrt{2}C_1/b_1^3$. Lemma~\ref{lem:varwinsestim} shows that with probability at least~$1-24/n$ it holds that
\begin{equation*}
\min_{j = 1, \hdots, d} \tilde{\sigma}_{n,j} \geq \frac{b_1}{\sqrt{2}} \quad \text{ and } \quad  \max_{j=1,\hdots,d}\envert[1]{\tilde{\sigma}_{n,j}^{-1}-\sigma_{2,n,j}^{-1}(P)}\leq C_1'\cdot \mathfrak{d}(n,d, m, \overline{\eta}_n)
\end{equation*}

Next, write~$\sigma_{n,j}=\sigma_{2,n,j}$, note that~$S_{n,W}
=
T_{n,W}+\sqrt{n}(\tilde{D}_n^{-1}-D_n^{-1})\mu_n$ (recalling the definition of~$\tilde{D}_n^{-1}$ right after Theorem~\ref{thm:41kp}), and write 
\begin{equation*}
\enVert[1]{S_{n,W}-T_{n,W}}_{\infty}
=
\enVert[1]{\sqrt{n}(\tilde{D}_n^{-1}-D_n^{-1})\mu_n}_{\infty}
\leq
\sqrt{n}\|\tilde{D}_n^{-1}-D_n^{-1}\|_{\infty}\|\mu_n\|_\infty.
\end{equation*}
Thus, abbreviating~$\overline{A}_n:=C_1'\sqrt{n}\mathfrak{d}(n, d, m, \overline{\eta}_n)\|\mu_n\|_\infty$, 
it holds that~
\begin{equation}\label{eq:boundSnW}
P_n\del[1]{\enVert[1]{S_{n,W}-T_{n,W}}_{\infty}> \overline{A}_n}\leq 24/n.
\end{equation}
Using~\eqref{eq:boundTnW} and~\eqref{eq:boundSnW}, we can now apply Lemma~\ref{lem:auxlem} with~$H = H^*$,~$U = T_{n,W} - \sqrt{n}D_n^{-1} \mu_n$, $V = S_{n,W} - \sqrt{n}D_n^{-1} \mu_n$, ~$\Xi = \Sigma_{0, n}$,~$\xi = 1$,~$a = a_n$, $b = \overline{A}_n$ and~$c = 24/n$ to conclude 
\begin{align*}
a_n + \left(\sqrt{2 \log(d)} + 4\right)\overline{A}_n + 24/n  &\geq \left| P(S_{n,W} - \sqrt{n}D_n^{-1} \mu_n \in H^*) - P(\Sigma_{0, n}^{1/2} Z_d \in H^*) \right| \\
&= \left| P(S_{n,W} \in H) - P(\Sigma_{0, n}^{1/2} Z_d + \sqrt{n}D_n^{-1}\mu_n \in H) \right|.
\end{align*}
To obtain~\eqref{eqn:SUMGASW}, note that there exists a constant~$C$, say, depending only on~$b_1,b_2, \lambda_1,\lambda_2$, $\lambda_1',\lambda_2'$, and~$m$, such that the upper bound~$a_n + (\sqrt{2 \log(d)} + 4)\overline{A}_n + 24/n$ is dominated by
\begin{equation}
C \cdot \left[\mathfrak{e}(n, d, m, \overline{\eta}_n) + \mathfrak{d}(n, d, m, \overline{\eta}_n)\sqrt{n \log(d)}\|\mu_n\|_{\infty}\right] + \frac{24}{n}.
\end{equation}

To conclude, we recall from Theorem 4.1 in \cite{kock2025high} (in the version stated in Theorem~\ref{thm:41kp} above) that under the conditions imposed in~\eqref{eqn:divndeta} it holds that~$\mathfrak{e}(n', d', m, \overline{\eta}_{n'}) \to 0$. We finish  by using~\eqref{eq:extrasslocgauss}, also noting that~$\eps_n,\eps_n'\in(0,1/2)$ is satisfied eventually.
\end{proof}

\subsection{Asymptotic behavior of~$S_{n,W}$ when~$\enVert[0]{\mu_n(P)}_\infty$ is not ``small''}

We consider the asymptotic behavior of~$S_{n,W}$ when~\eqref{eq:extrasslocgauss} is not satisfied.

\begin{theorem}\label{thm:SumGAWinsLarge}
Let~$b_1\in(0,\infty)$,~$b_2\in(b_1,\infty)$, and~$m > 2$.  Fix~$\lambda_1,\lambda_1'\in(1,\infty)$, and~$\lambda_2,\lambda_2'\in(0,\infty)$. Let a sequence~$P_{n} \in \mc{P}(b_1,b_2,m)$ and a subsequence~$n'$ of~$n$ be given, such that~$d' \geq 2$,~\eqref{eqn:divndeta} holds, and 
\begin{equation}\label{eq:extrasslocgaussLarge}
\sqrt{n'\log(d')} \mathfrak{d}(n',d',m,\overline{\eta}_{n'})\|\mu_{n'}(P_{n'})\|_\infty	
\to c \in (0, \infty].
\end{equation}
Then, for every sequence~$\mathcal{H}^*_n$ satisfying $\emptyset \neq \mathcal{H}^*_n \subseteq \mathcal{H}_n$ along~$n'$ and such that $$r_{n'} := \sup \left \{\|x\|_{\infty}: x \in H,~ H \in \mathcal{H}^*_{n'}\right \}= O(\sqrt{\log(d')}),$$ it holds that
\begin{equation}\label{eqn:SUMGA2W}
\sup_{H\in\mc{H}^*_{n'} }\envert[1]{P_{n'}\del[1]{S_{n', W}\in H}-P_{n'}\del[1]{\Sigma_{0,n'}^{1/2}(P_{n'})Z_{d'}+\sqrt{n'}D_{n'}^{-1}\mu_{n'}(P_{n'})\in H}} \to 0;
\end{equation}
more specifically, for every sequence~$H_n$ such that~$H_n \in \mathcal{H}^*_n$ along~$n'$, we then have 
\begin{equation}\label{eqn:SUMGA3W}
P_{n'}\del[1]{S_{n', W}\in H_{n'}} \to 0 \quad \text{and} \quad P_{n'}\del[1]{\Sigma_{0,n'}^{1/2}(P_{n'})Z_{d'}+\sqrt{n'}D_{n'}^{-1}\mu_{n'}(P_{n'})\in H_{n'}} \to 0.
\end{equation}
\end{theorem}
\begin{remark}\label{rem:wecond}
Inspection of the proof of Theorem \ref{thm:SumGAWinsLarge} shows that if instead of~\eqref{eq:extrasslocgaussLarge} one imposes~$\frac{\sqrt{\log(d')}}{\sqrt{n'}\|\mu_{n'}(P_{n'})\|_\infty}\to 0$, the so-obtained version of Theorem \ref{thm:SumGAWinsLarge} is a correct statement. 
\end{remark}
\begin{proof}
It suffices to establish the statement in~\eqref{eqn:SUMGA3W}. Write~$n$ instead of~$n'$, and drop the dependence of several quantities on~$P_n$ when possible. Let~$S_{n,W,c}\in\R^d$ be defined as in~\eqref{eqn:snwcdef} (note that eventually, and hence without loss of generality, it holds that~$\eps_n,\eps_n'\in(0,1/2)$, with~$\eps_n$ as in~\eqref{eq:epsfam} and~$\eps_n'$ as in~\eqref{eq:epsprime}). We can use the same argument as in the beginning of the proof of Theorem~\ref{thm:SumGAnonloc}, but now replacing $W_n$ by $S_{n,W,c}$ and replacing the reference to Lemma~\ref{lem:GAStanfardizedSums} by a reference to Theorem 4.1 in~\cite{kock2025high} (in the version stated in Theorem~\ref{thm:41kp} above), to show that
\begin{equation*}
P_n\del[1]{S_{n,W} \notin H_n}  \geq P_n\left(\|\sqrt{n}D_n^{-1}\mu_n+\sqrt{n}(\tilde{D}_n^{-1}-D_n^{-1})\mu_n\|_\infty > r_n + \|S_{n,W,c}\|_\infty\right).
\end{equation*}
and that
\begin{equation}\label{eqn:rhslognW}
r_n + \|S_{n,W,c}\|_\infty = O_{P_n}\del[1]{\sqrt{\log(d)}} \quad \text{ and } \quad r_n + \|\Sigma_{0,n}^{1/2}Z_d\|_{\infty} = O_{P_n}\del[1]{\sqrt{\log(d)}}.
\end{equation}
Note that the conditions~\eqref{eq:epscond'} and $\mathfrak{d}(n, d, m, \overline{\eta}_n) <\frac{b_1^2}{2C}$, imposed in Lemma~\ref{lem:varwinsestim}, hold eventually: First,  $\mathfrak{d}(n, d, m, \overline{\eta}_n)$, which was defined in~\eqref{eqn:ddef}, converges to~$0$ under~\eqref{eqn:divndeta}. Second,~$\eps_n'$ also converges to 0 under~\eqref{eqn:divndeta}. Thus, these conditions can be assumed to hold without loss of generality, and it follows from Lemma~\ref{lem:varwinsestim} that for a constant~$C'=\sqrt{2}C/b_1^3$ depending only on~$b_1, b_2, \lambda_1'$, $\lambda_2'$ and~$m$ 
\begin{equation*}
P_n\left(\min_{j = 1, \hdots, d} \tilde{\sigma}_{n,j} \geq \frac{b_1}{\sqrt{2}} \text{ and } \enVert[0]{\tilde{D}_n^{-1}-D_n^{-1}}_\infty \leq C'\mathfrak{d}(n, d, m, \overline{\eta}_n)\right)
\geq 1-\frac{24}{n}.
\end{equation*}
Arguing as around~\eqref{eqn:mulow*} one obtains that (eventually) with probability at least~$1-24/n$ 
\begin{equation*}
\enVert[2]{\sqrt{n}D_n^{-1}\mu_n+\sqrt{n}(\tilde{D}_n^{-1}-D_n^{-1})\mu_n}_\infty
\geq \sqrt{n}\enVert[0]{\mu_n}_\infty/(2b_2).
\end{equation*}
Hence, together with~\eqref{eqn:rhslognW} and~Lemma~\ref{lem:auxL} below, we obtain (note that eventually~$\mu_n \neq 0$)
\begin{equation*}
P_n\del[1]{S_{n,W} \notin H_n} \geq P_n\left((2b_2)^{-1} > [r_n + \|S_{n,W,c}\|_\infty ]/(\sqrt{n}\enVert[0]{\mu_n}_\infty)\right) + o(1) \to 1,
\end{equation*}
Concerning the second convergence statement in~\eqref{eqn:SUMGA3W} it remains to note that (eventually)
\begin{align*}
P_{n}\del[1]{\Sigma_{0,n}^{1/2}Z_{d}+\sqrt{n}D_{n}^{-1}\mu_{n}\notin H_{n}} &\geq  P_n\del[1]{\enVert[0]{\Sigma_{0,n}^{1/2}Z_d+\sqrt{n}D_{n}^{-1}\mu_{n}}_\infty> r_n}\\
&\geq
P_n\del[1]{\sqrt{n} \enVert[0]{D_{n}^{-1}\mu_{n}}_\infty>r_n + \enVert[0]{\Sigma_{0,n}^{1/2}Z_d}_\infty} \\
&\geq P_n(b_2^{-1} \geq   [r_n + \enVert[0]{\Sigma_{0,n}^{1/2}Z_d}_\infty]/(\sqrt{n}\enVert[0]{\mu_n}_\infty)) \to 1,
\end{align*}
where we used~\eqref{eqn:mulow*}, Lemma~\ref{lem:auxL}, and~\eqref{eqn:rhslognW}. 
\end{proof}

For convenience, we recall a limit from Lemma B.4 in~\cite{kock2025high}.
\begin{lemma}[Equation~(B.9) of Lemma B.4 in~\cite{kock2025high}]\label{lem:covcompare}
Let~$m>2$ and assume that~$\sqrt{n\log(d)}\overline{\eta}_n^{1-\frac{1}{m}}\to 0$ as well as~$\log(d)/n^{\frac{m-2}{5m-2}}\to 0$. Then,
\begin{equation}\label{eq:conv3}
\log^2(d) \cdot \mathfrak{d}(n,d, m, \overline{\eta}_n) = \log^2(d)\sbr[3]{\overline{\eta}_n^{1-\frac{2}{m}}+\del[2]{\frac{\log(dn)}{n}}^{1-\frac{1}{(m/2)\wedge 2}}}  \to 0.
\end{equation}
\end{lemma}

\begin{lemma}\label{lem:auxL}
Under the conditions imposed in Theorem~\ref{thm:SumGAWinsLarge} we have~$\frac{\sqrt{\log(d')}}{\sqrt{n'}\|\mu_{n'}(P_{n'})\|_\infty}\to 0$.
\end{lemma}

\begin{proof}
Suppose first that~$c\in(0,\infty)$ in~\eqref{eq:extrasslocgaussLarge}. Write~$n$ and~$d$ instead of~$n'$ and~$d'$, abbreviate~$\mu_n(P_n)$ by~$\mu_n$, assume without loss of generality that~$\mu_n \neq 0$, and note (recalling~\eqref{eqn:ddef}) that condition~\eqref{eq:extrasslocgaussLarge} implies that eventually one of the following two conditions holds:
\begin{equation*}
(i) ~~\sqrt{n\log(d)}\overline{\eta}_{n}^{1-\frac{2}{m}}\|\mu_{n}\|_\infty	\geq c/4 =: C \quad \Rightarrow \quad \frac{\sqrt{\log(d)}}{\sqrt{n}\|\mu_n\|_\infty}
\leq
C^{-1}\log(d)\overline{\eta}_{n}^{1-\frac{2}{m}} \to 0,
\end{equation*}
the convergence following from (B.9) of Lemma B.4 in~\cite{kock2025high} (as restated in Lemma~\ref{lem:covcompare}); or
\begin{align*}
(ii) ~~\sqrt{n\log(d)}\del[2]{\frac{\log(dn)}{n}}^{1-\frac{1}{(m/2)\wedge 2}}\|\mu_{n}\|_\infty
\geq  C \quad \Rightarrow\\
\frac{\sqrt{\log(d)}}{\sqrt{n}\|\mu_n\|_\infty}
\leq
\begin{cases}
C^{-1}\frac{\log(dn)^{2-\frac{2}{m}}}{n^{1-\frac{2}{m}}} \to 0	\qquad &\text{if }m\in(2,4], \\
C^{-1}\frac{\log(dn)^{\frac{3}{2}}}{n^{1/2}} \to 0	\qquad &\text{if }m\in(4,\infty),
\end{cases}
\end{align*}
the convergence following from~\eqref{eqn:divndeta}. Finally, if~$c=\infty$ in \eqref{eq:extrasslocgaussLarge} then we note that for, e.g.,~$C = 1$ one of (i) or (ii) above must hold eventually, again giving the desired conclusion.
\end{proof}

\subsection{Proof of Theorem~\ref{thm:Wins}}
The proof strategy is identical to that of Theorem~\ref{thm:sums}. The proof carries over almost verbatim. We only have to replace~$S_n$ by~$S_{n,W}$ and the references to Theorems~\ref{thm:SumGA} and~\ref{thm:SumGAnonloc} by references to Theorems~\ref{thm:SumGAWins} and~\ref{thm:SumGAWinsLarge}, respectively. We provide the argument for completeness:

To establish the first part, let~$P_n$ be a sequence in~$\mc{P}^0(b_1,b_2,m)$, and suppress the dependence of all quantities on~$P_n$. Theorem~\ref{thm:SumGAWins} gives (note that~$\mu_n = 0$ implies~\eqref{eq:extrasslocgauss})
\begin{equation*}
\envert[1]{P_{n}\del[1]{\enVert[0]{S_{n,W}}_\infty>c_{n,1-\alpha}}-P_n\del[1]{\enVert[0]{\Sigma_{0,n}^{1/2}Z_d}_\infty>c_{n,1-\alpha}}}\to 0.	
\end{equation*}
By the Khatri-{\v{S}}id{\'a}k inequality (the inequality being an equality in case~$\Sigma_{0,n} = I_d$)
\begin{equation*}
P_n\del[1]{\enVert[0]{\Sigma_{0,n}^{1/2}Z_d}_\infty>c_{n,1-\alpha}}
\leq 
P_n\del[1]{\enVert[0]{Z_d}_\infty>c_{n,1-\alpha}} = \alpha,
\end{equation*}
which establishes the first two statements in Part~1 (respectively choosing~$P_n$ to approximate the supremum in~\eqref{eq:sizeWins}, or~$P_n$ satisfying~$\Sigma_{0, n}(P_n) = I_d$ for~$n$ large enough). The last statement follows from Lemma~\ref{lem:extreme}, which shows that if~$d \to \infty$, then for every sequence~$P_n$ in~$\mathcal{P}(b_1, b_2, m, \bm{r})$ we have~$P_n(\|\Sigma_{0,n}^{1/2}Z_d\|_\infty>c_{n,1-\alpha}) \to \alpha$ under the assumption that~$\log(l)r_l \to 0$ as $l\to \infty$.

To establish~\eqref{eq:GaussApproxWins}, we use the equivalence between Parts 1 and 3 in Lemma~\ref{lem:unifseq}. Let~$P_n$ be a sequence in~$\mc{P}(b_1,b_2,m)$, and let~$n'$ be a subsequence of~$n$. If~$$\mathfrak{d}(n',d',m,\overline{\eta}_{n'})\sqrt{n'\log(d')}\|\mu_{n'}\|_{\infty} \to 0,$$ then Theorem~\ref{thm:SumGAWins} shows that
\begin{equation}\label{eqn:locconvWI}
\envert[2]{P_{n'}\del[1]{\enVert[0]{S_{n',W}}_\infty>c_{n',1-\alpha}}-P_{n'}\del[1]{\enVert[0]{\Sigma_{0,n'}^{1/2}Z+\sqrt{n'}D_{n'}^{-1}\mu_{n'}}_\infty>c_{n',1-\alpha}}}
\to 0.
\end{equation}
Otherwise,~$n'$ possesses a subsequence~$n''$ satisfying~$$\mathfrak{d}(n'',d'',m,\overline{\eta}_{n''})\sqrt{n''\log(d'')}\|\mu_{n''}\|_{\infty} \to c \in (0, \infty].$$ Theorem~\ref{thm:SumGAWinsLarge} together with~$c_{n,1-\alpha} = O(\sqrt{\log(d)})$ then shows~\eqref{eqn:locconvWI} along~$n''$.

\subsection{Auxiliary results for Section~\ref{sec:necess} and the proof of Theorem~\ref{thm:EtaNecessaryshort}}\label{sec:AppLB}
\emph{Throughout this section, fix~$m > 2$, let~$0 < b_1 < b_2 < \infty$, and recall from Section~\ref{sec:necess} that~$a =(b_2-b_1)/2$.} Dependence on those quantities will not be explicitly highlighted whenever this causes no confusion. We first establish some properties of~$Q_{p,d}$, cf.~\eqref{eqn:defqpd} for its definition. In this context only, to lighten notation, for a distribution~$Q$ on~$\R^d$, we write~$\mu(Q)$ and~$\Sigma(Q)$ for its mean and covariance matrix,~$D(Q)$ for the diagonal matrix of coordinate standard deviations, and~$\Sigma_0(Q)$ for its correlation matrix~$D(Q)^{-1}\Sigma(Q)D(Q)^{-1}$ (assuming that the respective quantity exists).
\begin{lemma}\label{lem:HDDistConstruct}
Fix~$d\geq2$ and~$p\in[0, 1]$. Then,
\begin{align}
\mu(Q_{p,d})&=ap^{1-1/m}\bm{\iota}_d,
& D(Q_{p,d})&=sI_d,
& \Sigma_0(Q_{p,d})&=(1-\rho)I_d+\rho\bm{\iota}_d\bm{\iota}_d',
\label{eq:EtaQ1Moments}
\end{align}
where~$s$ is the positive square root of~$s^2$,
\begin{equation*}
s^2:=b_1^2+a^2p^{1-2/m}(1-p),
\qquad \text{ and } \qquad
\rho:=\frac{a^2p^{1-2/m}(1-p)}{s^2}.
\end{equation*}
Under~$Q_{p,d}$ each coordinate projection mapping~$v \mapsto v_j$,~$v \in \R^d$ has a finite absolute~$m$-th moment, a variance bounded from below by~$b_1^2$, and its centered~$m$-th absolute moment bounded from above by~$b_2^m$.
\end{lemma}

\begin{proof}
Note that~$\E(V)=0_d$ and
$\operatorname{Cov}(V)=I_d$. This establishes~\eqref{eq:EtaQ1Moments} in case~$p = 0$. In case~$p> 0$, write~$h := ap^{-1/m}$ and note that~$$\E(h B\bm{\iota}_d) = ap^{1-1/m}\bm{\iota}_d \quad \text{ and } \quad \operatorname{Cov}(h B\bm{\iota}_d) = p(1-p) h^2 \bm{\iota}_d \bm{\iota}_d'.$$ Using that~$V$ is independent of~$B$, we obtain~\eqref{eq:EtaQ1Moments} also in that case. The remaining statement concerning the coordinate projection mappings is trivial in case~$p = 0$. Hence, again, we can focus on the case~$p > 0$. Fix a coordinate~$j \in \{1, \hdots, d\}$. Clearly,~$b_1 V_j + hB$ possesses absolute moments of any order. Furthermore, we have already established in~\eqref{eq:EtaQ1Moments} that its variance equals~$s^2 \geq b_1^2$. We conclude by using Minkowski's inequality and~$\E|B-p|^m = (1-p)p^m+p(1-p)^m \leq 2p$ to obtain
\begin{equation*}
\del[1]{\E\envert[1]{b_1V_j+h(B-p)}^m}^{1/m}
\leq
b_1+h\del[1]{E|B-p|^m}^{1/m}
\leq b_1+2^{1/m}a
\leq b_2.
\end{equation*}
\end{proof}
For~$p \in [0, 1]$, denote 
\begin{equation*}
\mc{S}_d(p) := 
\begin{cases}
b_1 \cdot \{-1, 1\}^d & \text{ if } p = 0, \\
b_1 \cdot \{-1, 1\}^d + ap^{-1/m}\bm{\iota}_d & \text{ if } p > 0.
\end{cases}
\end{equation*}
Observe that for~$p \in (0, 1]$ it holds that
\begin{equation}\label{eqn:suppsep}
\mc{S}_d(p) \cap \mc{S}_d(0) =
\begin{cases}
\emptyset & \text{ if } ap^{-1/m} \neq 2b_1, \\
\{b_1 \bm{\iota}_d\} & \text{ if } ap^{-1/m} = 2b_1.
\end{cases}
\end{equation}
For every
$x=(x_1,\hdots,x_n)\in(\R^d)^n$ and~$i \in \{1, \hdots, n\}$, we write~$N_i(x):=\sum_{l=1}^i\mathds{1}_{\mc{S}_d(p)}(x_l)$, counting how many of the~$x_1, \hdots, x_i$ are contained in~$\mathcal{S}_d(p)$. For~$\eta\in[0,1]$, we now formally define the contamination scheme~$\mathsf{C}_{n,p,\eta}: (\R^d)^n \to (\R^d)^n$ via
$\mathsf{C}_{n,p,\eta}(x)=(\tilde x_1,\hdots,\tilde x_n)$ for
\begin{equation}\label{eq:EtaContamMap}
\tilde x_i
:=
x_i-ap^{-1/m}\bm{\iota}_d
\mathds{1}_{\mc{S}_d(p)}(x_i)
\mathds{1}_{\cbr[0]{N_i\leq\eta n}}(x),
\qquad i=1,\hdots,n,\quad \text{and}\quad p\in(0,1].
\end{equation}
By definition,~$\mathsf{C}_{n,p,\eta}\in\mathfrak{C}_n(\eta)$; cf.~the notation introduced in Section~\ref{sec:necess}. Intuitively, cf.~\eqref{eqn:defqpd}, if one samples i.i.d.~random vectors from~$Q_{p,d}$ ($p > 0$), the sample corrupted via~$\mathsf{C}_{n,p,\eta}$ then very likely ``looks like'' a sample from~$Q_{0,d}$. This is formalized in Lemma~\ref{lem:HDTVDist} by bounding the total variation distance between~$Q_{0, d}^{\otimes n}\circ\mathsf{C}_{n,p,\eta}^{-1}$ and~$Q_{p,d}^{\otimes n}\circ\mathsf{C}_{n,p,\eta}^{-1}
$ from above. We exploit the ``coupling'' of~$Q_{0, d}$ and~$Q_{p, d}$ used in their construction in~\eqref{eqn:defqpd}. 

As mentioned in the main text, a similar argument is used in the proof of Lemma 2.2 in the updated (second) version of \cite{resende2024robust}. The total variation distance between two probability measures~$R_0$ and~$R_1$ on the same measurable space is denoted by
\begin{equation*}
\mathsf{TV}(R_0,R_1)
=
\sup_A\envert[1]{R_0(A)-R_1(A)},
\end{equation*}
the supremum being taken over all measurable sets~$A$.
\begin{lemma}\label{lem:HDTVDist}
Let~$n\in\N$,~$d \in \N$,~$0 < 2p \leq \eta \leq 1$, and suppose
that~$ap^{-1/m}\neq2b_1$. Then, 
\begin{equation}\label{eq:HDTVDist}
\mathsf{TV}\del[1]{
Q_{0, d}^{\otimes n}\circ\mathsf{C}_{n,p,\eta}^{-1},
Q_{p,d}^{\otimes n}\circ\mathsf{C}_{n,p,\eta}^{-1}
}
\leq
\del[2]{\frac{e}{4}}^{np}.
\end{equation}
\end{lemma}

\begin{proof}
Let~$(V_i,B_i)$,
$i=1,\hdots,n$, be independent copies of~$(V,B)$ used in Equation~\eqref{eqn:defqpd}, and set
\begin{equation*}
X_i^0=b_1V_i,
\qquad
X_i^p=b_1V_i+ap^{-1/m}B_i\bm{\iota}_d.
\end{equation*}
Then~$\mathcal{L} (X_1^j,\hdots,X_n^j)=Q_{j,d}^{\otimes n}$ for every
$j\in\{0,p\}$. From~$ap^{-1/m}\neq2b_1$ and~\eqref{eqn:suppsep}, we obtain that~$\mc{S}_d(0)$, the support of the~$X_i^0$, is disjoint from~$\mc{S}_d(p)$. It hence follows immediately that~$$\mathsf{C}_{n,p,\eta}(X_1^0,\hdots,X_n^0) = (X_1^0,\hdots,X_n^0).$$ Furthermore, note that~$X_i^p \in \mc{S}_d(p)$ if and only if~$B_i = 1$. Note also that~$X_i^p = X_i^0$ if~$B_i = 0$. If $B_i = 1$ holds for at most~$\eta n$ observations among the~$X_1^p, \hdots, X_n^p$, the contamination scheme~$\mathsf{C}_{n,p,\eta}$ subtracts~$ap^{-1/m} \bm{\iota}_d$ from \emph{all} of those observations, turning each of them (and hence all the~$X_i^p$) into~$X_i^0$. That is,
\begin{equation}\label{eqn:equalif}
\mathsf{C}_{n,p,\eta}(X_1^p,\hdots,X_n^p) = (X_1^0,\hdots,X_n^0) = \mathsf{C}_{n,p,\eta}(X_1^0,\hdots,X_n^0) \quad \text{ if } \quad \sum_{i=1}^nB_i\leq\eta n.
\end{equation}
For every measurable $A\subseteq(\R^d)^n$, we therefore obtain
\begin{align*}
&\envert[1]{
Q_{0, d}^{\otimes n}\circ\mathsf{C}_{n,p,\eta}^{-1}(A) - 
Q_{p,d}^{\otimes n}\circ\mathsf{C}_{n,p,\eta}^{-1}(A)
} \\
&=
\left|\E\left[\mathds{1}_A(\mathsf{C}_{n,p,\eta}(X_1^0,\hdots,X_n^0)) - \mathds{1}_A(\mathsf{C}_{n,p,\eta}(X_1^p,\hdots,X_n^p)\right] \right|\\
&\leq
\E\left|\mathds{1}_A(\mathsf{C}_{n,p,\eta}(X_1^0,\hdots,X_n^0)) - \mathds{1}_A(\mathsf{C}_{n,p,\eta}(X_1^p,\hdots,X_n^p))\right| \\
&\leq
\P\left(\mathsf{C}_{n,p,\eta}(X_1^0,\hdots,X_n^0) \neq \mathsf{C}_{n,p,\eta}(X_1^p,\hdots,X_n^p)\right).
\end{align*}
By~\eqref{eqn:equalif} the last probability is at most
\begin{equation*}
\P\del[3]{\sum_{i=1}^nB_i>\eta n}
\leq
\P\del[3]{\sum_{i=1}^nB_i>2np}
\leq
\del[2]{\frac{e}{4}}^{np},
\end{equation*}
where we used that~$2p \leq \eta$, ~$\sum_{i=1}^nB_i\sim\operatorname{Bin}(n,p)$, and applied~(5) in
\cite{hagerup1990guided}, with~$\epsilon$ there equal to~$1$ and~$m$
there equal to~$np$. As~$A$ was arbitrary,
\eqref{eq:HDTVDist} follows.
\end{proof}
We now prove the following extended version of Theorem~\ref{thm:EtaNecessaryshort}. Recall the notation introduced in~\eqref{eqn:qkpdef}.

\begin{theorem}\label{thm:EtaNecessary}
Let~$b_1\in(0,\infty)$,~$b_2\in(b_1,\infty)$,~$m>2$, and~$\alpha\in(0,1)$. Suppose Assumption~\ref{ass:lbd} is satisfied,~$d\to\infty$, and~$\log(d)/n^{\frac{m-2}{5m-2}}\to0$. Define the infinite set~$$\N(c) := \{n \in \N: \sqrt{n\log(d \vee 2)}\overline{\eta}_n^{1-1/m} \geq c\} \subseteq \N.$$ The following holds:
\begin{enumerate}
\item Denoting~$\tau:=a\kappa^{1-1/m}/b_1$ and~$\theta := 1-(1-\alpha)^{\cosh(\sqrt{2}\tau)}>\alpha$, we have for every~$n > M$ that
\begin{equation}\label{eq:EtaBoundary1}
\mu_n(P_n^1)
=
\frac{b_1 \tau}{\sqrt{n\log(d \vee 2)}}\bm{\iota}_d,
\qquad
\mc{G}_n(P_n^0)=\alpha;
\qquad\text{and}\qquad 
\mc{G}_n(P_n^1) \to \theta.
\end{equation}
\item %
Setting~$\mathsf{C}_n := \mathsf{C}_{n,p_n,\overline{\eta}_n} \in\mathfrak{C}_n(\overline{\eta}_n)$ for every~$n > M$, it holds that
\begin{equation}\label{eq:EtaBoundary2}
\lim_{\substack{n\to\infty\\ n\in \N(c)}} \mathsf{TV}\left(
P_n^0 \circ [\mathsf{C}_n(X_{1,n}, \hdots, X_{n,n})]^{-1},
P_n^1\circ [\mathsf{C}_n(X_{1,n}, \hdots, X_{n,n})]^{-1}
\right) =
0.
\end{equation}
\item For every sequence of tests~$\psi_n:(\R^d)^n\to[0,1]$, the property 
\begin{equation}\label{eqn:asysize0}
\limsup_{\substack{n\to\infty\\ n\in\N(c)}}
E_{P_n^0}\sbr[1]{\psi_n(\mathsf{C}_n(X_{1,n}, \hdots, X_{n,n}))}
\leq\alpha
\end{equation}
implies 
\begin{equation}\label{eqn:asypow00}
\limsup_{\substack{n\to\infty\\ n\in\N(c)}}
E_{P_n^1}\sbr[1]{\psi_n(\mathsf{C}_n(X_{1,n}, \hdots, X_{n,n}))}
\leq\alpha,
\end{equation}
and thus
\begin{equation}\label{eqn:asypow0}
\liminf_{\substack{n\to\infty\\ n\in\N(c)}} \left[\mc{G}_n(P_n^1) - 
E_{P_n^1}\sbr[1]{\psi_n(\mathsf{C}_n(X_{1,n}, \hdots, X_{n,n}))} \right]
\geq \theta - \alpha > 0.
\end{equation}
\end{enumerate}
\end{theorem}

\begin{proof}[Proof of Theorem~\ref{thm:EtaNecessary}]
We start with Part 1 and assume that~$n > M$: Assumption~\ref{ass:lbd} and Lemma~\ref{lem:HDDistConstruct} show that~$\mu_n(P_n^0)=0$ and~$\Sigma_{0,n}(P_n^0)=I_d$. Hence, it follows from the definition of~$c_{n, 1-\alpha}$, cf.~the discussion preceding~\eqref{eq:CVDiag}, that
\begin{equation*}
\mc{G}_n(P_n^0) = P_n^0\left(\enVert[0]{Z_d}_\infty>c_{n,1-\alpha}\right) =\alpha.
\end{equation*} 
Similarly, by Assumption~\ref{ass:lbd} and Lemma~\ref{lem:HDDistConstruct},~$\mu_n(P_n^1)=ap_n^{1-1/m}\bm{\iota}_d$, which, recalling~\eqref{eqn:qkpdef}, establishes the first equality in~\eqref{eq:EtaBoundary1}. To finish the proof of~\eqref{eq:EtaBoundary1}, it remains to verify~$\mc{G}_n(P_n^1) \to \theta > \alpha$. While the inequality is trivial ($\tau > 0$,  hence~$\cosh(\sqrt{2}\tau)>1$), the convergence is more involved. Let~$s_n$ and~$\rho_n$ be the quantities~$s$ and~$\rho$ in
Lemma~\ref{lem:HDDistConstruct} evaluated at~$p=p_n$, i.e., 
\begin{equation*}
s_n^2 =b_1^2+a^2p_n^{1-2/m}(1-p_n)
\qquad \text{ and } \qquad
\rho_n= a^2p_n^{1-2/m}(1-p_n)s_n^{-2}.
\end{equation*}
Hence, by definition (and Assumption~\ref{ass:lbd} and Lemma~\ref{lem:HDDistConstruct}), we can write
\begin{equation*}
\mc{G}_n(P_n^1) = P_n^1\left(\enVert[0]{\Sigma_{0,n}^{1/2}(P_n^1)Z_d+ t_n \bm{\iota}_d}_\infty>c_{n,1-\alpha}\right), \quad \text{ where } \quad	t_n:=\sqrt{n}ap_n^{1-1/m}s^{-1}_n.
\end{equation*} 
We next determine some auxiliary limits: the assumptions~$\log(d)/n^{\frac{m-2}{5m-2}}\to0$ and~$d\to \infty$ deliver, after recalling that~$p_n = \kappa (n\log(d \vee 2))^{-\frac{m}{2(m-1)}}$,
\begin{align}
np_n
&=
\kappa\del[2]{\frac{n^{m-2}}{\log^m(d \vee 2)}}^{\frac{1}{2(m-1)}} 
\to\infty, \text{ and}
\label{eq:EtaRates0} \\
\log^2(d)p_n^{1-2/m}
&=
\kappa^{1-2/m}
\del[2]{\frac{\log^{3m-2}(d \vee 2)}{n^{m-2}}}^{\frac{1}{2(m-1)}}
\to0;
\label{eq:EtaRates}
\end{align}
moreover, by~$p_n \to 0$,~$0 \leq \rho_n\leq
a^2p_n^{1-2/m}/b_1^2$, and~\eqref{eq:EtaRates}
\begin{equation}\label{eq:EtaTandRho}
s_n\to b_1,
\qquad
t_n\sqrt{\log(d)}\to \frac{a\kappa^{1-1/m}}{b_1} = \tau,
\qquad\text{and}\qquad
\log^2(d)\rho_n\to0.
\end{equation}
Now, because~$\{x \in \R^d:\|x+t_n\bm{\iota}_d\|_\infty\leq c_{n,1-\alpha}\}$
is a hyperrectangle, and by Assumption~\ref{ass:lbd} and Lemma~\ref{lem:HDDistConstruct} 
$$\log^2(d)\|\Sigma_{0,n}(P_n^1)-I_d\|_{\max}
=\log^2(d)\rho_n\to0,$$ the
Gaussian-to-Gaussian comparison inequality in Proposition~2.1 of
\cite{chernozhuokov2022improved} (cf.~Proposition~2 in
\cite{chernozhukov2023high}) gives
\begin{equation}\label{eq:EtaGaussianCompare}
\envert[2]{
\mc{G}_n(P_n^1)
-
P_n^1\del[1]{\enVert[0]{Z_d+t_n\bm{\iota}_d}_\infty>c_{n,1-\alpha}}
}
\to0.
\end{equation}
Thus it suffices to study~$P_n^1\del[1]{\enVert[0]{Z_d+t_n\bm{\iota}_d}_\infty>c_{n,1-\alpha}}$. With~$\Phi$ denoting the cdf of the standard normal distribution and~$\overline{\Phi} := 1- \Phi$, we can write
\begin{equation*}
P_n^1\del[1]{\enVert[0]{Z_d+t_n\bm{\iota}_d}_\infty\leq c_{n,1-\alpha}}
=
\sbr[1]{
1-\overline{\Phi}(c_{n,1-\alpha}-t_n)-\overline{\Phi}(c_{n,1-\alpha}+t_n)
}^d =: \del[1]{1-v_n}^d.
\end{equation*}
Conclude from the definition of~$c_{n,1-\alpha}$, noting that~$1-\alpha = (1-2\overline{\Phi}(c_{n,1-\alpha}))^d$, and from~$d \to \infty$,~\eqref{eq:CVDiag} and~\eqref{eq:EtaTandRho} that
\begin{equation}\label{eqn:auxmill}
2d\overline{\Phi}(c_{n,1-\alpha})\to - \log(1-\alpha),
\qquad
\frac{c_{n,1-\alpha}}{\sqrt{\log(d)}}\to\sqrt{2}, \qquad c_{n,1-\alpha}t_n\to \sqrt{{2}}\tau.
\end{equation}
This, in turn, shows that~$c_{n,1-\alpha} \to \infty$ as well as~$c_{n,1-\alpha}\pm t_n\to\infty$, and hence delivers via a Mills' ratio argument ($\phi$ denoting the standard normal density, and~``$a_n\sim b_n$'' denoting $a_n/b_n \to 1$ for sequences~$a_n$ and~$b_n$, say, in this display) that
\begin{equation*}
\frac{\overline{\Phi}(c_{n,1-\alpha} \pm t_n)}{\overline{\Phi}(c_{n,1-\alpha})} \sim \frac{\phi(c_{n,1-\alpha} \pm t_n)}{\phi(c_{n,1-\alpha})} \frac{c_{n,1-\alpha}}{c_{n,1-\alpha} \pm t_n} \sim e^{-(c_{n,1-\alpha} \pm t_n)^2/2 + c^2_{n,1-\alpha}/2}
\to e^{\mp \sqrt{2}\tau}.
\end{equation*}
Consequently,~$$\frac{v_n}{u_n} := \frac{\overline{\Phi}(c_{n,1-\alpha}-t_n)+\overline{\Phi}(c_{n,1-\alpha}+t_n)}{2\overline{\Phi}(c_{n,1-\alpha})}\to (e^{\sqrt{2}\tau}+e^{-\sqrt{2}\tau})/2=\cosh(\sqrt{2}\tau).$$ Together with~\eqref{eq:EtaGaussianCompare}, 
$\mc{G}_n(P_n^1)\to\theta$ now follows from (recall the first convergence in~\eqref{eqn:auxmill})
\begin{equation*}
P_n^1\del[1]{\enVert[0]{Z_d+t_n\bm{\iota}_d}_\infty\leq c_{n,1-\alpha}}
=
\del[2]{1-\frac{1}{d}(du_n)\frac{v_n}{u_n}}^d 
\to 
\exp\cbr[1]{\log(1-\alpha)\cosh(\sqrt{2}\tau)} = 1-\theta.
\end{equation*}

To establish Part 2, recall that~$p_n\to 0$. Thus, for~$n$ large enough~$ap_n^{-1/m}>2b_1$. Furthermore, for every~$n\in\N(c)$, we have
\begin{equation*}
\frac{\overline{\eta}_n}{q_n}
=
\del[1]{\sqrt{n\log(d \vee 2)}\overline{\eta}_n^{1-1/m}}^{\frac{m}{m-1}}
\geq
c^{\frac{m}{m-1}}
=
2\kappa,
\end{equation*}
and hence~$2p_n=2\kappa q_n\leq\overline{\eta}_n$. Thus, by Assumption~\ref{ass:lbd} and Lemma~\ref{lem:HDTVDist}, for every~$n \in \N(c)$ large enough,
\begin{align*}
&\mathsf{TV}\del[1]{
P_n^0 \circ [\mathsf{C}_n(X_{1,n}, \hdots, X_{n,n})]^{-1},
P_n^1\circ [\mathsf{C}_n(X_{1,n}, \hdots, X_{n,n})]^{-1} } \\
& = \mathsf{TV}\del[1]{Q_{0, d}^{\otimes n} \circ \mathsf{C}^{-1}_n, Q_{p_n, d}^{\otimes n} \circ \mathsf{C}^{-1}_n} \leq \del[2]{\frac{e}{4}}^{np_n} \to 0,
\end{align*}
the convergence following from~\eqref{eq:EtaRates0}. This proves~\eqref{eq:EtaBoundary2}.

To establish the statements in Part 3., just note that for every sequence of tests~$\psi_n:(\R^d)^n\to[0,1]$, by~\eqref{eq:EtaBoundary2} and, e.g., Lemma 2.3 in~\cite{strasser1985mathematical}, along~$n \in \N(c)$ we have
\begin{equation}\label{eq:EtaTVExpectations}
\envert[1]{E_{P_n^1}\psi_n(\mathsf{C}_n(X_{1,n}, \hdots, X_{n,n}))-E_{P_n^0}\psi_n(\mathsf{C}_n(X_{1,n}, \hdots, X_{n,n}))}
\to 0.
\end{equation}
\end{proof}

The proof strategy underlying Theorem~\ref{thm:EtaNecessary} is explained in some more detail subsequently.
\begin{remark}\label{rem:necproof}
As can be seen from~\eqref{eqn:defqpd} together with Assumption~\ref{ass:lbd}, the construction underlying Theorem \ref{thm:EtaNecessary} exploits rare, dense
shifts. Under $P_n^1$, independently across~$i$, with probability $p_n$ the observation~$X_{i,n}$ receives a shift $h_n\iota_d$, say, common to all coordinates (and absent under $P_n^0$), where~$h_n \asymp p_n^{-1/m}$. This keeps the centered $m$th moments of the coordinates uniformly bounded, while (cf.~Lemma~\ref{lem:HDDistConstruct}) producing a coordinate-wise mean of order~$p_nh_n\asymp p_n^{1-1/m}$. Choosing~$p_n\asymp(n\log(d))^{-m/(2(m-1))}$ therefore places the mean at order~$1/\sqrt{n\log(d)}$, the local scale at which a dense common shift produces nontrivial power for the Gaussian max-test. If~$\sqrt{n\log(d)}\,\overline{\eta}_n^{\,1-1/m}\not\to0$, then the contamination budget is large enough to remove all rare shifts (along a subsequence) with probability tending to one, making the contaminated null and alternative samples asymptotically indistinguishable. The rare shifts also induce a non-identity equicorrelation matrix under~$P_n^1$; however, its off-diagonal correlation~$\rho_n$ satisfies~$\log^2(d)\rho_n\to0$, so a Gaussian comparison argument shows that~$\mathcal{G}_n(P_n^1)$ is asymptotically equivalent to that in the independent Gaussian case.
\end{remark}

\begin{proof}[Proof of Theorem~\ref{thm:EtaNecessaryshort}:] 
Note that all assumptions required in Theorem~\ref{thm:EtaNecessary} are maintained. Define~$\mathsf{C}_n$ as in that theorem for~$n>M$ (and, e.g., equal to the identity map for~$n \leq M$). Consider a sequence of tests~$\psi_n:(\R^d)^n\to[0,1]$ satisfying~\eqref{eqn:asysizeshort}. Then, also~\eqref{eqn:asysize0} is satisfied. Theorem~\ref{thm:EtaNecessary} hence  delivers~\eqref{eq:EtaBoundary1} and~\eqref{eqn:asypow00}, which imply the first and second inequality in~\eqref{eqn:powshort}.
\end{proof}

\section{Proof of Results in Section~\ref{sec:BootstrapApprox}}

\subsection{The limiting distribution of~$\enVert[0]{\Sigma_{0,n}^{1/2}(P_n)Z_d}_\infty$}
Lemma~\ref{lem:extreme} provides the limiting distribution of~$\enVert[0]{\Sigma_{0,n}^{1/2}(P_n)Z_d}_\infty$ under the condition that~$d \to \infty$. 
\begin{lemma}\label{lem:extreme}
Let~$b_1 \in (0, \infty)$,~$b_2 \in (b_1, \infty)$, and~$\bm{r} = (r_l)_{l \in \N}$ be a sequence in~$[0, 1)$, such that~$\log(l)r_l\to 0$ as $l \to \infty$. Assume that~$\mc{P}(b_1,b_2,2,\bm{r})\neq \emptyset$ and~$d \to \infty$. Then,
\begin{equation}\label{eq:extreme1}
\lim_{n\to\infty} \sup_{P \in \mc{P}(b_1,b_2,2,\bm{r})} \sup_{x\in\R} \envert[2]{P\del[2]{a_d\sbr[2]{\enVert[1]{\Sigma_{0,n}^{1/2}(P)Z_d}_\infty-b_d}\leq x}
-\exp(-2\exp(-x))}
=
0,
\end{equation}
where 
\begin{equation*}
a_d=\sqrt{2\log(d)}\quad\text{and}\quad b_d=\sqrt{2\log(d)}-\frac{\log\log(d)+\log(4\pi)}{2\sqrt{2\log(d)}}. 
\end{equation*}
As a consequence, for every~$\alpha\in(0,1)$, every sequence~$s^*_d\to 0$, and every sequence~$P_n \in \mc{P}(b_1,b_2,2,\bm{r})$, it holds that
\begin{equation}\label{eq:CVs}
c_{n,1-(\alpha \pm s^*_d)}(\Sigma_{0,n}(P_n))
=
c_{n,1-\alpha}(\Sigma_{0,n}(P_n))+o\del[1]{1/\sqrt{\log(d)}}
=
c_{n,1-\alpha}+o\del[1]{1/\sqrt{\log(d)}}
;
\end{equation}
where~$c_{n,1-(\alpha \pm s^*_d)}(\Sigma_{0,n})$ is eventually well-defined.
\end{lemma}
\begin{proof}
Because~$r_l \in [0, 1)$ and~$r_l \to 0$ as~$l \to \infty$, there exists a~$\delta = \delta(\bm{r}) < 1$ such that~$r_l \in [0, \delta)$ for every~$l \in \N$. Fix a sequence~$P_n \in \mc{P}(b_1,b_2,2,\bm{r})$. Fix~$x\in\R$, and note that for~$u_d=u_d(x)=\frac{x}{a_d}+b_d$ (which is eventually positive) it holds that
\begin{equation*}
P_n\del[2]{a_d\sbr[2]{\enVert[1]{\Sigma_{0,n}^{1/2}(P_n)Z_d}_\infty-b_d}\leq x}
=
P_n\del[4]{\bigcap_{j=1}^d\cbr[2]{-u_d \leq(\Sigma_{0,n}^{1/2}(P_n)Z_d)_j\leq u_d}}.
\end{equation*}
Because~$\max_{1\leq i<j\leq d}|\Sigma_{0,n,i,j}|\leq \delta$, it follows from~(2.56) in~\cite{gine2016mathematical} that there exists a constant~$C$ depending only on~$\delta$ such that
\begin{align*}
&\envert[4]{P_n\del[4]{\bigcap_{j=1}^d\cbr[2]{-u_d\leq(\Sigma_{0,n}^{1/2}(P_n)Z_d)_j\leq u_d}}-P_1\del[4]{\bigcap_{j=1}^d\cbr[2]{-u_d\leq (Z_d)_j\leq u_d}}} \\
& \leq 
C\cdot d\sum_{j=1}^dr_j\exp\del[2]{-\frac{u_d^2}{1+r_j}} =: r(\delta, d, \bm{r}, x);
\end{align*}
we have exploited that the distribution of~$Z_d$ under~$P_n$ is the same as under~$P_1$. By Theorems~1.5.1 and~1.5.3 in~\cite{leadbetter} it holds that~$d\del[1]{1-\Phi(u_d)} \to e^{-x}$. Therefore, by Lemma 4.3.2 in the same reference~$r(\delta, d, \bm{r}, x) \to 0$ as~$n \to \infty$. 
Because~$P_1\del[1]{a_d\sbr[1]{\enVert[0]{Z_d}_\infty-b_d}\leq x} \to \exp(-2\exp(-x))$, e.g., by Equation~(55) in~\cite{kock2023consistency}, we obtain
\begin{equation}\label{eqn:prePol}
P_n\del[2]{a_d\sbr[2]{\enVert[1]{\Sigma_{0,n}^{1/2}(P_n)Z_d}_\infty-b_d}\leq x} \to \exp(-2\exp(-x)).
\end{equation}
P\'{o}lya's theorem, cf.,~e.g., Lemma 2.11 in~\cite{van2000asymptotic}, shows that the convergence in~\eqref{eqn:prePol} is uniform in~$x \in \R$. Because~$P_n$ was arbitrary~\eqref{eq:extreme1} follows (via the equivalence of~1 and~2 in Lemma~\ref{lem:unifseq}). 

To establish~\eqref{eq:CVs}, note that by definition of~$c_{n,1-\alpha\pm s_d^*}(\Sigma_{0,n}(P_n))$ it holds (eventually, i.e., once~$1-(\alpha \pm s^*_d) \in (0, 1)$) that
\begin{align*}
1-(\alpha \pm s^*_d)
&=
P_n\del[2]{\enVert[1]{\Sigma_{0,n}^{1/2}(P_n)Z_d}_\infty\leq c_{n,1-(\alpha \pm s^*_d)}(\Sigma_{0,n})}\\
&=
P_n\del[2]{a_d\sbr[2]{\enVert[1]{\Sigma_{0,n}^{1/2}(P_n)Z_d}_\infty-b_d}
\leq a_d\sbr[1]{c_{n,1-(\alpha \pm s^*_d)}(\Sigma_{0,n}(P_n))-b_d}}.
\end{align*}
Thus, it follows from~\eqref{eq:extreme1} and~$s_d^*\to 0$ that
\begin{equation*}
a_d\sbr[1]{c_{n,1-(\alpha \pm s^*_d)}(\Sigma_{0,n}(P_n))-b_d}=-\log(-\log(1-\alpha)/2)+o(1),
\end{equation*}
which implies
\begin{equation*}
c_{n,1-(\alpha \pm s^*_d)}(\Sigma_{0,n}(P_n))
=
\sqrt{2\log(d)}-\frac{\log\log(d)+\log(4\pi)}{2\sqrt{2\log(d)}}-\frac{\log\del[0]{-\log(1-\alpha)/2}+o(1)}{\sqrt{2\log(d)}}.
\end{equation*}
Applied to~$s_d^* \equiv 0$ this delivers the first equality in~\eqref{eq:CVs}. The second follows from~\eqref{eq:CVDiag}.

\end{proof}

\subsection{Bounding~$c_{n,1-\alpha}^B$}
The following controls~$c_{n,1-\alpha}^B$ irrespective of the structure of~$\Sigma_{0,n}$.
\begin{lemma}\label{lem:BSAlwaysvalid}
Let~$b_1\in(0,\infty)$,~$b_2\in(b_1,\infty)$,~$m > 2$,~$\lambda_1' \in (1,\infty)$, and~$\lambda_2'\in(0,\infty)$. Assume that~$d \geq 2$. There exists a (positive) constant~$C^*$, say, depending only on~$b_1, b_2, \lambda_1', \lambda_2'$, and~$m$, such that
\begin{equation}\label{eqn:sndef}
\mathfrak{s}_n:=C^*\log(d)\mathfrak{d}^{1/2}(n,d, m, \overline{\eta}_n)
\end{equation}
satisfies for every~$\alpha\in(0,1)$, and for every~$n \in \N$ such that~$\mathfrak{s}_n <\alpha \wedge (1-\alpha)$, $\mathfrak{d}(n, d, m, \overline{\eta}_n) <\frac{b_1^2}{2C}$ (for~$C$ as in Lemma~\ref{lem:varwinsestim}), and~$\varepsilon_n'$ as in~\eqref{eq:epsprime} obeying~\eqref{eq:epscond'}, the inequality
\begin{equation}\label{eq:critvalcontrol}
\inf_{P \in \mc{P}(b_1,b_2,m) }P\del[1]{c_{n,1-\alpha-\mathfrak{s}_n}(\Sigma_{0,n}(P))\leq c^B_{n,1-\alpha}\leq  c_{n,1-\alpha+\mathfrak{s}_n}(\Sigma_{0,n}(P))}
\geq 
1-24/n.
\end{equation}

\end{lemma}

\begin{proof}
Fix~$P\in\mc{P}(b_1,b_2,m)$, $\alpha \in (0, 1)$,~$n \in \N$ such that~$\mathfrak{s}_n < \alpha \wedge (1-\alpha)$, $\mathfrak{d}(n, d, m, \overline{\eta}_n) <\frac{b_1^2}{2C}$ (for~$C$ as in Lemma~\ref{lem:varwinsestim}), and~$\varepsilon_n'$ satisfying~\eqref{eq:epscond'}. Without loss of generality we can assume that~$n > 24$, so that the assumptions in Lemma~\ref{lem:varwinsestim} are satisfied. By the Gaussian-to-Gaussian comparison inequality as in Proposition~2.1 of~\cite{chernozhuokov2022improved} (cf.~Proposition~2 in~\cite{chernozhukov2023high}) for an absolute constant~$C_1$ (since~$\Sigma_{0,n,j,j}=1$ for all~$j=1,\hdots,d$), we have
\begin{align*}
A_n &:=\sup_{H\in\mc{H}}\envert[2]{P\del[1]{\Sigma_{0,n}^{1/2}(P)Z_d\in H}-P\del[1]{\tilde{\Sigma}_{0,n}^{1/2}Z_d\in H\mid\tilde{X}_1,\hdots,\tilde{X}_n}} \\
&\leq 	C_1\del[2]{\log^2(d)\max_{1\leq j,k\leq d}\envert[1]{\tilde{\Sigma}_{0,n,j,k}-\Sigma_{0,j,k}(P)}}^{1/2}.
\end{align*} 
Let~$C''$ be as in Lemma~\ref{lem:varwinsestim} and set~$C^* = C_1 \times \sqrt{C''}$, which only depends on~$b_1,b_2,\lambda_1',\lambda_2'$, and~$m$. By~\eqref{eqn:correst} in Lemma~\ref{lem:varwinsestim}, and with~$\mathfrak{s}_n$ as in~\eqref{eqn:sndef} (for the just-defined constant~$C^*$), we have~$$P\left(\min_{j = 1, \hdots, d} \tilde{\sigma}_{n,j} \geq \frac{b_1}{\sqrt{2}} \text{ and }A_n\leq \mathfrak{s}_n\right) 
\geq 1-24/n.$$  Hence, not showing the dependence of~$\Sigma_{0, n}$ on~$P$ for brevity, with probability at least~$1-24/n$ it holds that~$\min_{j = 1, \hdots, d} \tilde{\sigma}_{n,j} \geq \frac{b_1}{\sqrt{2}}$,
\begin{equation*}
P\del[2]{\enVert[1]{\tilde{\Sigma}_{0,n}^{1/2}Z_d}_\infty\leq c_{n,1-\alpha+\mathfrak{s}_n}(\Sigma_{0,n})\mid \tilde{X}_1,\hdots, \tilde{X}_n}
\geq
P\del[2]{\enVert[1]{\Sigma_{0,n}^{1/2}Z_d}_\infty\leq c_{n,1-\alpha+\mathfrak{s}_n}(\Sigma_{0,n})}-\mathfrak{s}_n,
\end{equation*}
and
\begin{equation*}
P\del[2]{\enVert[1]{\tilde{\Sigma}_{0,n}^{1/2}Z_d}_\infty\leq c_{n,1-\alpha-\mathfrak{s}_n}(\Sigma_{0,n})\mid \tilde{X}_1,\hdots, \tilde{X}_n}
\leq
P\del[2]{\enVert[1]{\Sigma_{0,n}^{1/2}Z_d}_\infty\leq c_{n,1-\alpha-\mathfrak{s}_n}(\Sigma_{0,n})}+\mathfrak{s}_n.
\end{equation*}
The lower and upper bound just established both equal~$1-\alpha$. Recalling that (by definition)
\begin{equation*}
P\del[2]{\enVert[1]{\tilde{\Sigma}_{0,n}^{1/2}Z}_\infty \leq c^B_{n,1-\alpha}\mid \tilde{X}_1,\hdots, \tilde{X}_n} = 1-\alpha,
\end{equation*}
it follows that
\begin{equation*}
P\del[1]{c_{n,1-\alpha-\mathfrak{s}_n}(\Sigma_{0,n})\leq c^B_{n,1-\alpha}\leq  c_{n,1-\alpha+\mathfrak{s}_n}(\Sigma_{0,n})}
\geq 
1-24/n.
\end{equation*}
\end{proof}

\subsection{Proof of Theorem~\ref{thm:BSWins}}
Let~$P_n\in\mc{P}(b_1,b_2,m)$ for every~$n\in\N$ and let~$\mathfrak{s}_n$ be as in~\eqref{eqn:sndef}. We start with some preliminary observations: We note that~$\mathfrak{s}_n \to 0$ (and hence~$\mathfrak{d}(n, d, m, \overline{\eta}_n) \to 0$) by (B.9) of Lemma B.4 in \cite{kock2025high} (as restated in~\eqref{eq:conv3} in Lemma~\ref{lem:covcompare}). We suppress the dependence on~$P_n$ of several quantities notationally. For~$n$ large enough, which we assume from now on without loss of generality, it holds that~$\mathfrak{s}_n <\alpha \wedge (1-\alpha)$, $\mathfrak{d}(n, d, m, \overline{\eta}_n) <\frac{b_1^2}{2C}$,~$\varepsilon_n\in (0, 1/2)$, and~$ \varepsilon_n'$ satisfies~\eqref{eq:epscond'} (recall their definitions from~\eqref{eq:epsfam} and~\eqref{eq:epsprime}, respectively, and note that both sequences converge to~$0$). It follows from Lemma~\ref{lem:BSAlwaysvalid} that
\begin{equation}
\begin{aligned}\label{eqn:initlinfs}
P_n\del[1]{\enVert[0]{S_{n,W}}_\infty > c^B_{n,1-\alpha}}
&\leq
P_n\del[1]{\enVert[0]{S_{n,W}}_\infty> c_{n,1-\alpha-\mathfrak{s}_n}(\Sigma_{0,n})}+\frac{24}{n}, \\
P_n\del[1]{\enVert[0]{S_{n,W}}_\infty > c^B_{n,1-\alpha}}
&\geq
P_n\del[1]{\enVert[0]{S_{n,W}}_\infty> c_{n,1-\alpha+\mathfrak{s}_n}(\Sigma_{0,n})}-\frac{24}{n}.
\end{aligned}
\end{equation}
If~\eqref{eq:extrasslocgauss} is satisfied along a given subsequence~$n'$, it hence follows from Theorem~\ref{thm:SumGAWins} that
\begin{equation}
\begin{aligned}\label{eq:limsup}
\limsup_{n'\to\infty}\bigg[&P_{n'}\del[1]{\enVert[0]{S_{n',W}}_\infty > c^B_{n',1-\alpha}}\\ &~~~~ -P_{n'}\del[1]{\enVert[0]{\Sigma_{0,n'}^{1/2}Z_{d'}+\sqrt{n'}D^{-1}_{n'}\mu_{n'}}_\infty>c_{n',1-\alpha-\mathfrak{s}_{n'}}(\Sigma_{0,n'})}\bigg]
\leq 
0,
\end{aligned}
\end{equation}
and
\begin{equation}
\begin{aligned}\label{eq:liminf}
\liminf_{n'\to\infty}\bigg[&P_{n'}\del[1]{\enVert[0]{S_{n',W}}_\infty > c^B_{n',1-\alpha}}\\ &~~~~-P_{n'}\del[1]{\enVert[0]{\Sigma_{0,n'}^{1/2}Z_{d'}+\sqrt{n'}D^{-1}_{n'}\mu_{n'}}_\infty>c_{n',1-\alpha+\mathfrak{s}_{n'}}(\Sigma_{0,n'})}\bigg]
\geq 
0.
\end{aligned}
\end{equation}
Because~$\mathfrak{s}_n\to 0$, it holds that  
\begin{equation}\label{eqn:GsiW}
P_n\del[1]{\enVert[0]{\Sigma_{0,n}^{1/2}Z_{d}}_\infty>c_{n,1-(\alpha\pm \mathfrak{s}_n)}(\Sigma_{0,n})}
=
\alpha \pm \mathfrak{s}_n
\to
\alpha.
\end{equation}
Now, we can establish \eqref{eq:BSsize}: Note that if~$P_n\in\mc{P}^0(b_1,b_2,m)$ for every~$n$, then~\eqref{eq:extrasslocgauss} clearly holds (along~$n' \equiv n$). Hence,~\eqref{eq:limsup},~\eqref{eq:liminf}, and~\eqref{eqn:GsiW} deliver~$P_n(\enVert[0]{S_{n,W}}_\infty > c^B_{n,1-\alpha}) \to \alpha$. This establishes \eqref{eq:BSsize} because~$P_n\in\mc{P}^0(b_1,b_2,m)$ was arbitrary (via the equivalence of~1 and~2 in Lemma~\ref{lem:unifseq}).

We next establish~\eqref{eq:BSPower}: To this end, fix a sequence~$\bm{r}$ as in the statement of the theorem and suppose, furthermore, that~$P_n\in \mc{P}(b_1,b_2,m,\bm{r})$. From~\eqref{eq:CVs} of Lemma~\ref{lem:extreme} we recall that
$$
c_{n,1-(\alpha\pm \mathfrak{s}_n)}(\Sigma_{0,n}(P_n))
=
c_{n,1-\alpha}(\Sigma_{0,n}(P_n))+o\del[1]{1/\sqrt{\log(d)}};
$$ 
so that with~$\bm{\iota}_d = (1, \hdots, 1)'\in \R^d$ we obtain from the Gaussian anti-concentration inequality (cf., e.g., Lemma E.1 in \cite{kock2025high}, as restated in Lemma~\ref{lem:anticonc}) that
\begin{equation}\label{eqn:antiL}
\begin{aligned}
&P_n\del[1]{\enVert[0]{\Sigma_{0,n}^{1/2}Z_{d}+\sqrt{n}D^{-1}_n\mu_n}_\infty\leq c_{n,1-(\alpha\pm \mathfrak{s}_n)}(\Sigma_{0,n})}\\
&=P_n\left( 
\Sigma_{0,n}^{1/2}Z_{d} \in [- c_{n,1-(\alpha\pm \mathfrak{s}_n)}(\Sigma_{0,n}) \bm{\iota}_d - \sqrt{n}D^{-1}_n\mu_n,  c_{n,1-(\alpha\pm \mathfrak{s}_n)}(\Sigma_{0,n}) \bm{\iota}_d - \sqrt{n}D^{-1}_n\mu_n]
\right)
\\
&=P_n\left( 
\Sigma_{0,n}^{1/2}Z_{d} \in [- c_{n,1-\alpha }(\Sigma_{0,n}) \bm{\iota}_d - \sqrt{n}D^{-1}_n\mu_n,  c_{n,1-\alpha }(\Sigma_{0,n}) \bm{\iota}_d - \sqrt{n}D^{-1}_n\mu_n]
\right) + o(1)
\\
&=P_n\del[1]{\enVert[0]{\Sigma_{0,n}^{1/2}Z_{d}+\sqrt{n}D^{-1}_n\mu_n}_\infty\leq c_{n,1-\alpha}(\Sigma_{0,n})}+o(1).	
\end{aligned}
\end{equation}
Let~$n'$ be a subsequence of~$n$. In case~\eqref{eq:extrasslocgauss} is satisfied,~\eqref{eq:limsup},~\eqref{eq:liminf}, and~\eqref{eqn:antiL} yield
\begin{equation}\label{eq:BSlocal}
\envert[2]{P_{n'}\del[1]{\enVert[0]{S_{n',W}}_\infty > c^B_{n',1-\alpha}}-P_{n'}\del[1]{\enVert[0]{\Sigma_{0,n'}^{1/2}Z_{d'}+\sqrt{n'}D^{-1}_{n'}\mu_{n'}}_\infty>c_{n',1-\alpha}(\Sigma_{0,n'})}}
\to 
0.
\end{equation}
In case~\eqref{eq:extrasslocgauss} is not satisfied, then~\eqref{eq:extrasslocgaussLarge} is satisfied along a subsequence~$n''$~of~$n'$. Due to the Khatri-{\v{S}}id{\'a}k inequality (cf.~Theorem \ref{thm:ksi})~$c_{n,1-\alpha^*}(\Sigma_{0,n})\leq c_{n,1-\alpha^*} = O(\sqrt{\log(d)})$ for every~$\alpha^* \in (0, 1)$. Hence, by Theorem~\ref{thm:SumGAWinsLarge}, we have
\begin{equation}
\begin{aligned}\label{eq:BSconv}
&P_{n''}\del[1]{\enVert[0]{S_{n'',W}}_\infty> c_{n'',1-\alpha/2}(\Sigma_{0,n''})}\to 1, \quad \text{ and } \\ &P_{n''}\del[1]{\enVert[0]{\Sigma_{0,n''}^{1/2}Z_{d''}+\sqrt{n''}D^{-1}_{n''}\mu_{n''}}_\infty>c_{n'',1-\alpha}(\Sigma_{0,n''})} \to 1.
\end{aligned}
\end{equation}
By~\eqref{eqn:initlinfs} eventually (recall that~$\mathfrak{s}_n \to 0$)
\begin{equation}\label{eqn:c2lat}
\begin{aligned}
P_{n''}\del[1]{\enVert[0]{S_{n'',W}}_\infty > c^B_{n'',1-\alpha}}
\geq
P_{n''}\del[1]{\enVert[0]{S_{n'',W}}_\infty> c_{n'',1-\alpha/2}(\Sigma_{0,n''})}-\frac{24}{n''}\to
1.
\end{aligned}
\end{equation}
Summarizing the two cases, every subsequence~$n'$ has a subsequence~$n''$, say, such that
\begin{equation*}
P_{n''}\del[1]{\enVert[0]{S_{n'',W}}_\infty > c^B_{n'',1-\alpha}} - 
P_{n''}\del[1]{\enVert[0]{\Sigma_{0,n''}^{1/2}Z_{d''}+\sqrt{n''}D^{-1}_{n''}\mu_{n''}}_\infty>c_{n'',1-\alpha}(\Sigma_{0,n''})} \to 0,
\end{equation*}
which establishes~\eqref{eq:BSPower}, because the sequence~$P_n\in \mc{P}(b_1,b_2,m,\bm{r})$ was arbitrary (and where we used the equivalence between Parts 1 and 3 in Lemma~\ref{lem:unifseq}).

\subsection{Proof of Theorem~\ref{thm:comparison}}

We start with~\eqref{eq:NoGain}. Fix a sequence~$\bm{r}$ as in the statement of  Theorem~\ref{thm:comparison}. Let~$P_n$ be a sequence in~$\mc{P}(b_1,b_2,m,\bm{r})$. We suppress the dependence of several quantities on~$P_n$ throughout this proof. From the second part, respectively, of Theorems~\ref{thm:Wins} and~\ref{thm:BSWins} it follows that
\begin{equation}\label{eqn:FS51}
\begin{aligned}
&\limsup_{n\to\infty}\big|P_n\del[1]{\enVert[0]{S_{n,W}}_\infty > c_{n,1-\alpha}}-P_n\del[1]{\enVert[0]{S_{n,W}}_\infty > c^B_{n,1-\alpha}}\big| \\
= &
\limsup_{n\to\infty} \big| P_n\del[1]{\enVert[0]{\Sigma_{0,n}^{1/2}Z_d+\sqrt{n}D_n^{-1}\mu_n}_\infty>c_{n,1-\alpha}} \\
& \hspace{3cm} -
P_n\del[1]{\enVert[0]{\Sigma_{0,n}^{1/2}Z_d+\sqrt{n}D_n^{-1}\mu_n}_\infty>c_{n,1-\alpha}(\Sigma_{0,n})}\big|.
\end{aligned}
\end{equation}
Recall that~$c_{n,1-\alpha}(\Sigma_{0,n}(P_n))
=
c_{n,1-\alpha}+o\del[1]{1/\sqrt{\log(d)}}$ from~\eqref{eq:CVs} of Lemma~\ref{lem:extreme}, and argue as in~\eqref{eqn:antiL} to obtain 
\begin{equation*}
P_n\del[1]{\enVert[0]{\Sigma_{0,n}^{1/2}Z_d+\sqrt{n}D^{-1}_n\mu_n}_\infty\leq c_{n,1-\alpha}(\Sigma_{0,n})}
=
P_n\del[1]{\enVert[0]{\Sigma_{0,n}^{1/2}Z_d+\sqrt{n}D^{-1}_n\mu_n}_\infty\leq c_{n,1-\alpha}}+o(1).
\end{equation*}
Together with~\eqref{eqn:FS51} and~$P_n\in\mc{P}(b_1,b_2,m,\bm{r})$ being arbitrary, this implies~\eqref{eq:NoGain} (via the equivalence of~1 and~2 in Lemma~\ref{lem:unifseq}).

To prove~\eqref{eq:NoGainMeanCond} let~$P_n$ be a sequence in~$\mc{P}(b_1,b_2,m)$ satisfying ~$\sqrt{n}\|\mu_n\|_\infty/\sqrt{\log(d)}\to \infty$. By Theorem~\ref{thm:SumGAWinsLarge} together with Remark~\ref{rem:wecond} and~\eqref{eq:CVDiag},~$P_n\del[1]{\enVert[0]{S_{n,W}}_\infty>c_{n,1-\alpha}}\to 1$. Arguing as in~\eqref{eqn:c2lat}, we obtain~$P_n\del[1]{\enVert[0]{S_{n,W}}_\infty > c^B_{n,1-\alpha}} \to 1$. 

To prove~\eqref{eq:Gain}, let~$P_n$ be as indicated in Part 3 of Theorem~\ref{thm:comparison}. Clearly,~$\|\mu_n\|_\infty=\log^{1/4}(d)/\sqrt{n}$ for~$n$ large enough, so that~\eqref{eq:extrasslocgauss} is satisfied along~$n$ by (B.9) of Lemma B.4 in~\cite{kock2025high} (as restated in~\eqref{eq:conv3} in Lemma~\ref{lem:covcompare}). Therefore, Theorem~\ref{thm:SumGAWins} implies that
\begin{equation}\label{eqn:ridB}
\limsup_{n\to\infty}P_n\del[1]{\enVert[0]{S_{n,W}}_\infty > c_{n,1-\alpha}}
=
\limsup_{n\to\infty}P_n\del[1]{\enVert[0]{\Sigma_{0,n}^{1/2}Z_d+\sqrt{n}D_n^{-1}\mu_n}_\infty>c_{n,1-\alpha}}.
\end{equation}
Arguing as in~\eqref{eqn:initlinfs} (using Lemma~\ref{lem:BSAlwaysvalid} and~$\mathfrak{s}_n\to 0$)
\begin{equation*}
P_n\del[1]{\enVert[0]{S_{n,W}}_\infty > c^B_{n,1-\alpha}}
\geq
P_n\del[1]{\enVert[0]{S_{n,W}}_\infty > c_{n,1-\alpha+\mathfrak{s}_n}(\Sigma_{0,n})}-\frac{24}{n},
\end{equation*}
such that another application of Theorem~\ref{thm:SumGAWins} implies that
\begin{equation}\label{eqn:ridB2}
\liminf_{n\to\infty}P_n\del[1]{\enVert[0]{S_{n,W}}_\infty > c^B_{n,1-\alpha}}	
\geq
\liminf_{n\to\infty}P_n\del[1]{\enVert[0]{\Sigma_{0,n}^{1/2}Z_d+\sqrt{n}D_n^{-1}\mu_n}_\infty>c_{n,1-\alpha/2}(\Sigma_{0,n})}.
\end{equation}
Next, note that~$\sqrt{n}D_n^{-1}\mu_n=\log^{1/4}(d)/b_1\cdot \bm{\iota}_d$, that~$\Sigma_{0,n}$ is the~$d\times d$ matrix with all entries equal to one, and that~$\Sigma_{0,n}^{1/2}$ is the~$d\times d$ matrix with all entries equal to~$d^{-1/2}$. Therefore, for every~$c\in\R$
\begin{equation}\label{eq:reduction}
P_n\del[2]{\enVert[0]{\Sigma_{0,n}^{1/2}Z_d+\sqrt{n}D_n^{-1}\mu_n}_\infty>c}
=
P_n\del[1]{|Z_1+\log^{1/4}(d)/b_1|>c}.
\end{equation}
Since~$c_{n,1-\alpha}\geq \sqrt{\log(d)}$ for~$d$ sufficiently large,~\eqref{eqn:ridB} together with~\eqref{eq:reduction} implies
\begin{equation*}
\limsup_{n\to\infty}P_n\del[1]{\enVert[0]{S_{n,W}}_\infty > c_{n,1-\alpha}}
=
\limsup_{n\to\infty}P_n\del[1]{|Z_1+\log^{1/4}(d)/b_1|>c_{n,1-\alpha}}
=
0.	
\end{equation*}
Since~$c_{n,1-\alpha/2}(\Sigma_{0,n})=\Phi^{-1}(1-\alpha/4)$ [with~$\Phi$ denoting the standard normal cdf on~$\R$],~\eqref{eqn:ridB2} together with~\eqref{eq:reduction} implies
\begin{equation*}
\liminf_{n\to\infty}P_n\del[1]{\enVert[0]{S_{n,W}}_\infty > c_{n,1-\alpha}^B}
\geq
\liminf_{n\to\infty}P_n\del[1]{|Z_1+\log^{1/4}(d)/b_1|>\Phi^{-1}(1-\alpha/4)}
=
1.	
\end{equation*}

\section{Further numerical checks}\label{sec:RobCheck}
In this section we investigate the sensitivity of the empirical rejection frequencies of
\begin{equation*}
\varphi_{n,W}
=
\mathds{1}\cbr[1]{\enVert[0]{S_{n,W}}_\infty>c_{n,1-\alpha}}\qquad \text{and}\qquad \varphi_{n,W,B}=\mathds{1}\cbr[1]{\enVert[0]{S_{n,W}}_\infty>c_{n,1-\alpha}^B}
\end{equation*}
to
\begin{itemize}
\item misspecifying the upper bound on the fraction (or, equivalently, number) of contaminated observations. 
\item the choice of~$\lambda_2'$.
\end{itemize} 
Finally, we also repeat the size experiments in Tables~\ref{tab:1} and~\ref{tab:2} of the main text but for the more heavy-tailed setting of~$\nu=3.01$.

To understand the mechanism underlying the sensitivity checks, recall first that,
\begin{equation*}
\eps_n'= \lambda_1'\cdot\overline{\eta}_n+\lambda_2'\cdot\frac{\log(d^2n)}{n},\qquad \lambda_1'\in(1,\infty)\text{ and }\lambda_2'\in (0,\infty)
\end{equation*}
and that the number of winsorized observations used to form the denominators~$\tilde{\sigma}_{n,j}$ of the statistics~$|S_{n,W,j}|=\frac{|S^{\dagger}_{n, W,j}|}{\tilde{\sigma}_{n,j}}$, which~$\varphi_{n,W}$ and~$\varphi_{n,W,B}$ are both based on, is increasing in~$\eps_n'$ for~$j=1,\hdots,d$. Thus,~$\tilde{\sigma}_{n,j}$ is decreasing in~$\eps_n'$ and hence~$\enVert[0]{S_{n,W}}_\infty$ is increasing in~$\eps_n'$. As a result, since the critical values~$c_{n,1-\alpha}$ do not depend on~$\eps_n'$, it is mechanically the case that larger~$\eps_n'$ lead to larger rejection frequencies for~$\varphi_{n,W}$. In particular, since~$\eps_n'$ is an increasing function of~$\lambda_2'$, this implies that the rejection frequencies of~$\varphi_{n,W}$ are increasing in~$\lambda_2'$. This monotonicity of the rejection frequencies need not hold for~$\varphi_{n,W,B}$ as also the bootstrap critical values~$c_{n,1-\alpha}^B$ change when~$\eps_n'$ changes.

\emph{Unless stated otherwise, the setup in the tables below is as in Section~\ref{sec:Sims} of the main text.} 

We first investigate the sensitivity of the rejection frequencies of~$\varphi_{n,W}$ and~$\varphi_{n,W,B}$ to specifying the upper bound~$\overline{\eta}_n$ on the fraction of contaminated observations. Equivalently, we specify an upper bound~$q_{\text{imp}}=\lfloor\overline{\eta}_nn\rfloor\in\cbr[0]{0,\hdots,5}$ on the number of contaminated observations. Throughout, the true number of contaminated observations is~$q_{\text{true}}=3$. As in the main text, the contaminated~$X_{i,n}$ are chosen randomly and changed to~$(500,\hdots,500)'$. Table~\ref{tab:q_mis} shows that for all values of~$\rho\in\cbr[0]{0,0.9,1}$, the rejection frequencies of~$\varphi_{n,W}$ and~$\varphi_{n,W,B}$ are increasing in~$q_{\text{imp}}$. It is encouraging that in the settings studied, overspecifying $q_{\text{imp}}$ does not lead to overly excessive size inflation. We also implemented the setting of Table~\ref{tab:q_mis} but lowering~$n$ to~$200$ (unreported): $\varphi_{n,W}$ and~$\varphi_{n,W,B}$ became more sensitive to the choice of~$q_{\text{imp}}$ for~$\rho\in\cbr[0]{0,0.9}$. The rejection frequencies remained increasing in $q_{\text{imp}}$.

\begin{table}[ht]
\centering
\begin{tabular}{@{}c cc cc cc@{}}
\toprule
& \multicolumn{2}{c}{$\rho=0$}
& \multicolumn{2}{c}{$\rho=0.9$}
& \multicolumn{2}{c}{$\rho=1$} \\
\cmidrule(lr){2-3}\cmidrule(lr){4-5}\cmidrule(lr){6-7}
$q_{\mathrm{imp}}$
& $\varphi_{n,W}$ & $\varphi_{n,W,B}$
& $\varphi_{n,W}$ & $\varphi_{n,W,B}$
& $\varphi_{n,W}$ & $\varphi_{n,W,B}$ \\
\midrule
$0$ & 0.000 & 0.000 & 0.000 & 0.000 & 0.000 & 0.000 \\
$1$ & 0.000 & 0.000 & 0.000 & 0.000 & 0.000 & 0.000 \\
$2$ & 0.000 & 0.000 & 0.000 & 0.000 & 0.000 & 0.000 \\
$3$ & 0.027 & 0.034 & 0.019 & 0.038 & 0.000 & 0.050 \\
$4$ & 0.049 & 0.051 & 0.036 & 0.070 & 0.000 & 0.052 \\
$5$ & 0.073 & 0.076 & 0.051 & 0.086 & 0.000 & 0.065 \\
\bottomrule
\end{tabular}
\caption{Empirical rejection frequencies when the true number of contaminated observations is~$q_{\mathrm{true}}=3$ and the tests are implemented using~$q_{\mathrm{imp}}\in\{0,\ldots,5\}$. Here,~$n=500$,~$d=1{,}000$, and~$\rho$ is as indicated.}
\label{tab:q_mis}
\end{table}
Next, Tables~\ref{tab:lambda_q0} and~\ref{tab:lambda_q3} investigate the sensitivity of the rejection frequencies of $\varphi_{n,W}$ and~$\varphi_{n,W,B}$ to the choice of~$\lambda_2'$. As noted above, the rejection frequency of the former is, \emph{by construction}, increasing in~$\lambda_2'$. Furthermore, recall from Section~\ref{sec:Sims} of the main text that all values of~$\lambda_2'\in(0,0.04]$ lead to numerically identical implementations of $\varphi_{n,W}$ and~$\varphi_{n,W,B}$. Both of these observations are confirmed by Tables~\ref{tab:lambda_q0} and~\ref{tab:lambda_q3}. In addition, the tables also reveal that the rejection frequency of~$\varphi_{n,W,B}$ is increasing in~$\lambda_2'$, thus motivating our recommendation of choosing~$\lambda_2'$ ``small'' in order to avoid excessive null rejection frequencies. The exact small value of~$\lambda_2'$ matters less. It is seen that choosing~$\lambda_2'\geq 0.2$ is not advisable. However, across all combinations of~$q_{\text{true}}=q_{\text{imp}}\in\cbr[0]{0,3}$ and~$\rho\in\cbr[0]{0,0.9,1}$ it is the case that~$\lambda_2'\leq 0.1$ leads to null rejection frequencies below~$0.1$.
\begin{table}[ht]
\centering
\begin{tabular}{@{}c cc cc cc@{}}
\toprule
& \multicolumn{2}{c}{$\rho=0$}
& \multicolumn{2}{c}{$\rho=0.9$}
& \multicolumn{2}{c}{$\rho=1$} \\
\cmidrule(lr){2-3}\cmidrule(lr){4-5}\cmidrule(lr){6-7}
$\lambda_2'$
& $\varphi_{n,W}$ & $\varphi_{n,W,B}$
& $\varphi_{n,W}$ & $\varphi_{n,W,B}$
& $\varphi_{n,W}$ & $\varphi_{n,W,B}$ \\
\midrule
$0.010$ & 0.021 & 0.022 & 0.021 & 0.037 & 0.000 & 0.051 \\
$0.025$ & 0.021 & 0.022 & 0.021 & 0.037 & 0.000 & 0.051 \\
$0.050$ & 0.031 & 0.030 & 0.032 & 0.050 & 0.000 & 0.059 \\
$0.075$ & 0.031 & 0.030 & 0.032 & 0.050 & 0.000 & 0.059 \\
$0.100$ & 0.042 & 0.047 & 0.037 & 0.059 & 0.000 & 0.065 \\
$0.200$ & 0.064 & 0.070 & 0.052 & 0.081 & 0.000 & 0.072 \\
$0.500$ & 0.159 & 0.170 & 0.106 & 0.155 & 0.000 & 0.087 \\
$1.000$ & 0.404 & 0.409 & 0.239 & 0.328 & 0.000 & 0.113 \\
\bottomrule
\end{tabular}
\caption{Empirical rejection frequencies for different choices of~$\lambda_2'$ without contamination ($q_{\text{true}}=q_{\text{imp}}=0$). Here,~$n=500$,~$d=1{,}000$, and~$\rho$ is as indicated.}
\label{tab:lambda_q0}
\end{table}

\begin{table}[ht]
\centering
\begin{tabular}{@{}c cc cc cc@{}}
\toprule
& \multicolumn{2}{c}{$\rho=0$}
& \multicolumn{2}{c}{$\rho=0.9$}
& \multicolumn{2}{c}{$\rho=1$} \\
\cmidrule(lr){2-3}\cmidrule(lr){4-5}\cmidrule(lr){6-7}
$\lambda_2'$
& $\varphi_{n,W}$ & $\varphi_{n,W,B}$
& $\varphi_{n,W}$ & $\varphi_{n,W,B}$
& $\varphi_{n,W}$ & $\varphi_{n,W,B}$ \\
\midrule
$0.010$ & 0.027 & 0.036 & 0.019 & 0.037 & 0.000 & 0.046 \\
$0.025$ & 0.027 & 0.036 & 0.019 & 0.037 & 0.000 & 0.046 \\
$0.050$ & 0.053 & 0.057 & 0.037 & 0.073 & 0.000 & 0.056 \\
$0.075$ & 0.053 & 0.057 & 0.037 & 0.073 & 0.000 & 0.056 \\
$0.100$ & 0.082 & 0.085 & 0.056 & 0.093 & 0.000 & 0.066 \\
$0.200$ & 0.134 & 0.135 & 0.088 & 0.140 & 0.001 & 0.070 \\
$0.500$ & 0.308 & 0.310 & 0.188 & 0.264 & 0.001 & 0.099 \\
$1.000$ & 0.625 & 0.630 & 0.403 & 0.518 & 0.003 & 0.126 \\
\bottomrule
\end{tabular}
\caption{Empirical rejection frequencies for different choices of~$\lambda_2'$ under contamination ($q_{\text{true}}=q_{\text{imp}}=3$). Here,~$n=500$,~$d=1{,}000$, and~$\rho$ is as indicated.}
\label{tab:lambda_q3}
\end{table}

Finally, Tables~\ref{tab:nu301_q0} and~\ref{tab:nu301_q3} repeat the size experiments in Tables~\ref{tab:1} and~\ref{tab:2} of the main text with~$\nu=3.01$ degrees of freedom instead of~$\nu=4.01$. The qualitative conclusions are unchanged. Absent contamination, the tests generally control size, the non-bootstrap tests become more conservative as~$\rho$ increases, and the bootstrap tests adapt to the dependence structure. Relative to the case~$\nu=4.01$, the rejection frequencies of the winsorized tests are typically somewhat lower, especially for~$\rho\in\cbr[0]{0,0.9}$, while those of~$\varphi_{n,B}$ and~$\varphi_{n,LL}$ are broadly similar. Under contamination,~$\varphi_n$ and~$\varphi_{n,B}$ again have rejection frequencies close to zero, the size distortions of~$\varphi_{n,W}$ and~$\varphi_{n,W,B}$ are most pronounced for small~$n$ and diminish as~$n$ increases, and~$\varphi_{n,LL}$ remains oversized, particularly for~$\rho=1$. Compared with~$\nu=4.01$, the winsorized tests are again generally more conservative, whereas the conclusions for the remaining tests are very similar.

\begin{landscape}
\begin{table}[p]
\centering
\setlength{\tabcolsep}{4pt}
\renewcommand{\arraystretch}{0.88}
\begin{minipage}[t]{0.40\linewidth}
\centering
{\footnotesize
\makebox[\linewidth][c]{$\rho=0$}\par\vspace{2pt}

\resizebox{\linewidth}{!}{%
\begin{tabular}{cccccccc}
\toprule
$n$&$d$& $\varphi_n$ & $\varphi_{n,B}$ & $\varphi_{n,W}$ & $\varphi_{n,W,B}$ & $\varphi_{n,LL}$ & $\varphi_{n,GSS}$ \\
\midrule
\multirow{4}{*}{$100$}
&$200$ & 0.043 & 0.054 & 0.023 & 0.027 & 0.053 & 0.003 \\
&$500$ & 0.040 & 0.053 & 0.020 & 0.023 & 0.053 & 0.002 \\
&$1000$ & 0.054 & 0.072 & 0.021 & 0.032 & 0.058 & 0.002 \\
&$2000$ & 0.042 & 0.065 & 0.020 & 0.026 & 0.059 & 0.000 \\
\midrule
\multirow{4}{*}{$200$}
&$200$ & 0.043 & 0.046 & 0.015 & 0.018 & 0.054 & 0.013 \\
&$500$ & 0.034 & 0.043 & 0.012 & 0.014 & 0.047 & 0.006 \\
&$1000$ & 0.048 & 0.064 & 0.014 & 0.015 & 0.041 & 0.010 \\
&$2000$ & 0.042 & 0.051 & 0.018 & 0.019 & 0.053 & 0.005 \\
\midrule
\multirow{4}{*}{$500$}
&$200$ & 0.039 & 0.047 & 0.017 & 0.020 & 0.047 & 0.033 \\
&$500$ & 0.029 & 0.038 & 0.016 & 0.016 & 0.060 & 0.020 \\
&$1000$ & 0.041 & 0.047 & 0.013 & 0.014 & 0.035 & 0.017 \\
&$2000$ & 0.048 & 0.054 & 0.012 & 0.018 & 0.042 & 0.025 \\
\midrule
\multirow{4}{*}{$1000$}
&$200$ & 0.044 & 0.047 & 0.016 & 0.018 & 0.040 & 0.032 \\
&$500$ & 0.041 & 0.044 & 0.017 & 0.019 & 0.042 & 0.032 \\
&$1000$ & 0.041 & 0.044 & 0.018 & 0.017 & 0.039 & 0.031 \\
&$2000$ & 0.046 & 0.052 & 0.016 & 0.015 & 0.045 & 0.031 \\
\bottomrule
\end{tabular}}
}
\end{minipage}
\hspace{0.04\linewidth}
\begin{minipage}[t]{0.40\linewidth}
\centering
{\footnotesize
\makebox[\linewidth][c]{$\rho=0.9$}\par\vspace{2pt}

\resizebox{\linewidth}{!}{%
\begin{tabular}{cccccccc}
\toprule
$n$&$d$& $\varphi_n$ & $\varphi_{n,B}$ & $\varphi_{n,W}$ & $\varphi_{n,W,B}$ & $\varphi_{n,LL}$ & $\varphi_{n,GSS}$ \\
\midrule
\multirow{4}{*}{$100$}
&$200$ & 0.030 & 0.056 & 0.011 & 0.034 & 0.061 & 0.003 \\
&$500$ & 0.029 & 0.061 & 0.016 & 0.028 & 0.041 & 0.001 \\
&$1000$ & 0.036 & 0.068 & 0.016 & 0.027 & 0.051 & 0.001 \\
&$2000$ & 0.031 & 0.066 & 0.013 & 0.029 & 0.058 & 0.000 \\
\midrule
\multirow{4}{*}{$200$}
&$200$ & 0.031 & 0.064 & 0.019 & 0.030 & 0.061 & 0.013 \\
&$500$ & 0.032 & 0.057 & 0.014 & 0.027 & 0.045 & 0.012 \\
&$1000$ & 0.030 & 0.061 & 0.011 & 0.021 & 0.047 & 0.003 \\
&$2000$ & 0.034 & 0.062 & 0.016 & 0.021 & 0.052 & 0.003 \\
\midrule
\multirow{4}{*}{$500$}
&$200$ & 0.027 & 0.049 & 0.014 & 0.030 & 0.037 & 0.016 \\
&$500$ & 0.022 & 0.047 & 0.013 & 0.020 & 0.040 & 0.014 \\
&$1000$ & 0.023 & 0.047 & 0.004 & 0.012 & 0.046 & 0.011 \\
&$2000$ & 0.020 & 0.036 & 0.010 & 0.012 & 0.044 & 0.009 \\
\midrule
\multirow{4}{*}{$1000$}
&$200$ & 0.036 & 0.059 & 0.015 & 0.034 & 0.053 & 0.027 \\
&$500$ & 0.022 & 0.040 & 0.009 & 0.022 & 0.038 & 0.016 \\
&$1000$ & 0.024 & 0.042 & 0.012 & 0.017 & 0.043 & 0.017 \\
&$2000$ & 0.039 & 0.059 & 0.011 & 0.021 & 0.038 & 0.024 \\
\bottomrule
\end{tabular}}
}
\end{minipage}
\par\vspace{0.5cm}
\begin{minipage}[t]{0.40\linewidth}
\centering
{\footnotesize
\makebox[\linewidth][c]{$\rho=1$}\par\vspace{2pt}

\resizebox{\linewidth}{!}{%
\begin{tabular}{cccccccc}
\toprule
$n$&$d$& $\varphi_n$ & $\varphi_{n,B}$ & $\varphi_{n,W}$ & $\varphi_{n,W,B}$ & $\varphi_{n,LL}$ & $\varphi_{n,GSS}$ \\
\midrule
\multirow{4}{*}{$100$}
&$200$ & 0.001 & 0.062 & 0.000 & 0.030 & 0.056 & 0.000 \\
&$500$ & 0.000 & 0.058 & 0.000 & 0.028 & 0.055 & 0.000 \\
&$1000$ & 0.000 & 0.056 & 0.000 & 0.031 & 0.056 & 0.000 \\
&$2000$ & 0.000 & 0.040 & 0.000 & 0.016 & 0.050 & 0.000 \\
\midrule
\multirow{4}{*}{$200$}
&$200$ & 0.001 & 0.046 & 0.000 & 0.030 & 0.049 & 0.000 \\
&$500$ & 0.000 & 0.044 & 0.000 & 0.023 & 0.042 & 0.000 \\
&$1000$ & 0.000 & 0.053 & 0.000 & 0.038 & 0.063 & 0.000 \\
&$2000$ & 0.000 & 0.047 & 0.000 & 0.025 & 0.053 & 0.000 \\
\midrule
\multirow{4}{*}{$500$}
&$200$ & 0.000 & 0.064 & 0.000 & 0.044 & 0.071 & 0.000 \\
&$500$ & 0.001 & 0.046 & 0.000 & 0.026 & 0.056 & 0.000 \\
&$1000$ & 0.000 & 0.046 & 0.000 & 0.030 & 0.057 & 0.000 \\
&$2000$ & 0.000 & 0.052 & 0.000 & 0.035 & 0.055 & 0.000 \\
\midrule
\multirow{4}{*}{$1000$}
&$200$ & 0.000 & 0.051 & 0.000 & 0.043 & 0.058 & 0.000 \\
&$500$ & 0.000 & 0.044 & 0.000 & 0.032 & 0.049 & 0.000 \\
&$1000$ & 0.000 & 0.058 & 0.000 & 0.038 & 0.055 & 0.000 \\
&$2000$ & 0.000 & 0.061 & 0.000 & 0.043 & 0.057 & 0.000 \\
\bottomrule
\end{tabular}}
}
\end{minipage}
\caption{\footnotesize Empirical rejection frequencies under the null absent contamination. Each tabular has~$\Sigma_{n,i,j}=\frac{\nu}{\nu-2}\rho^{|i-j|}$ for~$\nu=3.01$ and~$\rho$ as indicated.}
\label{tab:nu301_q0}
\end{table}
\end{landscape}

\begin{landscape}
\begin{table}[p]
\centering
\setlength{\tabcolsep}{4pt}
\renewcommand{\arraystretch}{0.88}
\begin{minipage}[t]{0.40\linewidth}
\centering
{\footnotesize
\makebox[\linewidth][c]{$\rho=0$}\par\vspace{2pt}

\resizebox{\linewidth}{!}{%
\begin{tabular}{cccccccc}
\toprule
$n$&$d$& $\varphi_n$ & $\varphi_{n,B}$ & $\varphi_{n,W}$ & $\varphi_{n,W,B}$ & $\varphi_{n,LL}$ & $\varphi_{n,GSS}$ \\
\midrule
\multirow{4}{*}{$100$}
&$200$ & 0.000 & 0.006 & 0.045 & 0.067 & 0.081 & 0.000 \\
&$500$ & 0.000 & 0.011 & 0.035 & 0.061 & 0.080 & 0.000 \\
&$1000$ & 0.000 & 0.007 & 0.040 & 0.067 & 0.088 & 0.000 \\
&$2000$ & 0.000 & 0.008 & 0.040 & 0.076 & 0.096 & 0.000 \\
\midrule
\multirow{4}{*}{$200$}
&$200$ & 0.000 & 0.009 & 0.035 & 0.041 & 0.099 & 0.000 \\
&$500$ & 0.000 & 0.005 & 0.035 & 0.050 & 0.107 & 0.000 \\
&$1000$ & 0.000 & 0.006 & 0.017 & 0.024 & 0.093 & 0.000 \\
&$2000$ & 0.000 & 0.009 & 0.027 & 0.032 & 0.100 & 0.000 \\
\midrule
\multirow{4}{*}{$500$}
&$200$ & 0.000 & 0.003 & 0.030 & 0.038 & 0.130 & 0.000 \\
&$500$ & 0.000 & 0.007 & 0.021 & 0.026 & 0.124 & 0.000 \\
&$1000$ & 0.000 & 0.002 & 0.021 & 0.032 & 0.128 & 0.000 \\
&$2000$ & 0.000 & 0.009 & 0.021 & 0.030 & 0.116 & 0.000 \\
\midrule
\multirow{4}{*}{$1000$}
&$200$ & 0.000 & 0.008 & 0.022 & 0.025 & 0.108 & 0.000 \\
&$500$ & 0.000 & 0.003 & 0.018 & 0.020 & 0.122 & 0.000 \\
&$1000$ & 0.000 & 0.004 & 0.023 & 0.027 & 0.118 & 0.000 \\
&$2000$ & 0.000 & 0.004 & 0.023 & 0.026 & 0.118 & 0.000 \\
\bottomrule
\end{tabular}}
}
\end{minipage}
\hspace{0.04\linewidth}
\begin{minipage}[t]{0.40\linewidth}
\centering
{\footnotesize
\makebox[\linewidth][c]{$\rho=0.9$}\par\vspace{2pt}

\resizebox{\linewidth}{!}{%
\begin{tabular}{cccccccc}
\toprule
$n$&$d$& $\varphi_n$ & $\varphi_{n,B}$ & $\varphi_{n,W}$ & $\varphi_{n,W,B}$ & $\varphi_{n,LL}$ & $\varphi_{n,GSS}$ \\
\midrule
\multirow{4}{*}{$100$}
&$200$ & 0.000 & 0.009 & 0.031 & 0.081 & 0.080 & 0.000 \\
&$500$ & 0.000 & 0.012 & 0.036 & 0.089 & 0.090 & 0.000 \\
&$1000$ & 0.000 & 0.010 & 0.032 & 0.091 & 0.094 & 0.000 \\
&$2000$ & 0.000 & 0.006 & 0.034 & 0.081 & 0.088 & 0.000 \\
\midrule
\multirow{4}{*}{$200$}
&$200$ & 0.000 & 0.011 & 0.017 & 0.047 & 0.107 & 0.000 \\
&$500$ & 0.000 & 0.003 & 0.027 & 0.053 & 0.089 & 0.000 \\
&$1000$ & 0.000 & 0.006 & 0.017 & 0.046 & 0.092 & 0.000 \\
&$2000$ & 0.000 & 0.009 & 0.015 & 0.053 & 0.088 & 0.000 \\
\midrule
\multirow{4}{*}{$500$}
&$200$ & 0.000 & 0.001 & 0.017 & 0.038 & 0.145 & 0.000 \\
&$500$ & 0.000 & 0.006 & 0.011 & 0.038 & 0.128 & 0.000 \\
&$1000$ & 0.000 & 0.001 & 0.009 & 0.023 & 0.123 & 0.000 \\
&$2000$ & 0.000 & 0.005 & 0.015 & 0.035 & 0.123 & 0.000 \\
\midrule
\multirow{4}{*}{$1000$}
&$200$ & 0.000 & 0.004 & 0.014 & 0.030 & 0.139 & 0.000 \\
&$500$ & 0.000 & 0.008 & 0.016 & 0.034 & 0.134 & 0.000 \\
&$1000$ & 0.000 & 0.009 & 0.016 & 0.029 & 0.122 & 0.000 \\
&$2000$ & 0.000 & 0.005 & 0.009 & 0.024 & 0.131 & 0.000 \\
\bottomrule
\end{tabular}}
}
\end{minipage}
\par\vspace{0.5cm}
\begin{minipage}[t]{0.40\linewidth}
\centering
{\footnotesize
\makebox[\linewidth][c]{$\rho=1$}\par\vspace{2pt}

\resizebox{\linewidth}{!}{%
\begin{tabular}{cccccccc}
\toprule
$n$&$d$& $\varphi_n$ & $\varphi_{n,B}$ & $\varphi_{n,W}$ & $\varphi_{n,W,B}$ & $\varphi_{n,LL}$ & $\varphi_{n,GSS}$ \\
\midrule
\multirow{4}{*}{$100$}
&$200$ & 0.000 & 0.008 & 0.001 & 0.048 & 0.146 & 0.000 \\
&$500$ & 0.000 & 0.009 & 0.001 & 0.038 & 0.143 & 0.000 \\
&$1000$ & 0.000 & 0.008 & 0.000 & 0.033 & 0.145 & 0.000 \\
&$2000$ & 0.000 & 0.007 & 0.000 & 0.044 & 0.148 & 0.000 \\
\midrule
\multirow{4}{*}{$200$}
&$200$ & 0.000 & 0.004 & 0.000 & 0.034 & 0.153 & 0.000 \\
&$500$ & 0.000 & 0.006 & 0.000 & 0.030 & 0.128 & 0.000 \\
&$1000$ & 0.000 & 0.009 & 0.000 & 0.032 & 0.157 & 0.000 \\
&$2000$ & 0.000 & 0.004 & 0.000 & 0.035 & 0.167 & 0.000 \\
\midrule
\multirow{4}{*}{$500$}
&$200$ & 0.000 & 0.013 & 0.000 & 0.039 & 0.177 & 0.000 \\
&$500$ & 0.000 & 0.009 & 0.000 & 0.043 & 0.153 & 0.000 \\
&$1000$ & 0.000 & 0.010 & 0.000 & 0.027 & 0.139 & 0.000 \\
&$2000$ & 0.000 & 0.007 & 0.000 & 0.032 & 0.157 & 0.000 \\
\midrule
\multirow{4}{*}{$1000$}
&$200$ & 0.000 & 0.017 & 0.000 & 0.039 & 0.160 & 0.000 \\
&$500$ & 0.000 & 0.016 & 0.000 & 0.030 & 0.154 & 0.000 \\
&$1000$ & 0.000 & 0.015 & 0.000 & 0.051 & 0.179 & 0.000 \\
&$2000$ & 0.000 & 0.016 & 0.000 & 0.037 & 0.181 & 0.000 \\
\bottomrule
\end{tabular}}
}
\end{minipage}

\caption{\footnotesize Empirical rejection frequencies under the null in the presence of contamination. The adversary replaces~$q=3$ randomly chosen~$X_{i,n}$ by~$\tilde{X}_{i,n}=(500,\hdots,500)'$ and leaves the remaining observations unchanged. Each tabular has~$\Sigma_{n,i,j}=\frac{\nu}{\nu-2}\rho^{|i-j|}$ for~$\nu=3.01$ and~$\rho$ as indicated.}
\label{tab:nu301_q3}
\end{table}
\end{landscape}

\end{appendix}

\end{document}